\documentclass[english]{article}
\usepackage[T1]{fontenc}
\usepackage[latin9]{inputenc}
\usepackage{geometry}
\usepackage{float}
\usepackage{booktabs}
\usepackage{mathtools}
\usepackage{amsmath}
\usepackage{amsthm}
\usepackage{amssymb}
\usepackage{graphicx}

\makeatletter

\providecommand{\tabularnewline}{\\}
\floatstyle{ruled}
\newfloat{algorithm}{tbp}{loa}
\providecommand{\algorithmname}{Algorithm}
\floatname{algorithm}{\protect\algorithmname}

\theoremstyle{remark}
\newtheorem*{rem*}{\protect\remarkname}
\theoremstyle{plain}
\newtheorem*{assumption*}{\protect\assumptionname}
\theoremstyle{plain}
\newtheorem{thm}{\protect\theoremname}
\theoremstyle{plain}
\newtheorem{cor}[thm]{\protect\corollaryname}

\usepackage{algorithm,algpseudocode}

\@ifundefined{showcaptionsetup}{}{%
 \PassOptionsToPackage{caption=false}{subfig}}
\usepackage{subfig}
\makeatother

\usepackage{babel}
\providecommand{\assumptionname}{Assumption}
\providecommand{\corollaryname}{Corollary}
\providecommand{\remarkname}{Remark}
\providecommand{\theoremname}{Theorem}

\begin{document}
\title{High-order Moment Closure Models with Random Batch Method for Efficient Computation of Multiscale Turbulent Systems}
\author{Di Qi\textsuperscript{a} and Jian-Guo Liu\textsuperscript{b}}
\date{\textsuperscript{a }Department of Mathematics, Purdue University,
150 North University Street, West Lafayette, IN 47907, USA\\ \textsuperscript{b }Department of Mathematics and Department of Physics, Duke University, Durham, NC 27708, USA\\}
\maketitle

\begin{abstract}
We propose a high-order stochastic-statistical moment closure model
for efficient ensemble prediction of leading-order statistical moments
and probability density functions in multiscale complex turbulent
systems. The statistical moment equations are closed by a precise
calibration of the high-order feedbacks using ensemble solutions of
the consistent stochastic equations, suitable for modeling complex
phenomena including non-Gaussian statistics and extreme events. To
address challenges associated with closely coupled spatio-temporal
scales in turbulent states and expensive large ensemble simulation
for high-dimensional systems, we introduce efficient computational
strategies using the random batch method (RBM). This approach significantly
reduces the required ensemble size while accurately capturing essential
high-order structures. Only a small batch of small-scale fluctuation
modes is used for each time update of the samples, and exact convergence
to the full model statistics is ensured through frequent resampling
of the batches during time evolution. Furthermore, we develop a reduced-order
model to handle systems with really high dimension by linking the
large number of small-scale fluctuation modes to ensemble samples
of dominant leading modes. The effectiveness of the proposed models
is validated by numerical experiments on the one-layer and two-layer
Lorenz '96 systems, which exhibit representative chaotic features
and various statistical regimes. The full and reduced-order RBM models
demonstrate uniformly high skill in capturing the time evolution of
crucial leading-order statistics, non-Gaussian probability
distributions, while achieving significantly lower computational cost
compared to direct Monte-Carlo approaches. The models provide effective
tools for a wide range of real-world applications in prediction, uncertainty
quantification, and data assimilation.
\end{abstract}

\section{Introduction\label{sec:Introduction}}

Turbulent dynamical systems encountered in science and engineering
\cite{majda2006nonlinear,pedlosky2013geophysical,nicholson1983introduction,kalnay2003atmospheric}
exhibit distinguished characteristics including a high-dimensional
state space with multiple spatio-temporal scales and strong internal
instabilities \cite{Frisch1995,majda2016introduction}. Forecasting
the intricate behaviors of such complex systems poses a significant
challenge in uncertainty quantification and data assimilation problems
\cite{reich2015probabilistic,majda2018strategies,majda2019linear,calvello2022ensemble}.
One key aspect involves accurately quantifying the multiscale interaction
between the large-scale mean state and the many interacting small-scale
fluctuations induced by internal instability. The interplay between
multiscale coupling and instability gives rise to a diverse array
of complex phenomena, such as bursting extreme structures and skewed
non-Gaussian distributions \cite{vanden2001non,cousins2015unsteady,tong2021extreme,qi2022anomalous}.
Small randomness in initial conditions and external stochastic effects
is also rapidly amplified and redistributed along the spectrum as
the model evolves in time. Capturing these unique features with efficient
algorithms remains a central issue in many practical problems \cite{chen2020predicting,gao2022master,gao2023data,gao2023transition}.
For example in data assimilation, accurate prediction of the mean
and covariance statistics is essential, while the evolution of these
low-order moments is closely linked to higher-order feedbacks due
to the nonlinear coupling. This requires the accurate and efficient
quantification for the extreme outliers and the related non-Gaussian
statistics.

In developing efficient algorithms for accurate statistical prediction
of complex turbulent systems featuring multiscale interactions and
strong nonlinearity, a probabilistic approach is usually needed to
quantify the uncertainty utilizing a probability density function
(PDF) of the model states. Ensemble forecasting, through a Monte Carlo
(MC) type method, is commonly used to estimate the PDF evolution by
independently sampling an ensemble of trajectories from an initial
distribution \cite{leutbecher2008ensemble}. Ensemble-based approaches
have been extensively applied to address uncertainties arising from
various sources in real world problems such as weather and climate
forecast \cite{houtekamer1996system,kalnay2003atmospheric,palmer2019ecmwf}.
However, achieving accurate numerical prediction is hampered by the
prohibitively high computational cost attributed to the strongly coupled
nonlinear interactions among different scales in a high dimensional
space, known as `curse of dimensionality' \cite{friedman1997bias,donoho2000high}.
With insufficient samples, ensemble approximation often suffer from
the collapse of samples, leading to an inadequate coverage of the
entire probability measure in the high-dimensional phase space. Various
strategies have been devised to increase the effective ensemble size
and sufficiently characterize the probability distribution \cite{bishop2001adaptive,giggins2019stochastically}.
Parameterization is often used for efficient approximation of the
effects from unresolved small-scale processes using information from
the resolved large-scale states \cite{sapsis2013statistically,majda2018strategies,maulik2019subgrid,luo2022stability,qi2023data}.
In the case of high-dimensional turbulent systems, model errors are
significantly amplified by internal instability and careful calibrations
of the feedbacks from the large number of unresolved processes are
often required. Consequently, this often leads to an inherent difficulty
for generalization to high dimensional and strongly turbulent systems
in realistic applications.

In this paper, we aim to develop a systematic modeling and computational strategy
for statistical prediction and data assimilation of turbulent systems
under a unified framework. The proposed formulation \eqref{eq:abs_formu}
is applicable to complex spatially extended nonlinear dynamical systems
widely studied in many fields \cite{lesieur1987turbulence,kalnay2003atmospheric,Diamond2010,majda2016introduction,majda2019statistical}.
In particular, we construct a coupled stochastic-statistical model
\eqref{coupled_model} that is suitable for both efficient prediction
of leading-order statistics (in the statistical equations) and explicit
quantification of higher-order non-Gaussian features (using the stochastic
equations) at the same time. The  stochastic and statistical equations
are seamlessly coupled through the nonlinear interaction terms and using 
an essential relaxation term for consistency. The statistical equations
admit a hierarchical structure requiring a closure form for the higher-order
moment feedback induced by the quadratic nonlinear term. Instead of
the parameterization methods which usually require exhausting procedure
of model calibration, we close the moment equations by explicitly
modeling the higher-order feedbacks through the ensemble solution
of the stochastic equations. This high-order moment closure through
the stochastic-statistical model provides a statistically consistent
formulation that is able to correctly represent the crucial features
involving non-Gaussian and nonlinear phenomena.

The computational demand for solving the fully coupled system \eqref{coupled_model}
remains substantial when the state variable resides in a high-dimensional
space. For example, for a system with $K$-dimensional state space,
one time step update of the covariance dynamics requires computational
cost of $O\left(K^{4}\right)$ by evolving $K\times K$ matrix-valued
differential equations with $K^{2}$ quadratic coupling terms in each
entry. The ensemble simulation for the stochastic dynamics further
requires the computational cost of $O\left(MK^{3}\right)$ by updating
the $K$ stochastic coefficients using $M$ samples (usually with
exponential dependent on $K$). To overcome this issue, we use the
effective random batch method (RBM) \cite{jin2020random,jin2021convergence}
for efficient computation of the large number of high-order coupling
terms reaching a much reduced computational cost. The proposed method
generalizes the idea developed for simple mean-fluctuation systems
in \cite{qi2023random} and provides efficient computational strategy
for the more universal formulation \eqref{coupled_model} valid for
statistical prediction in a large group of problems. Using the RBM
approximation, the high-dimensional modes are divided into small batches
randomly drawn in each time updating interval. Then the nonlinear
interaction is solely computed among a small subset of modes within
one batch so that the computational cost is maintained in a low level.
During the time evolution, the batches are frequently resampled at
the start of each time updating step. The accuracy is preserved by
the ergodicity of the fast mixing modes, ensuring an equivalent statistics
in high-order feedback. The resulting RBM model \eqref{rbm_model}
offers significant computational reduction with a cost of $O\left(p^{2}K^{2}\right)$
for the covariance equations and $O\left(M_{1}p^{2}K\right)$ for
the ensemble forecast, where a considerable smaller ensemble size
$M_{1}\ll M$ is sufficient only requiring sampling the batch subspace consisting of $p=O\left(1\right)$ modes in one batch. The convergence of the mean and variance
using RBM approximation is proved through detailed error estimates
dependent on the discrete time step size. Furthermore, for really
high dimensional systems with an extended spectrum, a reduced-order
model \eqref{eq:rbm_stoch2} is proposed to further reduce the computational
cost to $O\left(M_{1}\left(K_{1}^{3}+p^{3}\right)\right)$, making it independent
of the full dimension $K$ by focusing on the first $K_{1}$ leading
modes and using RBM to approximation the rest $K-K_{1}$ small-scale
fluctuation modes. 

The performance of the stochastic-statistical model with RBM approximations
is examined under the one-layer and two-layer Lorenz '96 (L-96) systems
\cite{lorenz1996predictability,wilks2005effects}. The Lorenz '96
systems have been widely used as prototype models for the atmosphere
involving rich chaotic phenomena with direct link to realistic systems
\cite{smith2001disentangling,orrell2003model,giggins2019stochastically,chen2020predicting}.
In particular, the L-96 systems maintain a slowly decay variance spectrum
involving a large number of unstable modes (as illustrated in Figure
\ref{fig:Illustration-l96} and \ref{fig:Growth-rates-2layer}), setting
up a very challenging testing case for effective ensemble prediction.
The direct MC approach requires a very large ensemble size of order
$O\left(10^{5}\right)$ to resolve the highly non-Gaussian statistics.
The strong internal instability leads to additional problems for easy
divergence away from the final equilibrium state (shown in Figure
\ref{fig:instability}). Using the efficient RBM model, it is found
that a small sample size of $O\left(100\right)$ is sufficient to
fully recover the evolutions of statistical mean and variance as well
as the non-Gaussian PDFs in the 40-dimensional one-layer L-96 system
(see Figure \ref{fig:1D-marginal-PDFs} and \ref{fig:2D-joint-PDFs}).
In the genuinely high dimensional two-layer L-96 system with full
dimension $K=264$. The RBM model becomes especially effective to
capture the highly non-Gaussian statistics using at most $M_{1}=500$
samples (see Figure \ref{fig:1D-marginal-PDFs-2layer} and \ref{fig:2D-joint-PDFs-2layer}).
The computational cost is further reduced in the reduced-order model
focusing on the leading $K_{1}=8$ large-scale modes enabling an even
smaller ensemble for all the small-scale modes.

In the remainder part of this paper, we introduce the general formulation
of the stochastic-statistical model for multiscale turbulent systems
using high-order moment closure in Section \ref{sec:stat_stoc_model}.
The efficient algorithms for solving the coupled high-dimensional
equations using the RBM approximation for ensemble prediction are
developed in Section \ref{sec:Computational-methods-rbm} together
with the theoretical convergence analysis of the scheme. The performance
of the methods is evaluated on the concrete examples of the Lorenz
'96 systems in Section \ref{sec:Numerical-performance}. The paper
is closed with a summary and discussions on future directions in Section \ref{sec:Summary}.

\section{A statistically consistent modeling framework for multiscale turbulent systems\label{sec:stat_stoc_model}}

Turbulent dynamical systems are characterized by multiscale nonlinear
interactions, which redistribute energy across a broad spectrum of
stable and unstable modes, ultimately leading to a complicated statistical
equilibrium. The general formulation of complex turbulent systems
can be introduced in the following canonical equations about the state
variable $\mathbf{u}\in\mathbb{R}^{d}$ in a high-dimensional phase
space (with $d\gg1$)
\begin{equation}
\frac{\mathrm{d}\mathbf{u}}{\mathrm{d}t}=\Lambda\mathbf{u}+B\left(\mathbf{u},\mathbf{u}\right)+\mathbf{F}\left(t\right)+\boldsymbol{\sigma}\left(t\right)\dot{\mathbf{W}}\left(t\right).\label{eq:abs_formu}
\end{equation}
The model state starts from an initial distribution $\mathbf{u}\left(0\right)\sim\mu_{0}\left(\mathbf{u}\right)$
representing initial uncertainty. On the right hand side of the equation
\eqref{eq:abs_formu}, the first component, $\Lambda=L-D$, represents
linear dispersion and dissipation effects, where the dispersion $L^{*}=-L$
is an energy-conserving skew-symmetric operator; and the dissipation
$D<0$ is a negative-definite operator. The model \eqref{eq:abs_formu}
emphasizes the important role of nonlinear interactions in a bilinear
quadratic form, $B\left(\mathbf{u},\mathbf{u}\right)$. This typical
structure of nonlinear interactions is inherited from a discretization
of the continuous full system (for example, a spectral projection
of the nonlinear advection in fluid models). The nonlinear interaction
ensures the energy conservation invariance, such that $\mathbf{u}\cdot B\left(\mathbf{u},\mathbf{u}\right)\equiv0$
with the inner-product defined according to the conserved quantity.
In addition, the system is subject to external forcing effects that
are decomposed into a deterministic component, $\mathbf{F}\left(t\right)$,
and a stochastic component represented by a Gaussian random process,
$\boldsymbol{\sigma}\left(t\right)\dot{\mathbf{W}}\left(t\right)$,
used to model the unresolved processes. 

The evolution of the model state $\mathbf{u}$ depends on sensitivity
to the randomness in initial conditions and stochastic forcing effects.
These uncertainties will be amplified in time by the inherent internal
instability due to the nonlinear coupling term in \eqref{eq:abs_formu}.
The associated Fokker-Planck equation (FPE) \cite{varadhan2007stochastic}
\begin{equation}
\frac{\partial p_{t}}{\partial t}=\mathcal{L}_{\mathrm{FP}}p_{t}\coloneqq-\mathrm{div}_{\mathrm{u}}\left[\Lambda\mathbf{u}+B\left(\mathbf{u,u}\right)+\mathbf{F}\right]p_{t}+\frac{1}{2}\mathrm{div}_{\mathbf{u}}\nabla\left(\boldsymbol{\sigma}\boldsymbol{\sigma}^{T}p_{t}\right),\label{eq:FPE}
\end{equation}
describes the time evolution of the probability density function (PDF)
$p_{t}\left(\mathbf{u}\right)=e^{t\mathcal{L}_{\mathrm{FP}}\left(\mathbf{u}\right)}\mu_{0}$
starting from an initial distribution $p_{t=0}\left(\mathbf{u}\right)=\mu_{0}\left(\mathbf{u}\right)$.
However, it remains a challenging task in directly solving the FPE
\eqref{eq:FPE} as a high dimensional PDE system. As an alternative
approach, ensemble forecast by tracking the Monte-Carlo (MC) solutions
estimates the model statistics through empirical averages among a
group of independent samples drawn from the initial distribution $\mathbf{u}^{\left(i\right)}\left(0\right)\sim\mu_{0}\left(\mathbf{u}\right),i=1,\cdots,M$
at the starting time $t=0$. The PDF solution $p_{t}\left(\mathbf{u}\right)$
and the associated statistical expectation of any function $\varphi\left(\mathbf{u}\right)$
at each time instant $t>0$ are then approximated by the empirical
ensemble representation of the $M$ samples, that is,
\begin{equation}
p_{t}\left(\mathbf{u}\right)\simeq p_{t}^{\mathrm{MC}}\left(\mathbf{u}\right)\coloneqq\frac{1}{M}\sum_{i=1}^{M}\delta\left(\mathbf{u}-\mathbf{u}^{\left(i\right)}\left(t\right)\right),\quad\left\langle \varphi\left(\mathbf{u}\right)\right\rangle _{p_{t}}\simeq\left\langle \varphi\right\rangle _{p_{t}^{\mathrm{MC}}}\coloneqq\frac{1}{M}\sum_{i=1}^{M}\varphi\left(\mathbf{u}^{\left(i\right)}\left(t\right)\right),\label{eq:pdf_MC}
\end{equation}
where $\delta$ is the Dirac delta function and $\left\langle \cdot\right\rangle _{p}$
is the expectation about the probability measure $p$. Still, the
\emph{curse of dimensionality} \cite{donoho2000high,daum2003curse}
arises in systems of even moderate dimension $d$ since the model
errors grow significantly as the system dimension increases, while
only a small ensemble size $M$ is allowed in practical numerical
methods due to the limited computational resources. Clearly, efficient
strategies and algorithms are still needed to effectively reduce the
computational cost and maintain high accuracy in sampling the high
dimensional systems using a small number of samples.
\begin{rem*}
Many complex turbulent systems from nature and engineering can be
categorized into the general mathematical framework in \eqref{eq:abs_formu}.
The high-dimensional state $\mathbf{u}$ can be viewed as a finite
dimensional truncation of the corresponding continuous field with
sufficient numerical resolution. One major group of examples comes
from the fluid flows including the Navier-Stokes equation and turbulence
at high Reynolds number \cite{Frisch1995,lesieur1987turbulence} and
applications to the geophysical models in coupled atmosphere and ocean
systems involving rotation, stratification and topography \cite{pedlosky2013geophysical,majda2006nonlinear,kalnay2003atmospheric}
and controlled fusion in magnetically confined plasma systems \cite{Diamond2010,dewar2007zonal}.
In particular, we will consider the prototype Lorenz '96 systems in
\eqref{eq:l96_homo} and \eqref{eq:l96_2layer} \cite{lorenz1996predictability,arnold2013stochastic,wilks2005effects}
that admit all representative dynamical structures in \eqref{eq:abs_formu}
as the main test model in this paper.
\end{rem*}

\subsection{Statistical and stochastic formulations for multiscale systems}

One major difficulty in complex turbulent systems is the fully coupled
nonlinear interactions across scales. The multiscale interactions
involve a large-scale mean state, which can destabilize the smaller
scales, while the excited fluctuation energy contained in numerous
small-scale modes can inversely impact the development of the coherent
structure at largest scale. Thus, disregarding contributions from
small-scale modes through a simple high wavenumber truncation is not
a viable approach. To address this central issue of coupled interactions
with mixed scales, we start with a mean-fluctuation decomposition
for the model state $\mathbf{u}$, so that interactions between different
scales can be identified in detail. To achieve this, we view $\mathbf{u}$
as a random field (denoted by $\omega$ due to randomness in initial
state and stochastic forcing) and separate it into the composition
of a statistical mean state $\bar{\mathbf{u}}$ and a wide spectrum
of stochastic fluctuations $\mathbf{u}^{\prime}$ in a finite $K$-dimensional
projected representation under a fixed-in-time, orthonormal basis
$\left\{ \mathbf{v}_{k}\right\} _{k=1}^{K}$ (usually with $K=d$
for the full model and $K<d$ for the reduced-order model)
\begin{equation}
\mathbf{u}\left(t;\omega\right)=\bar{\mathbf{u}}\left(t\right)+\mathbf{u}^{\prime}\left(t;\omega\right)=\bar{\mathbf{u}}\left(t\right)+\sum_{k=1}^{K}Z_{k}\left(t;\omega\right)\mathbf{v}_{k}.\label{eq:decomp}
\end{equation}
Above, $\bar{\mathbf{u}}=\left\langle \mathbf{u}\right\rangle _{p_{t}}$
represents the \emph{statistical mean field} usually capturing the
dominant large-scale structure; and $Z_{k}\left(t;\omega\right)$
is the \emph{stochastic coefficient} measuring the uncertainty in
multiscale fluctuation processes $\mathbf{u}^{\prime}$ projected
on the eigenmode $\mathbf{v}_{k}$. The state decomposition \eqref{eq:decomp}
provides a convenient way to identify different components in the
multiscale interactions, thus can be used to derive new effective
multiscale models.

One way to avoid the high computational cost in directly solving the
FPE \eqref{eq:FPE} as well as the large ensemble MC simulation of
the full SDE \eqref{eq:abs_formu} for the probability distribution
$p_{t}\left(\mathbf{u}\right)$ is to seek a hierarchical statistical
description of its moments $\ensuremath{\left\langle \varphi\left(\mathbf{u}\right)\right\rangle _{p_{t}}}$
as the expectation with respect to the time-dependent probability
measure $p_{t}$. In most situations, the primary interest lies in
tracking the time evolution of the leading moments quantifying the
most essential statistical characteristics. We can first derive the
dynamics for the mean state $\bar{\mathbf{u}}=\left\langle \mathbf{u}\right\rangle _{p_{t}}$
and the covariance among fluctuation modes $R_{kl}=\left\langle Z_{k}Z_{l}^{*}\right\rangle _{p_{t}}$
governed by the following set of \emph{deterministic statistical equations}
\addtocounter{equation}{0}\begin{subequations}\label{stat} 
\begin{align}
\frac{\mathrm{d}\bar{\mathbf{u}}}{\mathrm{d}t} & =\Lambda\bar{\mathbf{u}}+B\left(\bar{\mathbf{u}},\bar{\mathbf{u}}\right)+\sum_{k,l=1}^{K}R_{kl}B\left(\mathbf{v}_{k},\mathbf{v}_{l}\right)+\mathbf{F},\label{eq:stat_mean}\\
\frac{\mathrm{d}R_{kl}}{\mathrm{d}t} & =\sum_{m=1}^{K}L_{v,km}\left(\bar{\mathbf{u}}\right)R_{ml}+R_{km}L_{v,kl}^{*}\left(\bar{\mathbf{u}}\right)+Q_{\sigma,kl}\label{eq:stat_cov}\\
 & +\sum_{m,n=1}^{K}\left\langle Z_{m}Z_{n}Z_{l}\right\rangle _{p_{t}}B\left(\mathbf{v}_{m},\mathbf{v}_{n}\right)\cdot\mathbf{v}_{k}+\left\langle Z_{m}Z_{n}Z_{k}\right\rangle _{p_{t}}B\left(\mathbf{v}_{m},\mathbf{v}_{n}\right)\cdot\mathbf{v}_{l},\nonumber 
\end{align}
\end{subequations}for all the wavenumbers $1\leq k,l\leq K$. In
\eqref{eq:stat_cov}, the operator $L_{v,kl}=\left[\Lambda\mathbf{v}_{l}+B\left(\bar{\mathbf{u}},\mathbf{v}_{l}\right)+B\left(\mathbf{v}_{l},\bar{\mathbf{u}}\right)\right]\cdot\mathbf{v}_{k}$
characterizes quasilinear coupling between the statistical mean and
stochastic modes; the positive-definite operator $Q_{\sigma,kl}=\sum_{m}\left(\mathbf{v}_{k}\cdot\boldsymbol{\sigma}_{m}\right)\left(\boldsymbol{\sigma}_{m}\cdot\mathbf{v}_{l}\right)$
expresses energy injection from the stochastic forcing. Notably, the
nonlinear flux term involving all the third-order moments $\left\langle Z_{m}Z_{n}Z_{k}\right\rangle _{p_{t}}$
enters the equation for the second-order covariance $R_{kl}$ describing
nonlinear energy exchanges among fluctuation modes, ending up with
a still unclosed set of equations. 

Accordingly, the stochastic coefficients $\left\{ Z_{k}\right\} _{k=1}^{K}$
in the decomposition \eqref{eq:decomp} satisfy the associated \emph{stochastic
fluctuation equations}
\begin{equation}
\begin{aligned}\frac{\mathrm{d}Z_{k}}{\mathrm{d}t}= & \sum_{m=1}^{K}L_{v,km}\left(\bar{\mathbf{u}}\right)Z_{m}+\sigma\left(t\right)\dot{\mathbf{W}}\left(t\right)\cdot\mathbf{v}_{k}\\
 & +\sum_{m,n=1}^{K}\left(Z_{m}Z_{n}-R_{mn}\right)B\left(\mathbf{v}_{m},\mathbf{v}_{n}\right)\cdot\mathbf{v}_{k}.
\end{aligned}
\label{eq:dyn_fluc}
\end{equation}
The above equation is achieved by simply projecting the original equation
\eqref{eq:abs_formu} on each fluctuation mode $\mathbf{v}_{k}$ and
removing the mean dynamics \eqref{eq:stat_mean}. The second-order
moments equation \eqref{eq:stat_cov} is then derived by applying
It\^o's formula to $\varphi\left(Z_{k}\right)=\left|Z_{k}\right|^{2}$
using the stochastic equations \eqref{eq:dyn_fluc}. Therefore, the
above statistical and stochastic formulations are consistent for the
evolution of uncertainty in multiscale fluctuations. The detailed
derivation of the equations \eqref{stat} and \eqref{eq:dyn_fluc}
from first principle can be found in \cite{majda2016introduction,majda2018strategies}.

Both the dynamical moment representation \eqref{stat} and its stochastic
counterpart \eqref{eq:dyn_fluc} have their respective advantages,
they also both suffer several difficulties when applied to resolve
the key statistical quantities. The statistical moment equations \eqref{stat}
are easier to compute with its deterministic dynamics, while such
hierarchical equations lead to a non-closed system of infinite-dimensional
ODEs as each lower-order moment equation is coupled to the next higher-order
moment. On the other hand, the stochastic equations \eqref{eq:dyn_fluc}
provide a closed formulation to include all the higher-order information.
However, direct simulation of the SDE requires a MC approach of large
sample size exponentially dependent on the state dimension. In addition,
the computational cost of both statistical and stochastic models remains
prohibitive for the coupled high-dimensional systems characterized
by an extended wide spectrum of fluctuation modes $K\gg1$. The situation
becomes especially challenging when an ensemble approach is required
for accurate state estimation and data assimilation of extreme events
represented by the extended PDF tails.

\subsection{A stochastic-statistical closure model with explicit higher-order feedbacks}

We introduce a seamless high-order closure model that integrates the
statistical equations \eqref{stat} with the stochastic counterpart
\eqref{eq:dyn_fluc} to effectively close the original non-closed
equations. The resulting \emph{coupled stochastic-statistical equations
for the multiscale interacting model} becomes\addtocounter{equation}{0}\begin{subequations}\label{coupled_model} 
\begin{align}
\frac{\mathrm{d}\bar{\mathbf{u}}}{\mathrm{d}t} & =\Lambda\bar{\mathbf{u}}+B\left(\bar{\mathbf{u}},\bar{\mathbf{u}}\right)+\sum_{k,l=1}^{K}R_{kl}B\left(\mathbf{v}_{k},\mathbf{v}_{l}\right)+\mathbf{F},\label{eq:dyn_mean}\\
\frac{\mathrm{d}Z_{k}}{\mathrm{d}t} & =\sum_{m=1}^{K}L_{v,km}\left(\bar{\mathbf{u}}\right)Z_{m}+\sum_{m,n=1}^{K}\gamma_{mnk}\left(Z_{m}Z_{n}-R_{mn}\right)+\sigma_{k}\dot{W}_{k},\label{eq:dyn_stoc}\\
\frac{\mathrm{d}R_{kl}}{\mathrm{d}t} & =\sum_{m=1}^{K}L_{v,km}\left(\bar{\mathbf{u}}\right)R_{ml}+R_{km}L_{v,ml}^{*}\left(\bar{\mathbf{u}}\right)+Q_{\sigma,kl}\label{eq:dyn_cov}\\
 & +\sum_{m,n=1}^{K}\left[\gamma_{mnk}\left\langle Z_{m}Z_{n}Z_{l}\right\rangle _{p_{t}}+\gamma_{mnl}\left\langle Z_{m}Z_{n}Z_{k}\right\rangle _{p_{t}}\right]+\epsilon^{-1}\left(\left\langle Z_{k}Z_{l}\right\rangle _{p_{t}}-R_{kl}\right),\nonumber 
\end{align}
\end{subequations}with the coupling coefficients $L_{v,km}=\left[\Lambda\mathbf{v}_{m}+B\left(\bar{\mathbf{u}},\mathbf{v}_{m}\right)+B\left(\mathbf{v}_{m},\bar{\mathbf{u}}\right)\right]\cdot\mathbf{v}_{k}$
and $\gamma_{mnk}=B\left(\mathbf{v}_{m},\mathbf{v}_{n}\right)\cdot\mathbf{v}_{k}$
derived from the original equations. Above, the mean equation \eqref{eq:dyn_mean}
for the leading-order statistics $\bar{\mathbf{u}}$ is kept the same
aiming to capturing the dominant large-scale mean structures. The
covariance equation \eqref{eq:dyn_cov} for $R$ is closed by computing
the expectations of the cubic terms under the probability measure
$p_{t}$ discovered by the stochastic solution $\mathbf{Z}$ from
\eqref{eq:dyn_stoc}, so that the higher-order moment feedbacks are
explicitly modeled. In addition, a relaxation term is added with a
small control parameter $\epsilon$ to guarantee statistical consistency.
Importantly, both the statistical and stochastic equations become
indispensable for the modeling of the fully coupled multiscale system:
i) the stochastic equations for the fluctuation modes $Z_{k}$ are
introduced to provide exact closure for the covariance $R$; and ii)
the covariance dynamics for $R$ serves as an auxiliary equation to
facilitate the explicit interactions between the mean $\bar{\mathbf{u}}$
and stochastic modes $Z_{k}$ and can deal with the inherent instability
in the turbulent systems.

In developing effective strategy to compute the high-order expectation
in the covariance equation \eqref{eq:dyn_cov} according to the PDF
$p_{t}$ of the stochastic solutions in \eqref{eq:dyn_stoc}, the
stochastic equations for the random fluctuation modes $Z_{k}$ are
solved through an ensemble approach using $\mathbf{Z}^{\left(i\right)}=\left\{ Z_{k}^{\left(i\right)}\right\} $
with sample index $i=1,\cdots,M$. The higher-order feedbacks are
then approximated through the empirical average of the ensemble as
in \eqref{eq:pdf_MC}, $\left\langle \varphi\left(\mathbf{Z}\right)\right\rangle _{p_{t}}\approx\frac{1}{M}\sum_{i=1}^{M}\varphi\left(\mathbf{Z}^{\left(i\right)}\right)$.
Compared with the direct MC approach of the original system \eqref{eq:abs_formu},
the new coupled model \eqref{coupled_model} adopting the explicitly
coupled stochastic-statistical dynamics presenting an attractive equivalent
formulation that enjoys the flexibility of developing efficient reduced-order
models.

The new formulation provides a desirable framework that is suitable
for the development of efficient computational methods and reduced-order
models as described in the following sections of this paper. It can
deal with the inherent difficulties raised in the original formulation
with irreducible equations. Instead of adding \emph{ad hoc} approximations
for the unresolved higher moments (such as the data-driven model in
\cite{qi2023data}), the crucial third moments are captured explicitly
through the ensemble estimation of the stochastic modes. In addition,
we aim to control computational cost by only running a very small
ensemble for limited samples $M_{1}\ll M$ without sacrificing accuracy
(through the efficient random batch method introduced next in Section
\ref{sec:Computational-methods-rbm}) compared with the original direct
MC approach which demands a large sample size $M$. This combined
framework enables flexible modeling of both key leading-order moments
through the statistical equation, and achieves accurate non-Gaussian
higher-order statistics and fat-tailed PDFs in the extreme events.
It is also more advantageous than the mean-fluctuation model using
only \eqref{eq:dyn_mean} and \eqref{eq:dyn_stoc} (as proposed in
\cite{qi2023random}) which often leads to large numerical errors
and unstable dynamics due to the inherent instability in the turbulent
systems (see examples in Figure \ref{fig:instability} of Section
\ref{subsec:one-layer-L-96}). 

\section{Methods for efficient ensemble forecast using random batch approximation\label{sec:Computational-methods-rbm}}

We propose new efficient computational methods to address the inherent
difficulties in complex turbulent systems involving multiscale interaction
terms by solving the coupled stochastic-statistical equations \eqref{coupled_model}.
Here, we describe the general strategy in the new approach for efficient
statistical prediction by running a very small ensemble simulation
of the stochastic fluctuation coefficients enabled by the random batch
method (RBM) approximation. The ideas will be further illustrated
using concrete examples from the L-96 systems next in Section \ref{sec:Numerical-performance}.

\subsection{Random batch method for coupled turbulent systems with a wide spectrum}

It is realized that the most computational demanding part in solving
the coupled equations comes from getting accurate quantification for
the combined nonlinear feedbacks in the stochastic coefficients and
statistical covariance equations \eqref{eq:dyn_stoc} and \eqref{eq:dyn_cov},
which involves high-order coupling terms involving a wide spectrum
of stochastic modes modes $\mathbf{Z}=\left\{ Z_{k}\right\} _{k=1}^{K},K\gg1$
for each single trajectory. Furthermore, in order to resolve the necessary
high-order statistics with desirable accuracy, an extremely large
ensemble $\left\{ \mathbf{Z}^{\left(1\right)},\cdots,\mathbf{Z}^{\left(M\right)}\right\} \in\mathbb{R}^{K\times M}$
(with the sample size $M$) is required for the full high-dimensional
modes. Thus using the direct MC method for the stochastic equations
ends up with a computational cost of $O\left(MK^{3}\right)$ for the
$M$ samples and $K$ modes and cost $O\left(K^{4}\right)$ for the
covariance equation. The $K^{3}$ cost is due to the quadratic interactions
(of total number $K^{2}$ for each mode) and then for all $K$ modes.
In addition, the required ensemble size $M$ to maintain sufficient
accuracy in the empirical statistical estimation \eqref{eq:pdf_MC}
grows with an exponential rate dependent on the dimension $K$. This
is known as the curse of dimensionality \cite{daum2003curse,friedman1997bias}
and sets an inherent obstacle for effective ensemble prediction of
high-dimensional systems, especially when non-Gaussian statistics
amplifies the demand for an even larger ensemble to capture PDF tails.

Here, we propose to design a computational efficient model using random
batch method (RBM) developed in \cite{jin2020random,jin2021convergence}.
The idea is proposed for the special mean-fluctuation systems in \cite{qi2023random}
concerning only the coupling between the mean and fluctuation modes.
Now, we aim to develop an effective practical strategy for the fully
coupled multiscale states using the general stochastic-statistical
formulation \eqref{coupled_model}. The crucial issue in constructing
the RBM approximation lies in devising an efficient estimation of
the nonlinear cross-interaction terms $Z_{m}Z_{n}$ in the stochastic
and covariance equations \eqref{eq:dyn_stoc} and \eqref{eq:dyn_cov}
without running a very large ensemble. This becomes especially important
in modeling high dimensional turbulent systems where a key feature
is the nonlocal coupling of multiscale states involving a large number
of fluctuation modes $\mathbf{Z}$. 

In the main idea of the new RBM approach, we no longer compute the
expensive nonlinear interactions among the entire stochastic coefficients
$1\leq k\leq K$ in the stochastic and covariance equations. At the
start of each time updating step $t=t_{s}$, a partition $\left\{ \mathcal{I}_{q}^{s}\right\} $
of the mode index is introduced with $\cup_{q}\mathcal{I}_{q}^{s}=\left\{ k:1\leq k\leq K\right\} $,
where each batch only contains a small portion of $\left|\mathcal{I}_{q}^{s}\right|=p$
elements randomly drawn from the total $K$ indices (ending up with
$\left\lceil \frac{K}{p}\right\rceil $ batches). Accordingly, the
full spectrum of modes is also randomly divided into small batches
$\mathcal{Z}_{q}^{s}=\left\{ Z_{k}\left(t\right),k\in\mathcal{I}_{q}^{s}\right\} $
with $\cup_{q}\mathcal{Z}_{q}^{s}=\left\{ Z_{k}\left(t\right):1\leq k\leq K\right\} $
in the time interval $t\in\left(t_{s},t_{s+1}\right]$. Then, instead
of taking summation over all the wavenumbers $m,n$ in the summation
terms of \eqref{eq:dyn_stoc} and \eqref{eq:dyn_cov}, only a small
portion of the modes $Z_{k}\in\mathcal{Z}_{q}^{s}$ with indices in
the batch $k\in\mathcal{I}_{q}^{s}$ are used to update the mode $Z_{k}$
during the time updating interval. The exhausting procedure to resolve
all high-order feedbacks is effectively avoided through this simple
RBM decomposition. The entire high-dimensional system is then decomposed
into smaller subsystems for modes $\left\{ Z_{k},R_{kl}\right\} $
constrained inside each small batch $\left\{ k,l\right\} \in\mathcal{I}_{q}^{s}$
rather than the entire spectral space. The resulting \emph{RBM model
for the stochastic and statistical equations} during the time interval
$t\in\left(t_{s},t_{s+1}\right]$ becomes\addtocounter{equation}{0}\begin{subequations}\label{rbm_model} 
\begin{align}
\frac{\mathrm{d}Z_{k}^{\left(i\right)}}{\mathrm{d}t} & =\sum_{m,n\in\mathcal{I}_{q}^{s}}\tilde{L}_{v,km}\left(\bar{\mathbf{u}}\right)Z_{m}^{\left(i\right)}+\tilde{\gamma}_{mnk}\left(Z_{m}^{\left(i\right)}Z_{n}^{\left(i\right)}-R_{mn}\right)+\sigma_{k}\dot{W}_{k}^{\left(i\right)},\label{eq:rbm_stoch}\\
\frac{\mathrm{d}R_{kl}}{\mathrm{d}t} & =\sum_{m,n\in\mathcal{I}_{q}^{s}}\tilde{L}_{v,km}\left(\bar{\mathbf{u}}\right)R_{ml}+R_{km}\tilde{L}_{v,ml}^{*}\left(\bar{\mathbf{u}}\right)+Q_{\sigma,kl}\nonumber \\
 & +\sum_{m,n\in\mathcal{I}_{q}^{s}}\tilde{\gamma}_{mnk}\frac{1}{M_{1}}\sum_{i=1}^{M_{1}}Z_{m}^{\left(i\right)}Z_{n}^{\left(i\right)}Z_{l}^{\left(i\right)}+\epsilon^{-1}\left(\frac{1}{M_{1}}\sum_{i=1}^{M_{1}}Z_{k}^{\left(i\right)}Z_{l}^{\left(i\right)}-R_{kl}\right).\label{eq:rbm_stat}
\end{align}
\end{subequations}Above, all the wavenumbers $k,l,m,n\in\mathcal{I}_{q}^{s}$
belong to the same self-consistent random batch $\mathcal{Z}_{q}^{s}$.
Importantly, new coupling coefficients, $\tilde{L}_{v,km}=c_{p}^{1}L_{v,km},\tilde{\gamma}_{mnk}=c_{p}^{2}\gamma_{mnk}$
(with the original coefficients $L_{v,km}$ and $\gamma_{mnk}$ defined
in \eqref{coupled_model}), are introduced in the new equations to
guarantee consistent statistics with the scaling factors $c_{p}^{1}=\frac{K}{p}$
and $c_{p}^{2}=\frac{K\left(K-1\right)}{p\left(p-1\right)}$. It is
based on the idea that the random batch approximation should yield
an equivalent effect in the summation of nonlinear coupling terms
under probability expectation on the batch partition $\mathcal{I}^{s}$.
Detailed explanation for the rescaling coefficients $c_{p}^{1,2}$
can be found in the proof of Theorem \ref{thm:converg_coeff} and
\cite{qi2023random}. A much smaller ensemble size is used for the
stochastic modes $\mathbf{Z}^{\left(i\right)},i=1,\cdots,M_{1}$ in
\eqref{eq:rbm_stoch} to give empirical estimation of the higher-order
moments in \eqref{eq:rbm_stat}. Notice that in \eqref{eq:rbm_stoch}
within the time interval, each sample is updated independently. After
this time updating interval $\left(t_{s},t_{s+1}\right]$, the batches
are resampled at the start of the new time step $t=t_{s+1}$ to repeat
the same procedure, so the modes from different batches get mixed.

The statistical consistency of the RBM approximation \eqref{rbm_model}
is guaranteed by the ergodicity and fast mixing of the high wavenumber
modes. First, considering the typical property of the turbulent modes,
the energy inside the single small-scale mode $R_{k}=\left\langle \left|Z_{k}\right|^{2}\right\rangle $
decays fast as $k$ grows large and de-correlates rapidly in time.
Second, ergodicity of the stochastic fluctuation modes $\mathbf{Z}$
implies that updating the key statistics using fractional fluctuation
modes at each time step with consistent time-averaged feedback can
provide an equivalent total contribution. Therefore, the total $K$
spectral modes are divided into smaller batches to be updated individually
during each time interval. As a result, rather than running a large
ensemble of high dimensional solutions of the full fluctuation modes
$\mathbf{Z}^{\left(i\right)}$ as in the direct MC approach of \eqref{coupled_model},
only a small number of stochastic trajectories are needed as long
as it is sufficient to sample the $p$-dimensional modes inside each
batch rather than the full $K$-dimensional space. We show the rigorous
convergence in leading-order statistics of this RBM approximation
next in Section \ref{subsec:Error-analysis}.

As the first major reduction of the above RBM model, instead of using
the entire spectrum of modes $1\leq k\leq K$ to update each wavenumber
mode, we only consider the nonlinear interactions between modes in
the a very small subset $\mathcal{I}_{q}^{s}$ of size $p$. As shown
in the numerical examples in Section \ref{subsec:one-layer-L-96},
the batch size $p$ can be picked as a very small number $p=O\left(1\right)$.
This leads to the effective computational reduction from $O\left(K^{3}\right)$
in direct MC to $O\left(p^{2}K\right)$ for time evolution of each
single stochastic trajectory, and from $O\left(K^{4}\right)$ to $O\left(p^{2}K^{2}\right)$
for the covariance equation. In addition, through the RBM, only a
very small number of modes are contained in each batch for the estimation
of the nonlinear coupling. Consequently, we don't need to compute
all the cubic terms $Z_{m}Z_{n}Z_{l}$ but only the modes inside one
batch of size $p\ll K$ independent of the full dimension. This leads
to the further significant reduction in the required sample size from
$O\left(MK^{3}\right)$ to $O\left(M_{1}p^{2}K\right)$ in the ensemble
prediction. Especially, the ensemble size $M_{1}$ only needs to sample
the $p$-dimensional batch subspace in contrast to $M$ in the full
MC approach to sample the $K$-dimensional full space. The algorithm
using random batch method for ensemble simulation of high dimensional
turbulent system under the coupled stochastic-statistical closure
model is summarized in Algorithm \ref{alg:full-RBM}.

\begin{algorithm}[H]
\begin{algorithmic}[1] 
\Require{At initial time $t=0$, assign the initial mean and covariance for $\left\{ \bar{\mathbf{u}}_{0},R_{0}\right\}$ and draw samples $\mathbf{Z}_{0}$ from the initial distribution $\mu_{0}$.}

\For{$s = 1$ while $s \leq \left\lceil T/\Delta t\right\rceil$, at the start of the time interval $t\in\left(t_{s},t_{s+1}\right]$  with time step $\Delta t=t_{s+1}-t_{s}$}
	\State{Partition the $K$ modes into $S$ batches (with $pS=K$) randomly as $\mathcal{Z}_{q}^{s}=\left\{ Z_{k}\left(t_{s}\right),k\in\mathcal{I}_{q}^{s}\right\}$ with $\cup_{q=1}^{S}\mathcal{Z}_{q}=\mathbf{Z}$.}
	\State{Update $\left\{ Z_{k}\left(t_{s+1}\right),R_{kl}\left(t_{s+1}\right)\right\}$ for modes in batch $k,l \in \mathcal{I}^{s}_{q}$ independently according to \eqref{eq:rbm_stoch} and \eqref{eq:rbm_stat}.}
    \State{Update $\bar{\mathbf{u}}\left(t_{s+1}\right)$ according to the original mean equation \eqref{eq:dyn_mean} using all the batch outputs.}
\EndFor

\end{algorithmic}

\caption{Full RBM approximation for the coupled multiscale stochastic-statistical
model\label{alg:full-RBM}}
\end{algorithm}

\subsection{Reduced-order model for efficient ensemble simulations}

In the above model with the RBM approximation, we still need to run
an ensemble simulation for the entire spectral modes $\left\{ Z_{k}^{\left(i\right)}\right\} _{k=1}^{K}$
even though with a much small number of samples $i=1,\cdots,M_{1}$,
which reaches the final computational cost dependent on the full dimension
$K$. On the other hand, in practical situations we are mostly interested
in the statistics in a much smaller number of leading modes $1\leq k\leq K_{1}\ll K$
(such as the largest scales or the most energetic modes) rather than
the entire spectrum. This leads to the second major approximation
to introduce the effective reduced-order modeling strategy focusing
on the leading dominant modes while still taking into account the
contributions from the large number of small-scale fluctuation modes
through the effective RBM approximation allowing an even smaller sample
size.

To achieve the model reduction, we utilize the idea introduced in
\cite{qi2023random} and extend it to the general coupled stochastic-statistical
model \eqref{coupled_model}. The key idea is still to notice that
the large number of fast mixing small-scale modes make an equivalent
contribution through time average, thus we can decompose them into
batches for updating different ensemble samples of the central large-scale
modes. We further introduce an explicit large-small scale decomposition
for the original fluctuation modes $\mathbf{Z}=\left(\mathbf{X},\mathbf{Y}\right)$,
where $\mathbf{X}=\left\{ X_{k}\right\} $ represents the small number
of leading modes (such as $\left\{ Z_{k},k\leq K_{1}\right\} $) while
$\mathbf{Y}=\left\{ Y_{l}\right\} $ are all the rest large number
of smaller scale fluctuation modes (such as all the rest modes $\left\{ Z_{l},K_{1}<l\leq K\right\} $).
We aim to use an ensemble empirical approximation for the marginal
distribution of $\mathbf{X}$, that is,
\begin{equation}
p_{t}^{\mathrm{RBM}}\left(\mathbf{X}\right)\coloneqq\frac{1}{M_{1}}\sum_{i=1}^{M_{1}}\prod_{k=1}^{K_{1}}\delta\left(X_{k}-X_{k}^{\left(i\right)}\left(t\right)\right).\label{eq:pdf_rbm}
\end{equation}
By focusing on the leading modes in a much lower dimension $K_{1}\ll K$,
the required ensemble size $M_{1}$ is further reduced. Still, each
large-scale mode $X_{k}^{\left(i\right)}$ is coupled to the unresolved
small-scale modes $Y_{l}$ through the linear and nonlinear coupling
terms. Using the model reduction strategy to be combined with the
RBM approximation, we no longer run ensemble simulation for the large
number of small-scale fluctuation modes $\mathbf{Y}$. Instead, only
one (or at most a small number) of stochastic trajectory $\mathbf{Y}\left(t\right)$
is solved in time. The total $K-K_{1}$ spectral modes in $\mathbf{Y}$
are then divided into small batches to update different ensemble members
of $\mathbf{X}^{\left(i\right)},i=1,\cdots,M_{1}$. This leads to
$M_{1}$ batches from the small-scale modes with the relation $pM_{1}=K-K_{1}$.
The idea here is to use the large number of small-scale modes to update
different large-scale ensemble samples at each time step. 

In this way, we decompose the original stochastic equation of independent
samples \eqref{eq:rbm_stoch} into a coupled system with interacting
samples in large-scale states $\mathbf{X}^{\left(i\right)}=\left\{ Z_{k}^{\left(i\right)}\right\} _{k=1}^{K_{1}},i=1,\cdots,M_{1}$
and one single trajectory of the small-scale state $\mathbf{Y}=\left\{ Z_{l}\right\} _{l=K_{1}}^{K}$
(with $K_{1}\ll K$). At the start of time step $t=t_{s}$, each large-scale
sample $\mathbf{X}^{\left(i\right)}$ is grouped with one batch of
the small-scale modes $\mathcal{Y}_{i}^{s}=\left\{ Y_{l}\left(t_{s}\right)\right\} _{l\in\mathcal{I}_{i}^{s}}$,
where $\mathcal{I}_{i}^{s}$ is the RBM partition of the small-scale
modes with $\cup_{i}\mathcal{I}_{i}^{s}=\left\{ l:K_{1}\leq l\leq K\right\} $.
Therefore, we get the \emph{coupled reduced-order RBM equations} for
the $i$-th large-scale ensemble sample grouped with a batch of the
small-scale modes $\left\{ X_{k}^{\left(i\right)},Y_{l}\right\} _{k\leq K_{1},l\in\mathcal{I}_{i}^{s}}$
during the time updating interval $t\in\left(t_{s},t_{s+1}\right]$
\begin{equation}
\begin{aligned}\frac{\mathrm{d}X_{k}^{\left(i\right)}}{\mathrm{d}t} & =F_{k}\left(\mathbf{X}^{\left(i\right)},\mathbf{Y}\right)+\sum_{m,n\leq K_{1}}\tilde{L}_{v,km}\left(\bar{\mathbf{u}}\right)X_{m}^{\left(i\right)}+\tilde{\gamma}_{mnk}\left(X_{m}^{\left(i\right)}X_{n}^{\left(i\right)}-R_{mn}\right)+\sigma_{k}\dot{W}_{k}^{\left(i\right)},\quad1\leq k\leq K_{1},\\
\frac{\mathrm{d}Y_{l}}{\mathrm{d}t} & =G_{l}\left(\mathbf{X}^{\left(i\right)},\mathbf{Y}\right)+\sum_{m,n\in\mathcal{I}_{i}^{s}}\tilde{L}_{v,lm}\left(\bar{\mathbf{u}}\right)Y_{m}+\tilde{\gamma}_{mnl}\left(Y_{m}Y_{n}-R_{mn}\right)+\sigma_{l}\dot{W}_{l},\qquad l\in\mathcal{I}_{i}^{s}.
\end{aligned}
\label{eq:rbm_stoch2}
\end{equation}
Above, $F_{k}$ and $G_{l}$ contain the residual terms that represent
the cross-coupling between the large and small scale states. Usually,
the state decomposition and the coupling dynamics will become straightforward
according to the specific dynamical structure of the multiscale system.
The batches are then resampled each time at the start of the new time
updating cycle $t=t_{s+1}$ same as the previous RBM strategy. It
needs to be noticed that through the above model approximation \eqref{eq:rbm_stoch2},
different ensemble samples $\mathbf{X}^{\left(i\right)}$ are no longer
independent since they are linked by the shared small-scale trajectory
of $\mathbf{Y}$. This is a reasonable assumption from the common
observation in turbulent systems that the large number of small-scale
fluctuation modes can be viewed as almost independent random processes
with a rapidly decaying energy spectrum. Next in Section \ref{subsec:two-layer-L-96},
we will show the construction of reduced-order model through one explicit
example from the two-layer L-96 system.

Using the above decomposition, we can effectively reduce the computational
cost by only sampling a much smaller space of dimension $K_{1},$
thus enabling an even smaller sample size. The final computational
cost is then reduced to $O\left(M_{1}K_{1}^{3}+\left(K-K_{1}\right)p^{2}\right)=O\left(M_{1}\left(K_{1}^{3}+p^{3}\right)\right)$
(using the above partition relation for small-scale modes, $pM_{1}=K-K_{1}$)
rather than the cost $O\left(M_{1}Kp^{2}\right)$ in the full RBM
model involving a very large number of $K$. The ensemble size $M_{1}$
is then only sampling the very small large-scale modes of dimension
$K_{1}=O\left(1\right)$ rather than the full state dimension $K$,
together with the rest dimension $K-K_{1}$ only related to the small
batch size $p=O\left(1\right)$. In this way, the curse of dimensionality
is effectively avoided. Similarly, we summarize the reduced-order
RBM strategy in the following Algorithm \ref{alg:red-RBM}.

\begin{algorithm}[H]
\begin{algorithmic}[1] 
\Ensure{Decompose the stochastic modes $\mathbf{Z}=\mathbf{X}+\mathbf{Y}$ into large scales $\mathbf{X}\in\mathbb{R}^{K_1}$ and small scales $\mathbf{Y}\in\mathbb{R}^{K-K_1}$ with $K_{1}\ll K$.}

\Require{At initial time $t=0$, assign the initial mean and covariance for $\left\{ \bar{\mathbf{u}}_{0},R_{0}\right\}$. 
Draw samples for the larg-scale modes $\mathbf{X}^{\left(i\right)}_{0},i=1,\cdots,M_{1}$ and small-scale modes $\mathbf{Y}_{0}$ from the initial distribution.}

\For{$s = 1$ while $s \leq \left\lceil T/\Delta t\right\rceil$, at the start of the time interval $t\in\left(t_{s},t_{s+1}\right]$  with time step $\Delta t=t_{s+1}-t_{s}$}
\State{Partition the $K-K_1$ small-scale modes into $M_1$ batches randomly. The $i$-th large-scale sample is grouped with the small-scale modes in one batch as $\mathcal{Z}_{i}^{s}=\left\{ \mathbf{X}^{\left(i\right)},Y_{l}\right\}_{l\in\mathcal{I}_{i}^{s}}$.}
	\State{Update $\left\{ \mathbf{X}^{\left(i\right)}\left(t_{s+1}\right),Y_{l}\left(t_{s+1}\right)\right\}$ for each batch $l \in \mathcal{I}^{s}_{i}$ according to \eqref{eq:rbm_stoch2}.}
	\State{Update the statistical mean and covariance to $\bar{\mathbf{u}}\left(t_{s+1}\right), R\left(t_{s+1}\right)$ accordingly based on \eqref{eq:dyn_mean} and \eqref{eq:rbm_stat} using all batch outputs.}

\EndFor
\end{algorithmic}

\caption{Reduced-order RBM approximation for the coupled multiscale stochastic-statistical
model\label{alg:red-RBM}}
\end{algorithm}

\subsection{Error analysis for the RBM model approximation\label{subsec:Error-analysis}}

Here, we provide convergence analysis on the RBM model Algorithm \ref{alg:full-RBM}
in \eqref{rbm_model} compared with the solution obtained from the
full stochastic-statistical equations \eqref{coupled_model}. Through
this analysis, we rigorously demonstrate the effectiveness of the
proposed RBM as a precise approximation to the direct MC solutions.
Besides, it also provides error estimates and guidelines for selecting
appropriate model parameters in implementing the computational schemes.

First, we consider the convergence of the RBM approximation $\tilde{Z}_{k}$
of the stochastic equation \eqref{eq:rbm_stoch} to the full stochastic
equation \eqref{eq:dyn_stoc} of the coefficients $Z_{k}$\addtocounter{equation}{0}\begin{subequations}\label{stoch_analy}
\begin{align}
\frac{\mathrm{d}\tilde{Z}_{k}}{\mathrm{d}t} & =\sum_{m,n\in\mathcal{I}_{k}^{s}}\tilde{L}_{v,km}\left(\bar{u}\right)\tilde{Z}_{m}+\tilde{\gamma}_{mnk}\left(\tilde{Z}_{m}\tilde{Z}_{n}-R_{mn}\right)+\sigma_{k}\dot{\tilde{W}}_{k},\label{eq:coeff_rbm}\\
\frac{\mathrm{d}Z_{k}}{\mathrm{d}t} & =\sum_{m,n\leq K}L_{v,km}\left(\bar{u}\right)Z_{m}+\gamma_{mnk}\left(Z_{m}Z_{n}-R_{mn}\right)+\sigma_{k}\dot{W}_{k},\label{eq:coeff_full}
\end{align}
\end{subequations}with the RBM parameters $\tilde{L}_{v,km}=\frac{K}{p}L_{v,km}$
and $\tilde{\gamma}_{mnk}=\frac{K\left(K-1\right)}{p\left(p-1\right)}\gamma_{mnk}$
and $p$ the batch size. We need the following structure assumption
for the model parameters:
\begin{assumption*}
Suppose that the quasi-linear coupling coefficients in the stochastic
equations are uniformly bounded
\begin{equation}
\max_{k,m}L_{km}\left(\bar{u}\right)\leq C,\quad\max_{k,m,n}\gamma_{mnk}\leq C.\label{eq:assump1}
\end{equation}
And the quadratic coupling term satisfies the symmetry
\begin{equation}
B\left(\mathbf{u}_{k},\mathbf{u}_{k}\right)=0,\quad\left[B\left(\mathbf{u}_{k},\mathbf{u}_{l}\right)+B\left(\mathbf{u}_{l},\mathbf{u}_{k}\right)\right]\cdot\mathbf{u}_{k}=0,\quad\mathrm{for\:any\:}k,l\leq K.\label{eq:assump2}
\end{equation}
\end{assumption*}
In the above equations \eqref{eq:coeff_rbm} and \eqref{eq:coeff_full},
we assume that the full and RBM model have consistent mean and covariance
$\bar{u},R$ from the statistical equations. In this way, we are able
to focus on the approximation due to the random batches $\mathcal{I}_{k}^{s}$.
Conditional on the accurate statistical mean state and covariance
$\bar{u}$ and $R$, the RBM approximation adds additional randomness
through the indices of the modes. The model dynamics fit into the
SDE systems discussed in \cite{jin2020random,jin2021convergence},
where the interacting particles are replaced by the coupled spectral
modes in the above turbulent multiscale model. We can find the statistical
convergence in the stochastic coefficients following the similar argument
as in \cite{qi2023random}.
\begin{thm}
\label{thm:converg_coeff}Under the assumptions \eqref{eq:assump1},
the statistical estimation of the RBM model \eqref{eq:coeff_rbm}
with time step $\Delta t$ converges to the statistics of the full
model \eqref{eq:coeff_full} up to the final time $T$ as
\begin{equation}
\sup_{s\Delta t\leq T}\left|\frac{1}{K}\sum_{k=1}^{K}\mathbb{E}\varphi\left(\tilde{Z}_{k}^{s}\right)-\frac{1}{K}\sum_{k=1}^{K}\mathbb{E}\varphi\left(Z_{k}^{s}\right)\right|\leq C_{\varphi}\left(T\right)\Delta t,\label{eq:converg_coeff}
\end{equation}
with the test function $\varphi\in C_{b}^{2}$, and $\tilde{Z}_{k}^{s}=\tilde{Z}_{k}\left(t_{s}\right),Z_{k}^{s}=Z_{k}\left(t_{s}\right)$
the solutions at $t_{s}=s\Delta t$. $C_{\varphi}$ is the coefficient
independent of the model dimension $K$.
\end{thm}

The proof of Theorem \ref{thm:converg_coeff} compares the backward
equations for the two equations \eqref{eq:coeff_rbm} and \eqref{eq:coeff_full}
and uses the fact that the expectation on the random batch samples
gives back the original equation (which leads to the precise forms
of the rescaling factors $\tilde{L}_{v},\tilde{\gamma}$). The statistical
consistency in the stochastic and statistical equations are guaranteed
by the additional relaxation term. We give the proof of the theorem
in Appendix \ref{sec:Proof-of-theorems}. Especially, using the above
conclusion and choosing $\varphi\left(Z\right)=Z^{2}$, we find the
convergence of the total variance in the ensemble estimate of the
stochastic equations. 
\begin{cor}
The total variance estimation in the RBM model \eqref{rbm_model}
converges to the true statistics as
\begin{equation}
\frac{1}{N}\sum_{s=1}^{N}\left|\mathrm{tr}\tilde{R}\left(t_{s}\right)-\mathrm{tr}R\left(t_{s}\right)\right|\leq\sup_{s\Delta t\leq T}\left|\mathrm{tr}\tilde{R}_{s}-\mathrm{tr}R_{s}\right|\leq C_{v}\left(T\right)K\Delta t,\label{eq:coverg_var}
\end{equation}
where $\mathrm{tr}\tilde{R}\left(t_{s}\right)=\mathbb{E}\left|\tilde{Z}_{k}^{s}\right|$
and $\mathrm{tr}R\left(t_{s}\right)=\mathbb{E}\left|Z_{k}^{s}\right|$
are the trace of the variance from the second moments of the stochastic
coefficients.
\end{cor}

Second, we consider the error in the statistical mean state using
the simplified version of the scalar mean equation for RBM model and
the full model\addtocounter{equation}{0}\begin{subequations}\label{mean_analy} 
\begin{align}
\frac{\mathrm{d}\tilde{\bar{u}}}{\mathrm{d}t} & =-\lambda\tilde{\bar{u}}+B\left(\tilde{\bar{u}},\tilde{\bar{u}}\right)+\sum_{k,l\leq K}\tilde{R}_{kl}B\left(\mathbf{v}_{k},\mathbf{v}_{l}\right)+F,\label{eq:mean_rbm}\\
\frac{\mathrm{d}\bar{u}}{\mathrm{d}t} & =-\lambda\bar{u}+B\left(\bar{u},\bar{u}\right)+\sum_{k,l\leq K}R_{kl}B\left(\mathbf{v}_{k},\mathbf{v}_{l}\right)+F.\label{eq:mean_full}
\end{align}
\end{subequations}In the above statistical mean equations, we assume
only linear damping in the linear term and consider the error from
the RBM approximation of the second moments. $R$ and $\tilde{R}$
are the covariance matrix estimated in \eqref{eq:converg_coeff}.
The idea is to use the total statistical energy equation \cite{majda2018strategies}
which provides a balanced estimation for the statistical in both mean
and total variance, $E=\bar{u}^{2}+\frac{1}{K}\mathrm{tr}R$. The
assumption \eqref{eq:assump2} provides a very simple equation describing
the evolution of the total statistical energy $E$. Using the estimate
for the total variance, we find the error in the RBM prediction in
the mean state.
\begin{thm}
\label{thm:converg_mean}Under the assumptions \eqref{eq:assump2},
the statistical mean estimation for the RBM model \eqref{eq:mean_rbm}
converges to the full model solution \eqref{eq:mean_full} with time
step $\Delta t$ and final time $T$ as
\begin{equation}
\left\Vert \tilde{\bar{u}}-\bar{u}\right\Vert _{\left[0,T\right]}\leq\sup_{s\Delta t\leq T}\left|\tilde{\bar{u}}\left(t_{s}\right)-\bar{u}\left(t_{s}\right)\right|\leq C_{m}\left(T\right)\Delta t.\label{eq:converg_mean}
\end{equation}
Above, the error in mean state is taking over the time average $\left\Vert f\right\Vert _{\left[0,T\right]}^{2}=\frac{1}{N}\sum_{s=1}^{N}f^{2}\left(t_{s}\right)$
among solutions at each time evaluation step $t_{s}=s\Delta t$ up
to the final time $T=N\Delta t$.
\end{thm}

Again, we put the proof of the above theorem in Appendix \ref{sec:Proof-of-theorems}.
The estimates in \eqref{eq:coverg_var} and \eqref{eq:converg_mean}
quantifies the approximation errors dependent on the time step size
$\Delta t$. In practice, it is always easier to apply a smaller time
step to effectively increase the prediction accuracy. The results
in the mean and covariance convergence still rely on the consistency
in the stochastic and statistical equations \eqref{eq:dyn_stoc} and
\eqref{eq:dyn_cov} in the coupled model. This is implicitly guaranteed
by the essential relaxation term as $\epsilon\rightarrow0$. We leave
the complete analysis of the coupled stochastic-statistical equations
\eqref{coupled_model} as a mean field system \cite{calvello2022ensemble}
in future research.

\section{Numerical performance on prototype models: the Lorenz '96 systems\label{sec:Numerical-performance}}

In evaluating the performance of the proposed methods, we start with a simple prototype model which is still able to capture
key properties of the general turbulent systems. The Lorenz '96 (L-96)
systems provide a desirable testbed exhibiting a range of representative
statistical features such as interaction of variables from different
scales and non-Gaussian statistics. The model is originally introduced
to study the mid-latitude weather in a simpler one-layer form \cite{lorenz1996predictability},
and is later generalized to a two-layer form to include strong multiscale
spatiotemporal coupling \cite{wilks2005effects}. In this section,
we illustrate the application of the full and reduced-order RBM methods
for predicting leading statistics and non-Gaussian PDFs using the
one-layer and two-layer L-96 systems.

\subsection{The one-layer L-96 system\label{subsec:one-layer-L-96}}

\begin{figure}
\subfloat{\includegraphics[scale=0.48]{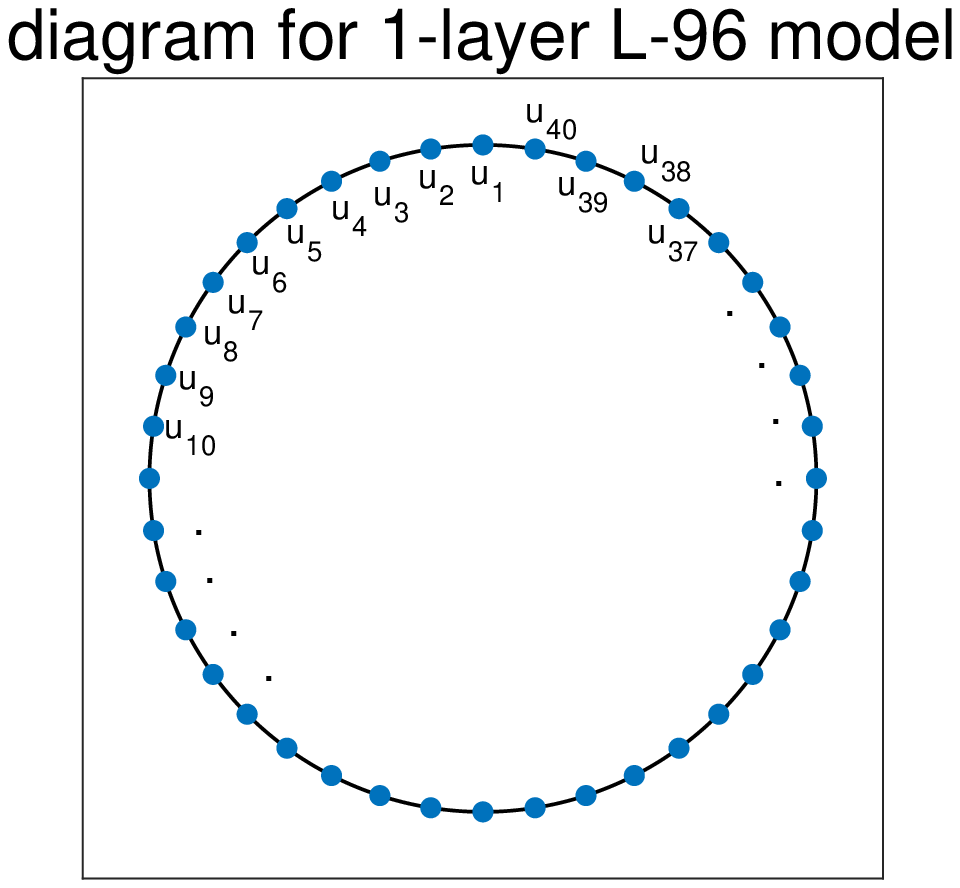}}

\vspace{-1em}\subfloat{\includegraphics[scale=0.48]{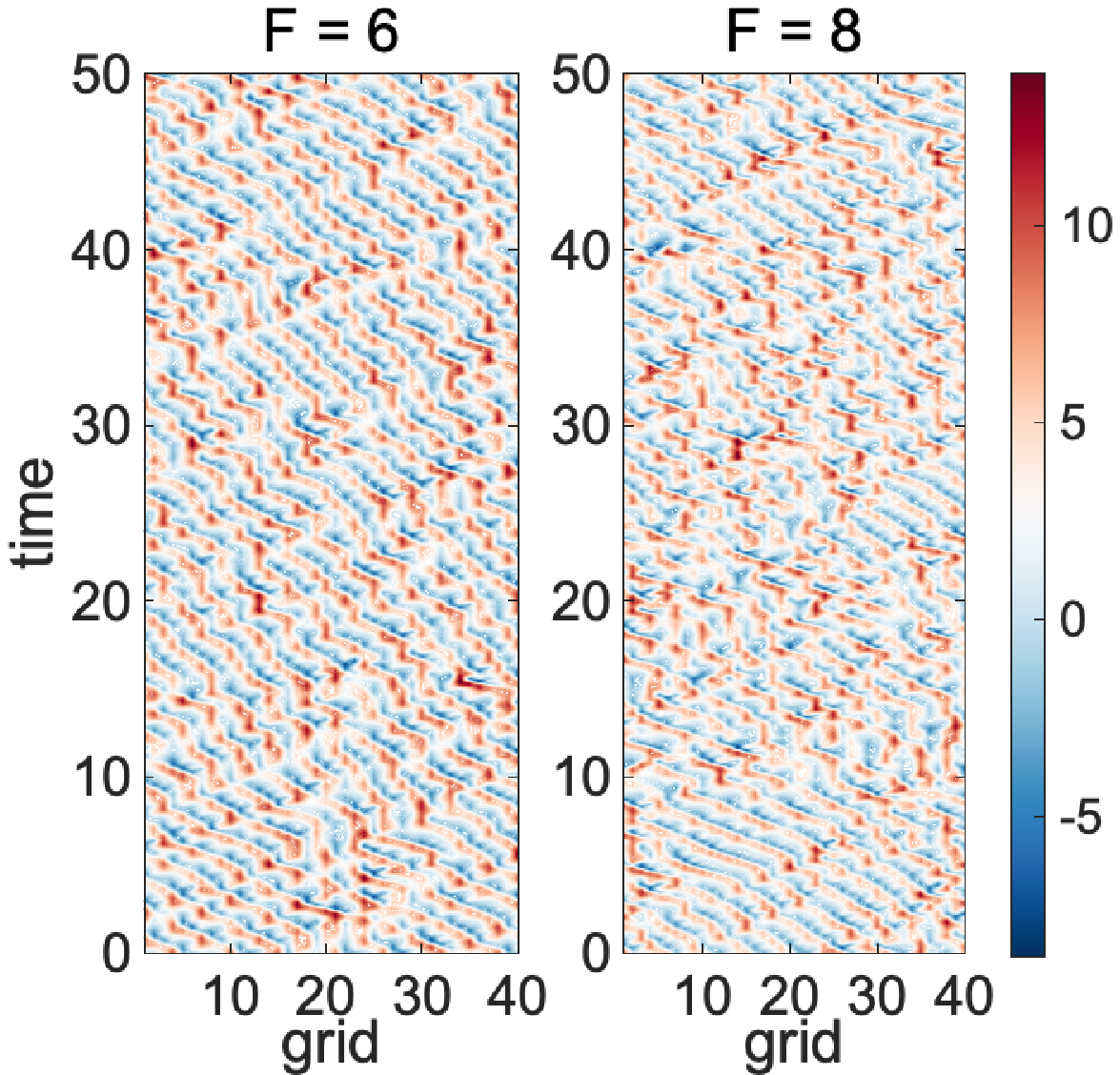}}

\caption{Illustration of the model structure and typical solutions of the one-layer
L-96 system \eqref{eq:l96_homo}.\label{fig:Illustration-l96}}

\end{figure}

First, we consider the one-layer L-96 system \cite{lorenz1996predictability}
that can be expressed as a 40-dimensional ODE system with homogeneous damping and forcing 
\begin{equation}
\frac{\mathrm{d}u_{j}}{\mathrm{d}t}=\left(u_{j+1}-u_{j-2}\right)u_{j-1}-u_{j}+F,\;j=1,\cdots,J=40,\label{eq:l96_homo}
\end{equation}
with $u_{J+1}=u_{1}$ and constant uniform damping and forcing. The model state is defined with periodic boundary condition $u_{j+J}=u_{j}$
mimicking geophysical waves at $J$ equally distributed locations
along a constant mid-latitude circle (see the diagram in Figure \ref{fig:Illustration-l96}).
The ODE system has a moderate dimension $J=40$.  Various representative statistical
features comparable with the real data from observations can be generated in the L-96 solutions \eqref{eq:l96_homo} by simply
changing the constant forcing $F$. A smaller forcing $F$ is related
to a weaker mean state and stronger non-Gaussian statistics, while
a large forcing value $F$ leads to stronger mixing and near-Gaussian
statistics \cite{majda2016improving}. As illustrated in the typical solution trajectories plotted
in Figure \ref{fig:Illustration-l96}, the one-layer L-96 solutions
display the distinctive dynamical transition from more regular wave
patterns (with $F=6$) to more chaotic flow features (with $F=8$).
The model also maintains a spectrum of large number of unstable and
energetic stochastic modes (see Figure \ref{fig:instability} and
\ref{fig:Time-evolutions}). Accordingly, the PDFs of spectral modes demonstrate
the statistical transition from non-Gaussian distributions to near-Gaussian ones
between the two regimes (see Figures \ref{fig:1D-marginal-PDFs} and \ref{fig:2D-joint-PDFs}
for the marginal and joint PDFs). Thus, it established a challenging and
important test case for accurate prediction of model statistics and
PDFs.

\subsubsection{Explicit formulation and RBM approximation for the one-layer L-96 system}

Here, we derive the explicit stochastic-statistical formulation for the one-layer L-96
system as an example to illustrate the main ideas and performance
in the proposed strategy. In the L-96 system \eqref{eq:l96_homo},
randomness comes only from the initial condition and is amplified
due to the internal instability. We seek the probability solution
of the model state $p_{t}\left(\left\{ u_{j}\right\} \right)$ starting
from an initial distribution $p_{0}\left(\left\{ u_{j}\right\} \right)$.
Clearly, the L-96 system fits into the general framework \eqref{eq:abs_formu}.
The periodic boundary condition implies the mean-fluctuation decomposition
as in \eqref{eq:decomp} by a Galerkin projection on the Fourier basis
\begin{equation}
u_{j}=\bar{u}_{j}+\frac{1}{J}\sum_{\left|k\right|\leq J/2}Z_{k}\left(t\right)e^{i2\pi k\frac{j}{J}},\label{eq:modes_1layer}
\end{equation}
where $\bar{u}_{j}=\left\langle u_{j}\right\rangle _{p_{t}}$ is the statistical
mean and $\left\{ Z_{k}\right\} $ are stochastic coefficients with
$\left\langle Z_{k}\right\rangle _{p_{t}}=0$. Due to the constant
forcing and damping terms, the resulting statistics in each moment
of the solution state is translation invariant \cite{majda2018strategies,qi2022machine}. It implies that
the mean state, $\left\langle u_{j}\right\rangle _{p_{t}}=\bar{u}\left(t\right)$, is uniform;
and the off-diagonal covariance entries are all vanishing, $R\left(t\right)=\left\langle \mathbf{Z}\mathbf{Z}^{*}\right\rangle _{p_{t}}=\mathrm{diag}\left\{ r_{k}\left(t\right)\right\} $.
Therefore, we can focus on the dynamical equations for the scalar mean $\bar{u}$ and variance
$r_{k}$ of each Fourier mode together with the stochastic coefficients $Z_{k}$
to recover the entire statistics in this system. 

The explicit full stochastic-statistical formulation \eqref{coupled_model}
and the corresponding RBM model \eqref{rbm_model} for the L-96 system
are listed in Appendix \ref{subsec:Explicit-RBM-1layer}. Due to the
still relatively low dimension $J=40$ of the one-layer model, we
apply the full RBM model in Algorithm \ref{alg:full-RBM} and test the random batch approximation for the full spectrum prediction involving a large number of highly unstable and energetic 
modes. We summarize several main features observed from using the
RBM model on the one-layer L-96 system before showing the detailed discussions on the numerical results in the next section:
\begin{itemize}
\item The computational cost for the full model \eqref{eq:full_l96} is
$O\left(J+MJ^{2}+MJ^{2}\right)=O\left(J\left(1+2MJ\right)\right)$
where $J$ is the full dimension of the system and $M$ is the ensemble
size to sufficiently sample the $J$-dimensional space. In the RBM
approximation \eqref{eq:rbm_l96}, the computational cost is effectively
reduced to $O\left(J\left(1+2M_{1}p\right)\right)$ with $p$ the
batch size and $M_{1}\ll M$ is the sample size only required to sample
the much smaller $p$-dimensional subspace in each batch. Especially,
$p$ does not increase as the dimensional $J$ increases.
\item The coupled formulation combining the variance equations for $r_{k}$
and stochastic equations for $Z_{k}$ with the relaxation factor $\epsilon$
is essential to reach the correct final equilibrium state. Using purely
the stochastic equations for $Z_{k}$ is insufficient to recover the
correct statistics when the sample size becomes small due to the strong model
instability.
\item The RBM approximation relies on the ergodicity and fast mixing of
the small-scale modes. The L-96 system possesses a wide spectrum
with large number of fluctuation modes showing distinctive time scales, making it a difficult test case for the RBM model. Still, accurate statistical
prediction is achieved even with an extremely small batch size $p=2$.
Prediction results can be further improved by reducing the time step
size $\Delta t$.
\end{itemize}

\subsubsection{Numerical tests on the one-layer L-96 model}

In the numerical tests, we take two typical regimes of the L-96 system
with $F=6$ (showing non-Gaussian statistics) and $F=8$ (showing
near-Gaussian statistics). Th standard 4th-order Runge-Kutta scheme
is adopted for the time integration with  time step size $\Delta t=1\times10^{-4}$.
A large ensemble size $M=1\times10^{5}$ is needed to capture the
true model statistics accurately from the direct MC simulation. In
the RBM prediction, different batch sizes $p=2,5,10$ are used to compute the high-order interactions 
compared with the full model using the entire $J=40$ modes. Therefore, only a very
small sample size $M_{1}=100$ becomes sufficient, allowing for efficient computation.

\begin{figure}
\includegraphics[scale=0.38]{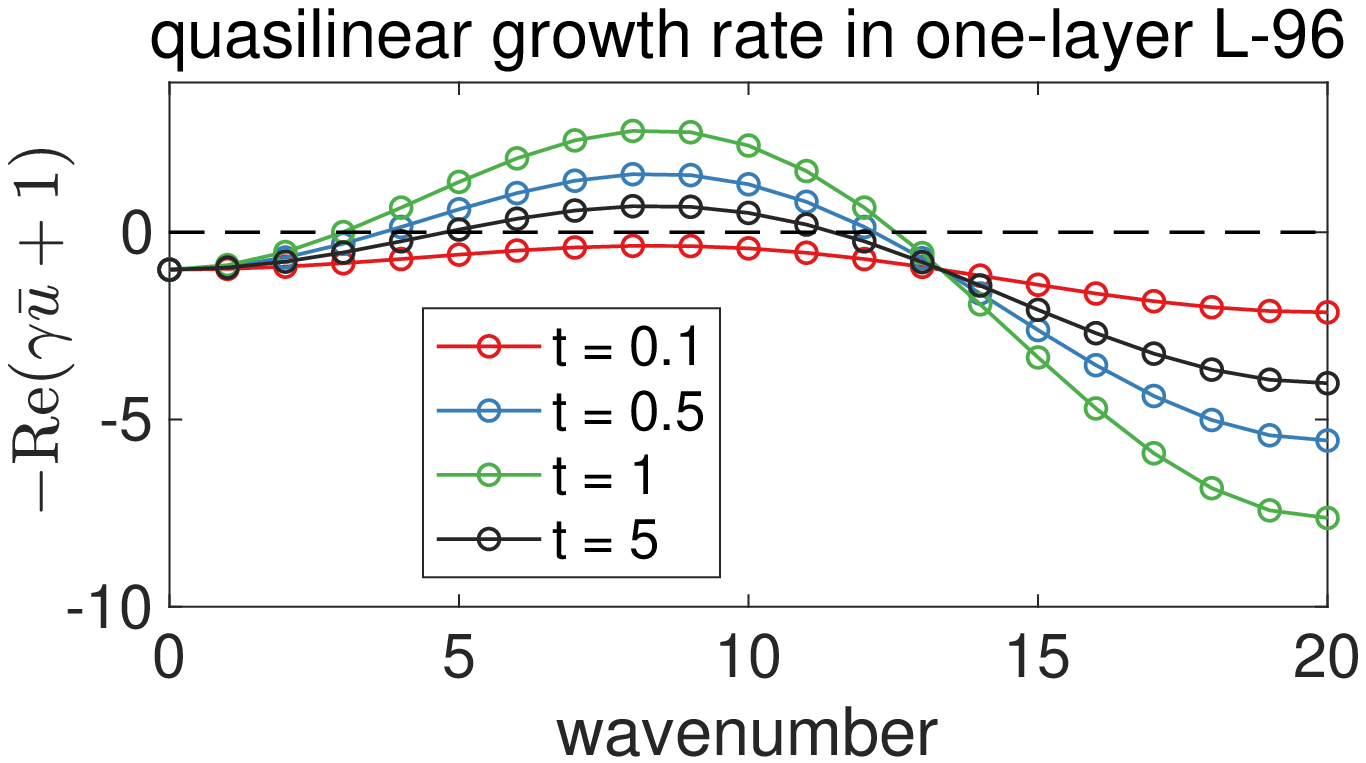}\includegraphics[scale=0.38]{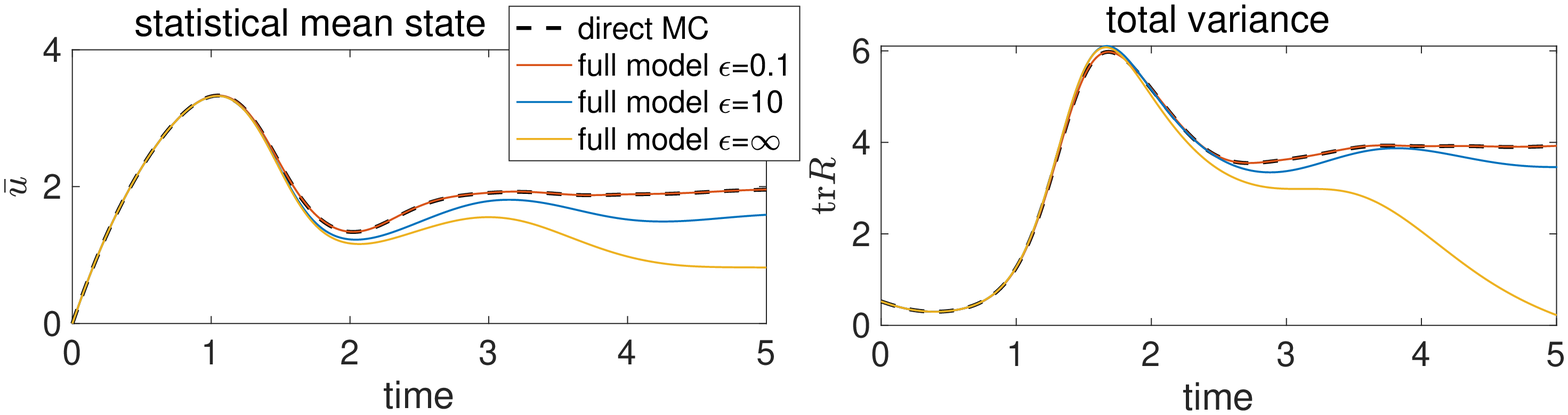}

\caption{Effects of instability in the one-layer L-96 system with $F=6$. Left:
quasilinear growth rate in each wavenumber $-\mathrm{Re}\left(\gamma_{k}^{*}\bar{u}+1\right)$
at several time instants; Right: direct solutions from the
full model with different relaxation strengths $\epsilon$.\label{fig:instability}}
\end{figure}

First, we illustrate the inherent difficulty in running ensemble prediction
for turbulent systems with instability. In the full model \eqref{eq:full_l96}, the
variance equations for $r_{k}$ (as well as the stochastic equations
for $Z_{k}$) are subject to internal instability due to the real part
of the coupling coefficient $-\left(\gamma_{k}^{*}\bar{u}+1\right)$
representing interaction with the mean state $\bar{u}$. As illustrated
in the left panel of Figure \ref{fig:instability}, positive growth
rates are induced in a large number of modes as the system evolves in 
time. This indicates the crucial role of the combined third order
moments $\left\langle Z_{m}Z_{n}^{*}Z_{k}^{*}\right\rangle $ acting
as a balancing factor for these internally unstable modes. On the
other hand, with insufficient sample size, large errors can be introduced
to the empirical estimation of the higher-order feedback term. The inherent instability in the turbulent model formulation can be also seen in the zeroth mode equation \eqref{eq:eqn_z0}.
It shows that the equation is only marginally stable thus small errors
in the ensemble estimation will lead to large errors.

In the practical simulation of the ensemble scheme, the internal instabilities
will amplify the small errors and lead to disastrous result. The situation
will become increasingly serious if we only want to use a small ensemble
size. As a typical example shown on the right panel of Figure \ref{fig:instability},
we get the truth from the direct MC simulation of \eqref{eq:l96_homo} with an extremely large
ensemble size $M=1\times10^{5}$. In comparison, we run the full model
\eqref{eq:full_l96} with several different values $\epsilon$
in the additional relaxation term $\epsilon^{-1}\left(\left\langle \left|Z_{k}\right|^{2}\right\rangle -r_{k}\right)$.
It shows that even with this moderate dimension $J=40$ and very large
ensemble size, the solution for the statistical mean and variance
will diverge without the relaxation term $\epsilon=\infty$, while
equilibrium consistency is guaranteed when smaller value of $\epsilon$
is added. It demonstrates the crucial role of the additional relaxation
term to guarantee equilibrium converge. It also shows that the
model is robust with consistent final statistics for a wide range of values of $\epsilon$ as long as it is not too large.
\begin{rem*}
As a simple example to illustrate the internal instability, the explicit
equation for the zeroth mode $Z_{0}$ in the full model \eqref{eq:full_l96}
gives
\begin{equation}
\frac{\mathrm{d}Z_{0}}{\mathrm{d}t}=\frac{1}{J}\sum_{k}\gamma_{k}\left(\left|Z_{k}\right|^{2}-r_{k}\right)-Z_{0},\label{eq:eqn_z0}
\end{equation}
with $\gamma_{k}=\cos\frac{4\pi k}{J}-\cos\frac{2\pi k}{J}$. The
consistent final equilibrium requires $\left\langle Z_{0}\right\rangle =\frac{1}{J}\sum_{k}\gamma_{k}\left(\left\langle \left|Z_{k}\right|^{2}\right\rangle -r_{k}\right)=0$,
while instability will be introduced through errors due to the empirical
average $\left\langle \left|Z_{k}\right|^{2}\right\rangle \sim r_{k}$.
This will lead to inherent difficulty even with a very large sample
size as shown in Figure \ref{fig:instability}.
\end{rem*}

\

\begin{figure}
\centering
\subfloat{\includegraphics[scale=0.38]{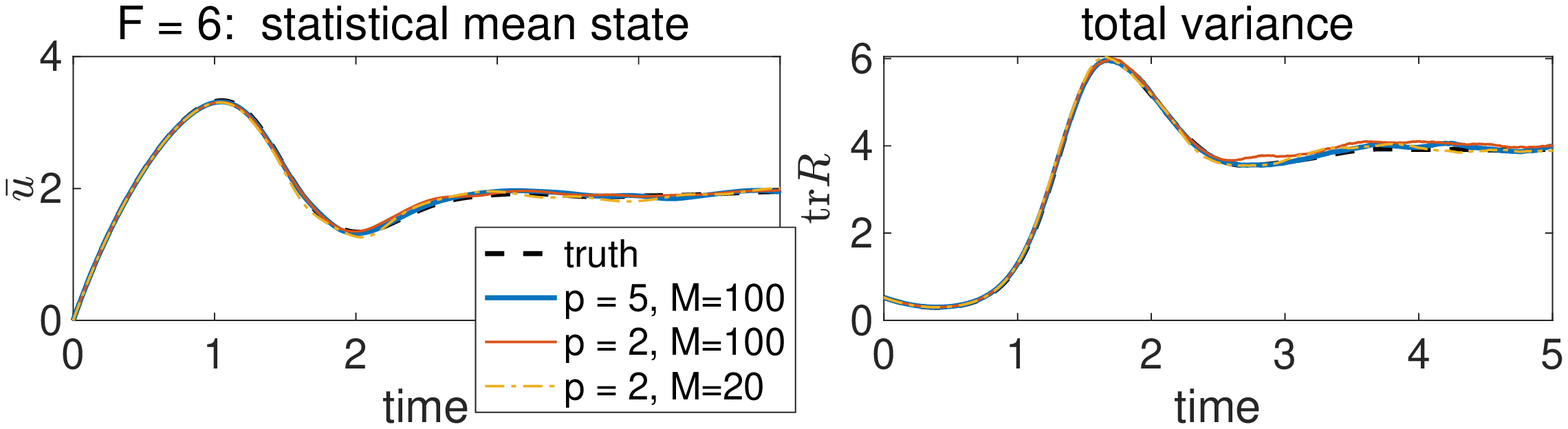}\includegraphics[scale=0.38]{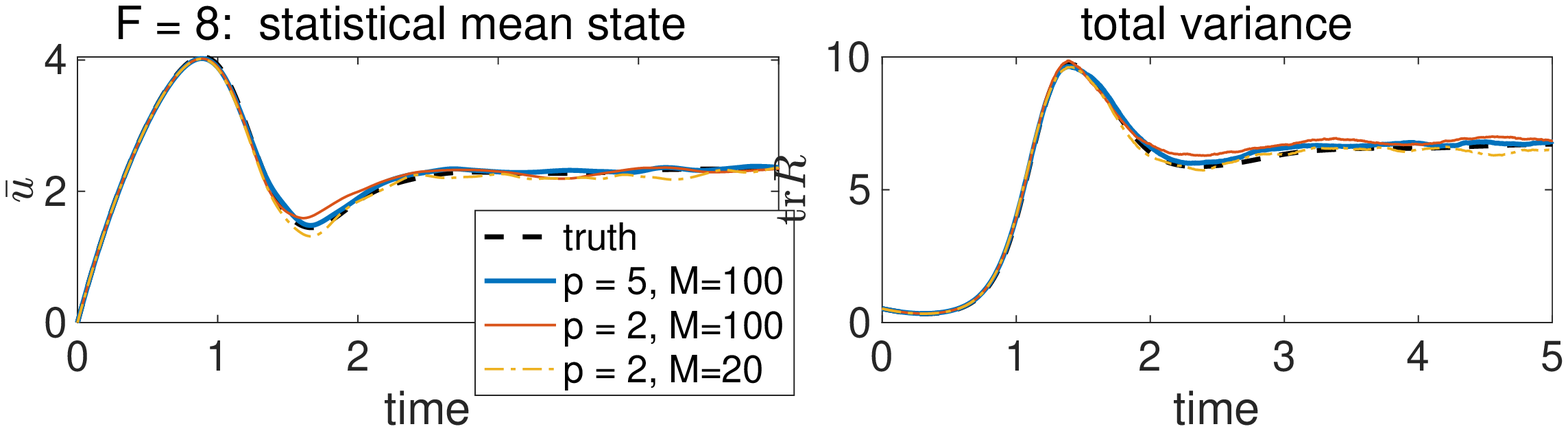}}

\subfloat{\includegraphics[scale=0.38]{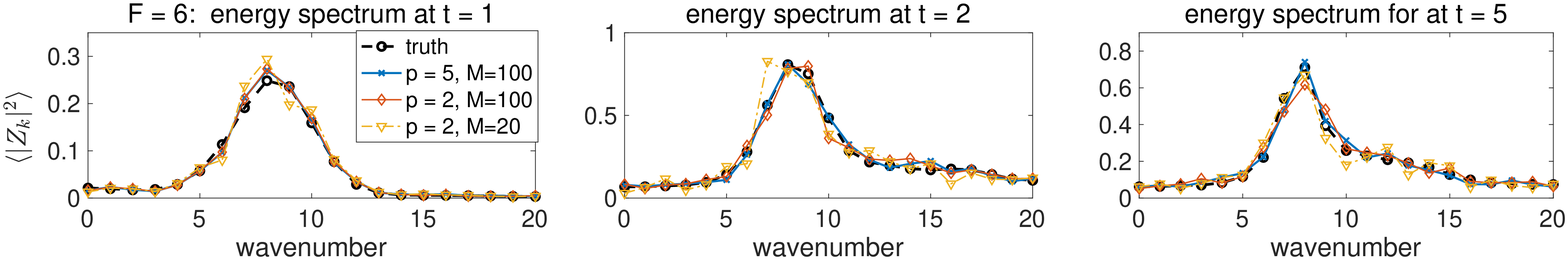}}

\subfloat{\includegraphics[scale=0.38]{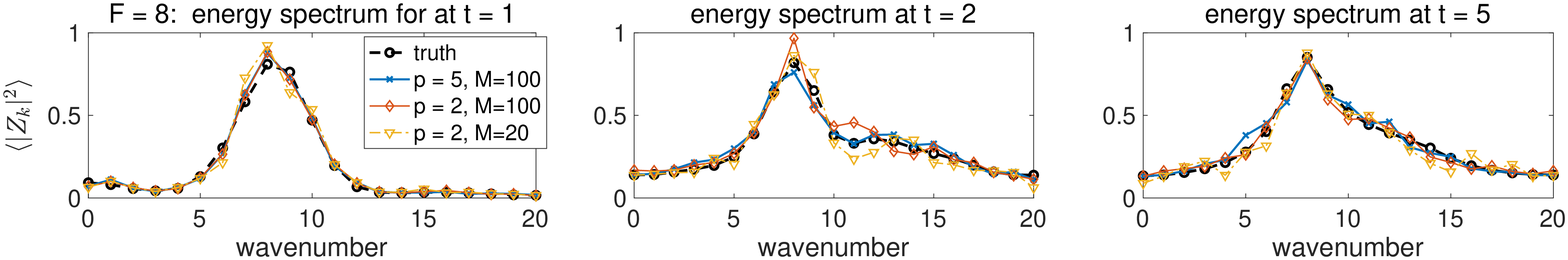}}

\caption{Time evolutions of the mean and total variance as well as the variance
spectra from the RBM prediction with batch size $p=2,5$ and ensemble size $M_{1}=100,20$.
Results from the two test regimes $F=6$ and $F=8$ are shown and
the truth is recovered from a direct MC simulation with a large sample size $M=1\times10^{5}$.\label{fig:Time-evolutions}}
\end{figure}

Next, we check the performance of the RBM prediction with different
batch sizes. In the one-layer L96 model with $J=40$, we will focus
on the fully resolved RBM model \eqref{eq:rbm_l96} described in Algorithm
\ref{alg:full-RBM}. In Figure \ref{fig:Time-evolutions}, the time-series
of the mean $\bar{u}$ and total variance $\mathrm{tr}R=\sum r_{k}$
as well as the prediction for the detailed  variance spectrum in each mode during
the time evolution are plotted. Two batch sizes $p=5,2$ with $M_{1}=100$
samples from the RBM model prediction are compared with the truth
from the direct MC simulation using $M=1\times10^{5}$ samples for
the full dimension $J=40$. Furthermore, we also test an extreme case
with only $p=2$ modes (that is, only using one term in the high-order feedback for the strongly unstable dynamics) in each batch and a very small sample size $M_{1}=20$ in computing the statistics. It shows that all the RBM results
accurately track the truth statistics with lines overlapping on each other while save a lot of computational
cost using the extremely small ensemble. In particular, the L-96 system
sets a challenging test model for the RBM model since it contains
a wider spectrum of energetic modes containing large degrees of instabilities. Still,
as shown in the energy spectra from the starting transient stage to
the final equilibrium, the variances in all the modes are captured
with high precision.

\begin{figure}
\centering
\subfloat{\includegraphics[scale=0.3]{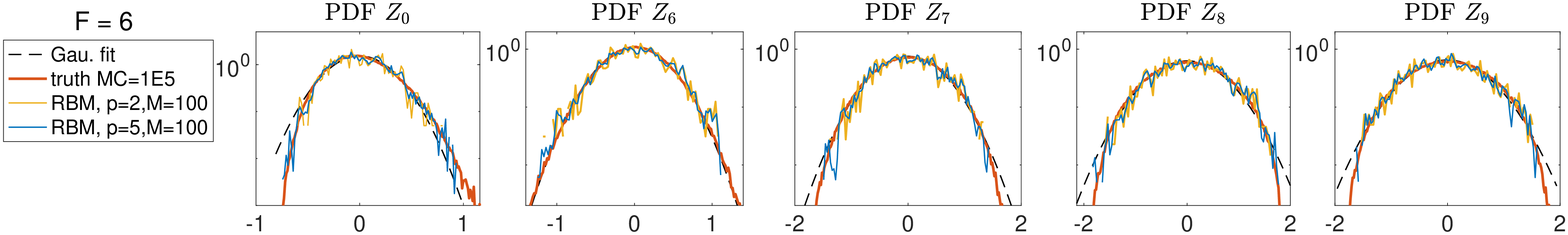}}

\subfloat{\includegraphics[scale=0.3]{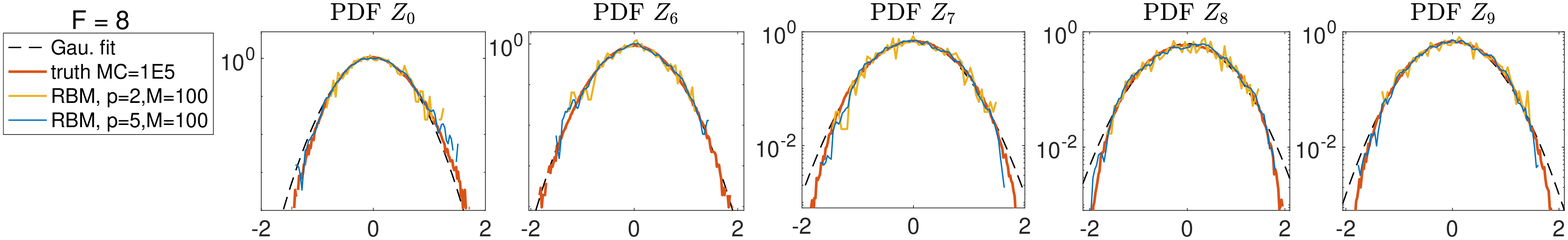}}

\caption{1D marginal PDFs of the stochastic coefficients $Z_{k}$ in the two
test regimes $F=6,8$ of the one-layer L-96 system. The RBM predictions with sample size $M_1=100$ 
are compared with the truth from direct MC simulation with $M=1\times10^5$. The Gaussian
fits of the PDFs with the same mean and variance are also plotted in
dashed lines.\label{fig:1D-marginal-PDFs}}

\end{figure}

\begin{figure}
\subfloat[MC $F=6$]{\includegraphics[scale=0.3]{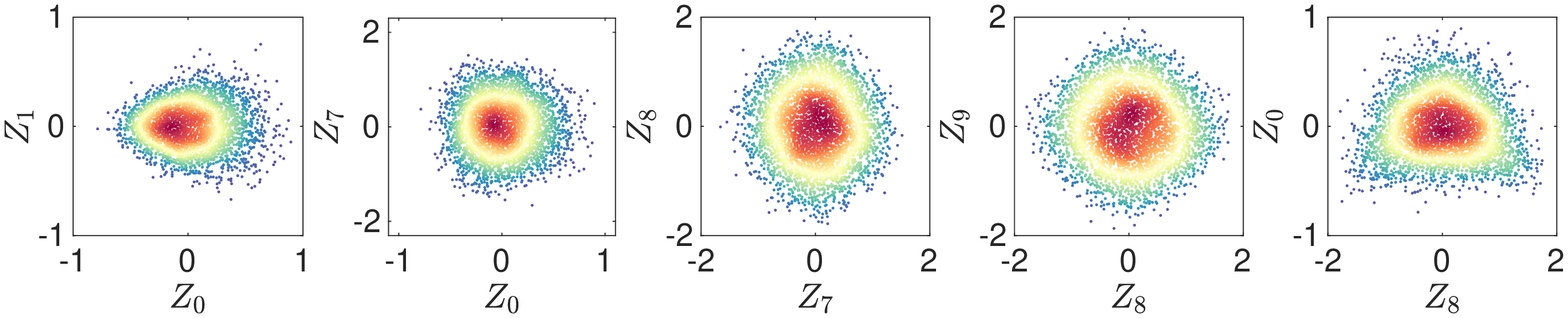}

}\hspace{.5em}\subfloat[MC $F=8$]{\includegraphics[scale=0.3]{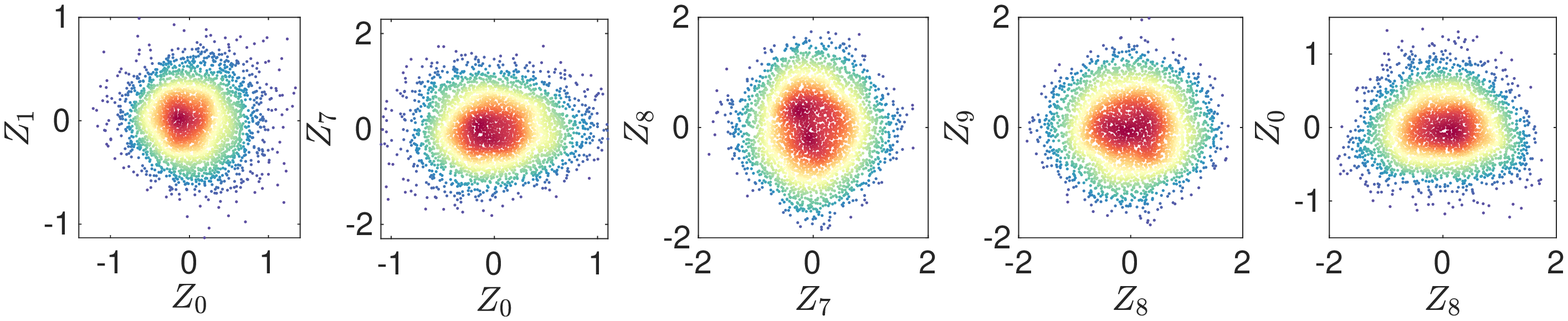}

}

\vspace{-1.em}

\subfloat[RBM $F=6$]{\includegraphics[scale=0.3]{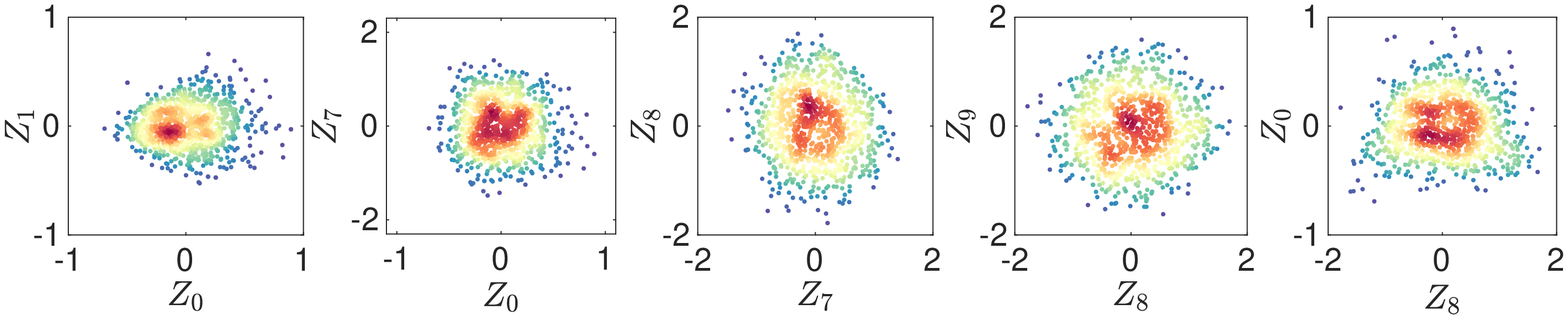}

}\hspace{.5em}\subfloat[RBM $F=8$]{\includegraphics[scale=0.3]{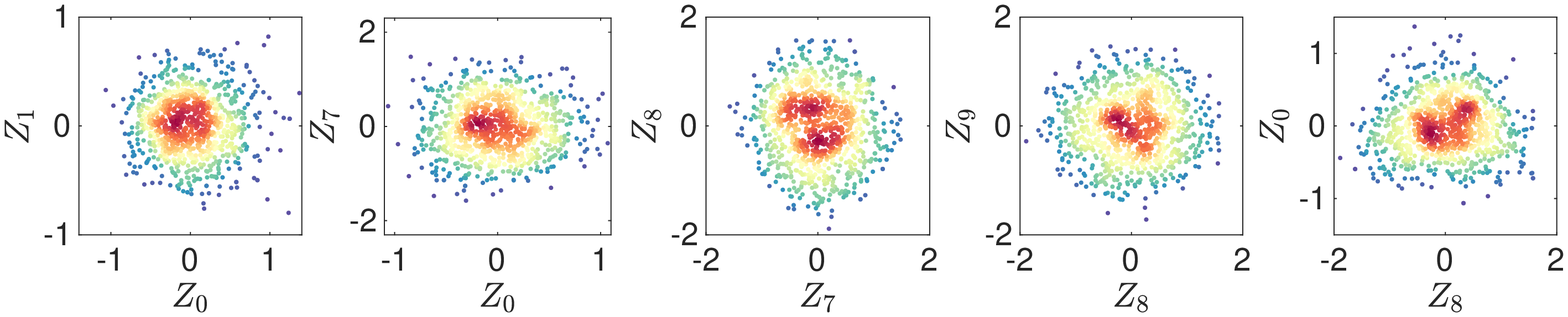}

}

\caption{2D joint PDFs of the leading modes in the two test regimes $F=6,8$
of the one-layer L-96 system. The RBM predictions with sample size $M_1=100$ are compared with
the truth from direct MC simulation with $M=1\times10^5$.\label{fig:2D-joint-PDFs}}

\end{figure}

Then, another important goal of the RBM model is to accurately characterize
the PDFs of the model states. This is especially important in
application to uncertainty quantification and data assimilation, where
an accurate estimate of PDFs using limited number of samples is crucial.
In the RBM approach, the PDFs are captured by the empirical distribution
\eqref{eq:pdf_MC} of the stochastic coefficients $Z_{k}$. Figure
\ref{fig:1D-marginal-PDFs} first plots the 1D marginally PDFs in
the leading Fourier modes $Z_{k}$ of the one-layer L-96 system in the two
test regimes. It shows that the $F=6$ regime generates stronger non-Gaussian
PDFs while the $F=8$ regime is closer to Gaussian but still contains
non-negligible non-Gaussian features. Usually, a very large ensemble is essential
to capture such non-Gaussian statistics in the direct MC approach.
Using the efficient RBM approximation, the shapes PDFs including the
skewed and sub-Gaussian structures are accurately characterized using a
very small ensemble size $M_{1}=100$. As a more precise calibration of the prediction of PDFs, we plot the 2D joint PDFs between
the most important modes in Figure \ref{fig:2D-joint-PDFs}. The non-Gaussian
structures can be observed more clearly in the joint distributions. The outliers
in the sampled PDFs play an important role to characterize the occurrence
of extreme events and have crucial impact in many practical applications
with limited samples. The RBM model successfully captures the joint
PDFs especially recovers the representative non-Gaussian shapes in
the outlier regions though using only a very small number of samples.

\begin{figure}
\centering
\includegraphics[scale=0.35]{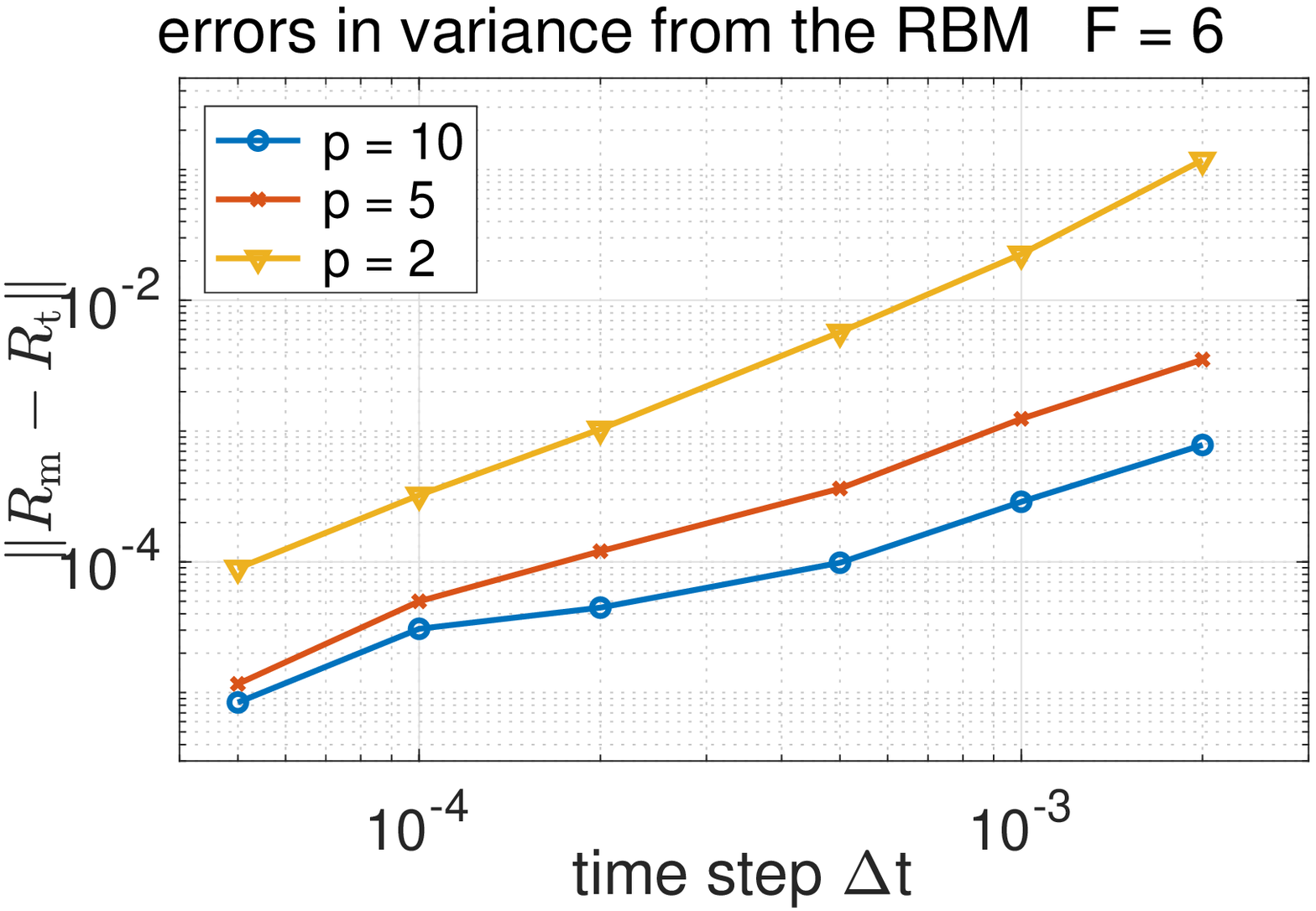}\includegraphics[scale=0.35]{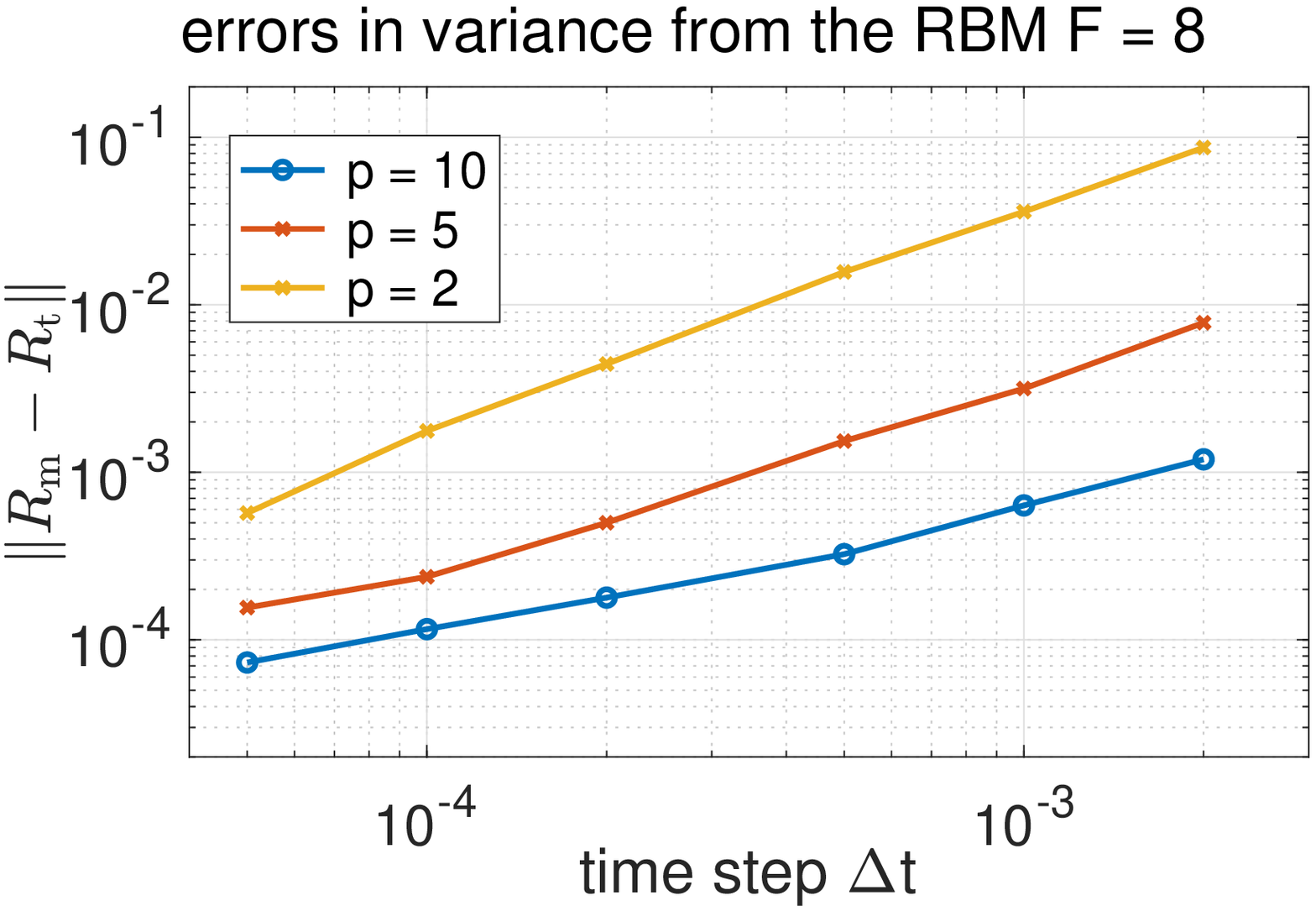}\caption{Errors in the total variance from the RBM prediction with different
batch sizes $p=2,5,10$ in the two test regimes $F=6,8$ of the one-layer L-96 system.\label{fig:Errors-comp}}
\end{figure}

Finally, we provide a quantitative quantification of the model errors
of the RBM prediction with different batch sizes $p$. In particular,
one of the central quantities to estimate in the RBM is the second
moments, thus we consider the time-averaged errors in the variance
prediction $\left\Vert R_{\mathrm{m}}-R_{\mathrm{t}}\right\Vert =\frac{1}{N}\sum_{s=1}^{N}\sum_{k}\left|r_{\mathrm{m},k}\left(t_{s}\right)-r_{\mathrm{t},k}\left(t_{s}\right)\right|$
where $R_{\mathrm{m}}$ is the RBM prediction and $R_{\mathrm{t}}$
is the truth. Figure \ref{fig:Errors-comp} plots the errors with
the batch sizes $p=2,5,10$ and sample size $M_{1}=500$. We use a
relatively larger sample size to reduce the fluctuation errors from
the small ensemble. The convergence of RBM estimates depends on the average on the random batches (thus equivalently on time average of the fast
mixing modes), therefore the errors grow as the time step size $\Delta t$ increases.
A smaller batch size will lead to larger errors due to the fewer high-order
terms explicitly modeled at each time updating step. When the time
step decreases to smaller values, the major source of errors is taken
over by the fluctuation errors from the small ensemble size, so
the decay of the errors starts to slow down. Also the larger batch
size $p=10$ shows a slower decay rate since it needs to sample a
relatively larger dimensional batch subspace in computing the nonlinear
coupling terms. Overall, the error plot shows that the RBM model remains uniformly high prediction
skill in accurately recovering the key statistics with a much lower affordable
computational cost. It also implies that it is a good strategy to
improve the prediction accuracy by taking a smaller time step $\Delta t$
without increasing much of the computational cost.

\subsection{The two-layer L-96 system\label{subsec:two-layer-L-96}}

\begin{figure}
\centering
\hspace{-1.5em}\subfloat{\includegraphics[scale=0.43]{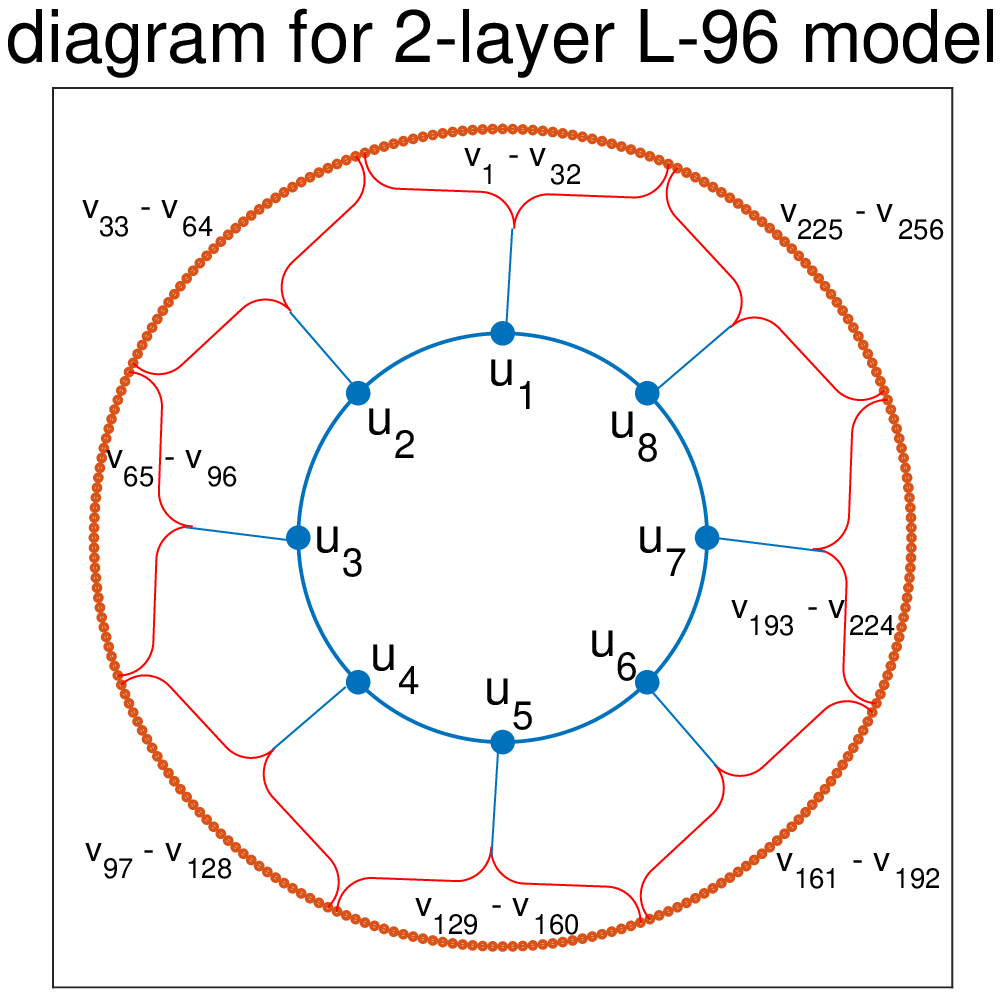}}\hspace{-2.em}\subfloat{\includegraphics[scale=0.43]{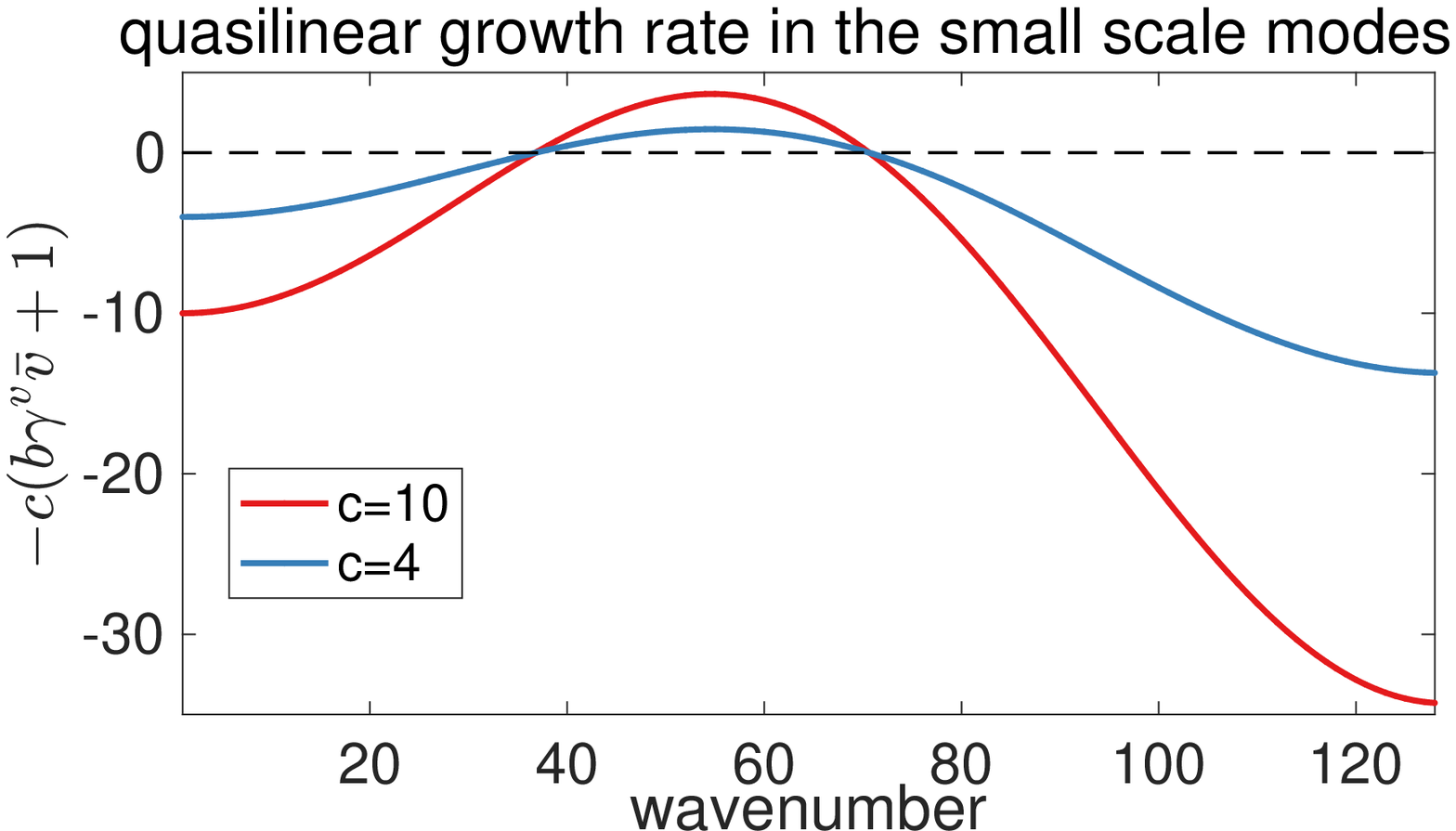}}

\subfloat{\includegraphics[scale=0.43]{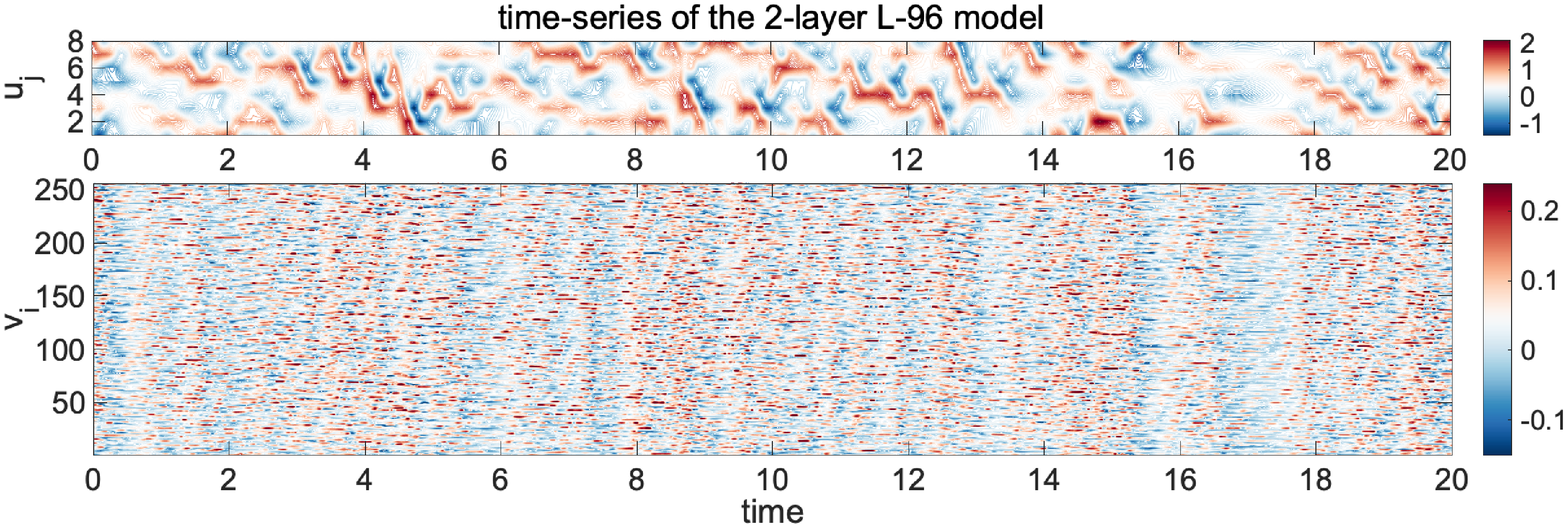}}

\caption{Illustration of the two-layer L-96 system with large-scale size $J=8$
and small-scale size $JL=256$. Upper: diagram for the model structure
and growth rates in the small-scale modes; Lower: a typical 
time-series solution for the large and small scale states with $c=4,b=10,h=1$.\label{fig:Growth-rates-2layer}}

\end{figure}

In the second test model, we examine the model effectiveness
in handling truly high-dimensional systems with multiscale structures by considering the more complicate two-layer L-96 system. As a further generalization
of the original one-layer system, the two-layer L-96 system \cite{arnold2013stochastic,wilks2005effects}
introduces an additional second layer state $v_{i},i=1,\cdots,JL$ to the
first layer state $u_{j},j=1,\cdots,J$ such that 
\begin{equation}
\begin{aligned}\frac{du_{j}}{dt} & =\left(u_{j+1}-u_{j-2}\right)u_{j-1}-u_{j}+F-\frac{hc}{b}\sum_{i=L\left(j-1\right)+1}^{jL}v_{i},\\
\frac{\mathrm{d}v_{i}}{\mathrm{d}t} & =-cb\left(v_{i+2}-v_{i-1}\right)v_{i+1}-cv_{i}+\frac{hc}{b}u_{\left[\frac{i-1}{L}\right]+1},
\end{aligned}
\label{eq:l96_2layer}
\end{equation}
both with periodic boundary conditions, $u_{j+J}=u_{j}$ and $v_{i+JL}=v_{i}$.
Above, $u$ is usually referred to as the large-scale slow variables and
$v$ as the small-scale fast variables. On the right hand sides of
\eqref{eq:l96_2layer}, the double layer states follow the same energy-conserving nonlinear
self-coupling structure as well as uniform linear damping effect as in the
one-layer case. The two states of different scales are then coupled
through three additional model parameters, $c,b,h$: $c$ signifies
the time-scale separation; $b$ controls the ratio between the amplitudes
of two layer states $u_{j}$ and $v_{i}$; and $h$ characterizes
the coupling strength between the two states.
We illustrate the coupling structure of the two-layer
system in the diagram of Figure \ref{fig:Growth-rates-2layer}.  In particular, the second
layer states $v_{i}$ for $i=1,\cdots,JL$ are locally coupled with
one corresponding first layer state $u_{j}$ with the index $j=\left[\frac{i-1}{L}\right]+1$
(where $\left[a\right]$ takes the integer part of $a$), and are
globally linked with each other by the nonlinear self-interactions. Inversely,
each first layer state $u_{j}$ receives the combined feedback from
a sequence of $L$ second layer states $v_{i}$, $i=L\left(j-1\right)+1,\cdots,jL$. This leads to
a fully coupled high-dimensional system including the two-level states
$\left\{ u_{j},v_{i}\right\} $ with a total of $J\left(L+1\right)$
state variables.  The multiscale structure of the
dynamical solution is demonstrated by a typical time-series of
the large and small scale processes. It is clear to observe the distinctive
time and spatial scales and close correlation between the two scales.

The model parameters used in the numerical tests are listed in Table~\ref{tab:Model-parameters-2layer}
taken the standard model setup as in \cite{wilks2005effects,arnold2013stochastic,giggins2019stochastically}.
Especially, we consider two typical parameter regimes with a mild time-scale separation
$c=4$ and a strong time-scale separation $c=10$. Again, the two-layer
L-96 system displays strong internal instability containing a large number of unstable
modes with positive growth rates as shown in Figure
\ref{fig:Growth-rates-2layer}. To recover the true model statistics
with sufficient accuracy, we need to take a small time step $\Delta t=1\times10^{-4}$
in order to completely resolve the extremely fast time scales in the
second layer variables $v_{i}$. A forth-order Runge-Kutta scheme
is adopted for the time integration in both the direct MC simulation
and the RBM approaches. A very large ensemble size $M=5\times10^{5}$
is required considering the very high-dimension $d=J\left(L+1\right)=264$
of the system. In contrast, it is sufficient to only use $M_{1}=O\left(100\right)$
samples in both the full and reduced-order RBM approaches.

\begin{table}
\centering
\begin{tabular}{cccccccccc}
\toprule 
$J$ & $L$ & $F$ & $c$ & $b$ & $h$ &  & $\Delta t$ & $M$ & $M_{1}$\tabularnewline
\midrule
\midrule 
8 & 32 & 20 & 4, 10 & 10 & 1 &  & $1\times10^{-4}$ & $5\times10^{5}$ & 500, 100\tabularnewline
\bottomrule
\end{tabular}

\caption{Model parameters used for the two-layer L-96 model \eqref{eq:l96_2layer}\label{tab:Model-parameters-2layer}}

\end{table}

\subsubsection{Stochastic-statistical formulation for the two-layer L-96 system}

Similar to the one-layer L-96 system, we can project the large and
small scale states to the Fourier modes due to the periodic boundary
condition for both $u_{j}$ and $v_{i}$ as
\begin{equation}
u_{j}=\bar{u}+\frac{1}{J}\sum_{\left|k\right|\leq J/2}Z_{k}\left(t\right)e^{i2\pi k\frac{j}{J}},\quad v_{i}=\bar{v}+\frac{1}{JL}\sum_{\left|l\right|\leq JL/2}Y_{l}\left(t\right)e^{i2\pi l\frac{i}{JL}},\label{eq:modes_2layer}
\end{equation}
where $Z_{k}$ and $Y_{l}$ are the stochastic coefficients for the
large and small scale states correspondingly. The two-layer L-96 system
\eqref{eq:l96_2layer} also accepts the general system structure \eqref{eq:abs_formu}.
Using the translation invariance in both large and small scales, the
statistical prediction aims to recover the homogeneous mean states
$\bar{u}\equiv\left\langle u_{j}\right\rangle _{p_{t}}$ and $\bar{v}\equiv\left\langle v_{i}\right\rangle _{p_{t}}$,
and the variances in large and small scale modes $r_{k}^{u}=\left\langle \left|Z_{k}\right|^{2}\right\rangle _{p_{t}},r_{l}^{v}=\left\langle \left|Y_{l}\right|^{2}\right\rangle _{p_{t}}$
according to the model probability measure $p_{t}$. The detailed
equations involving various linear and nonlinear coupling effects
between different scales become very complicated. We list the explicit mean and covariance equations involving complicated
higher-order coupling terms across the scales in Appendix \ref{subsec:RBM-equations-2layer}. 

The most important ideas concerning the RBM approximation of the two-layer L96 system \eqref{eq:l96_2layer} can be illustrated
in the stochastic equations for the large and small scale fluctuation
modes 
\begin{equation}
\begin{aligned}\frac{\mathrm{d}Z_{k}}{\mathrm{d}t} & =\frac{1}{J}\sum_{m-n=k}\left(Z_{m}Z_{n}^{*}-r_{m}^{u}\delta_{mn}\right)\gamma_{mn}^{u}-d_{u}Z_{k}-\frac{1}{L}\sum_{s=-L/2+1}^{L/2}\lambda_{k+sJ}^{*}Y_{k+sJ},\\
\frac{\mathrm{d}Y_{l}}{\mathrm{d}t} & =\frac{1}{JL}\sum_{p-q=l}\left(Y_{p}Y_{q}^{*}-r_{p}^{v}\delta_{pq}\right)\gamma_{pq}^{v}-d_{v}Y_{l}+\lambda_{l}Z_{\mathrm{mod}\left(l,J\right)}.
\end{aligned}
\label{eq:fluc_2layer}
\end{equation}
Above, we have the coupling coefficients $\gamma_{mn}^{u}=e^{2\pi i\frac{m+n}{J}}-e^{-2\pi i\frac{2m-n}{J}}$,
$\gamma_{pq}^{v}=cb\left(e^{-2\pi i\frac{p+q}{JL}}-e^{2\pi i\frac{2p-q}{JL}}\right)$,
$\lambda_{l}=\frac{hc}{b}\frac{1-e^{-2\pi i\frac{l}{J}}}{1-e^{-2\pi i\frac{l}{JL}}}$,
and the quasilinear operator for interaction with the mean states,
$d_{u}=1+\gamma_{k}^{u*}\bar{u}$, $d_{v}=c\left(1+b\gamma_{l}^{v*}\bar{v}\right)$ (see the detailed model parameters in \eqref{eq:cov_2layer}).
Notice that in the large-scale modes, the indices in the summation
for nonlinear coupling go through all the wavenumbers, $\left|m\right|\leq J/2$,
while in the small-scale modes the indices for nonlinear coupling include all the small-scale
wavenumbers, $\left|p\right|\leq JL/2$. This leads
to an extremely high computational overload considering the high dimension
of the small-scale modes. Finally, the large and small scale modes are coupled
through the last linear terms. Especially, the small-scale modes $Y_{k+sJ},\left|s\right|\leq L/2$
give a combined feedback to the large-scale mode $Z_{k}$, while each
$Z_{k}$ is acting on a sequence of small-scale modes $Y_{k+sJ}$.
It is realized that in \eqref{eq:fluc_2layer} the most computational
demanding part comes from the summation terms going through all the wavenumbers
representing both linear and nonlinear coupling between scales. 

Next, we present the performance of the RBM models on the two-layer
L-96 with genuinely high dimensional and multscale scale processes.
In particular, we first consider the full RBM model in Algorithm \ref{alg:full-RBM},
then further reduce the computational cost by applying the reduced-order
RBM model in Algorithm \ref{alg:red-RBM} exploiting the large number
of fast-mixing small-scale modes.

\subsubsection{Numerical results for the full RBM model}

First, we consider the full RBM approximation, that is, to introduce
random batch decomposition in both the large and small scale stochastic
modes $\left\{ Z_{k}^{\left(i\right)},Y_{l}^{\left(i\right)}\right\} $,
while still run independent samples for ensemble simulation of
the entire stochastic equations during each time updating interval. This
leads to the \emph{full RBM model} for the stochastic coefficients
of the two-layer L-96 system following the general formulation \eqref{rbm_model}
\begin{equation}
\begin{aligned}\frac{\mathrm{d}Z_{k}^{\left(i\right)}}{\mathrm{d}t} & =c_{p}\sum_{m\in\mathcal{I}_{k}^{s}}\left(Z_{m}^{\left(i\right)}Z_{m-k}^{\left(i\right)*}-r_{m}^{u}\delta_{m,m-k}\right)\gamma_{m,m-k}^{u}-d_{u}Z_{k}^{\left(i\right)}-c_{L}\sum_{k+sJ\in\mathcal{J}_{l}^{s}}\lambda_{k+sJ}^{*}Y_{k+sJ}^{\left(i\right)},\\
\frac{\mathrm{d}Y_{l}^{\left(i\right)}}{\mathrm{d}t} & =c_{q}\sum_{p\in\mathcal{J}_{l}^{s}}\left(Y_{p}^{\left(i\right)}Y_{p-l}^{\left(i\right)*}-r_{p}^{v}\delta_{p,p-l}\right)\gamma_{p,p-l}^{v}-d_{v}Y_{l}^{\left(i\right)}+\lambda_{l}Z_{\mathrm{mod}\left(l,J\right)}^{\left(i\right)},
\end{aligned}
\label{eq:fluc_2layer-rbm1}
\end{equation}
with independent samples $i=1,\cdots,M_{1}$. In the RBM approximation
at each time updating step $t=t_{s}$, we choose the batches $\cup_{k}\mathcal{I}_{k}^{s}=\left\{ k:\left|k\right|\leq\frac{J}{2}\right\} $
with batch size $p=\left|\mathcal{I}_{k}^{s}\right|$ for the large-scale
modes $Z_{k}$, and batches $\cup_{l}\mathcal{J}_{l}^{s}=\left\{ l:\left|l\right|\leq\frac{JL}{2}\right\} $
with $q=\left|\mathcal{J}_{l}^{s}\right|$ for the small-scale modes
$Y_{l}$. Then the original summation terms taking over the entire
spectrum space are reduced to the summation for modes restricted inside
one batch of very small size $\left(p,q\right)$ (with the pair elements corresponding to large and small scale batch sizes). The new normalization
factors can be found as $c_{p}=\frac{1}{J}\frac{J-1}{p-1}$, $c_{q}=\frac{1}{JL}\frac{JL-1}{q-1}$,
and $c_{L}=\frac{1}{L}\frac{L-1}{q-1}$ following the same principle
as the one-layer case. The corresponding RBM equations for the associated
covariance equations can be derived accordingly. We list the explicit
equations in Appendix \ref{subsec:RBM-equations-2layer}. In this
way, the computational cost for the nonlinear coupling terms, particularly
the small scales with a large number of modes, are greatly saved by
constraining the high-order interactions only inside the very small
batch. This leads to the effective reduction in computational
cost to $O\left(M_{1}J\left(p+Lq\right)\right)$ in solving the stochastic
equations compared to the original cost $O\left(MJ^{2}\left(1+L^{2}\right)\right)$
in the direct MC approach in computing the full ensemble of size $M\gg M_{1}$.

In the numerical test of the full RBM model, we take the model parameters
for the standard regime $c=10,b=10,h=1,F=20$ and $J=8,L=32$. The
resulting coupled large and small-scale processes form a very high
total dimension $d=J\left(L+1\right)=264$. In particular, the small-scale
variable $v_{i}$ has a much faster time scale and the large scale
state $u_{j}$. This leads to a more challenging multiscale problem
since we have to resolve the large number of small scale modes even
that we are only interested in the large-scale state. In order to
get accurate statistical prediction resolving all the multiscale features,
we use a very large ensemble size $M=5\times10^{5}$ for the direct
MC simulation to accurate recover the true reference solution. In
the RBM model, three different batch pairs $\left(p,q\right)=\left(4,16\right),\left(4,8\right),\left(2,4\right)$
are tested in contrast to the total number of modes $\left(J,JL\right)=\left(8,256\right)$.
A much smaller ensemble size $M_{1}=500$ is used and this is made
possible by sampling only the $q$-dimensional subspace instead of
the $JL$-dimensional full space of the small scales. 

\begin{figure}
\centering
\subfloat{\includegraphics[scale=0.4]{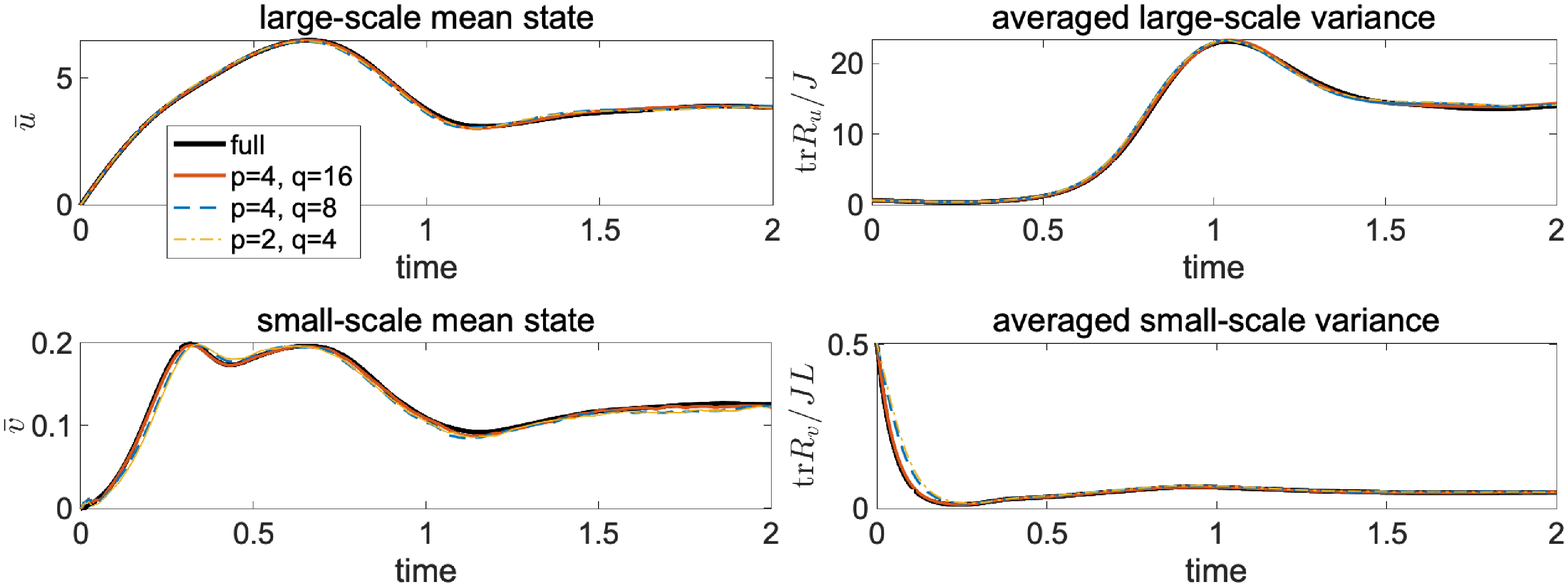}

}

\subfloat{\includegraphics[scale=0.4]{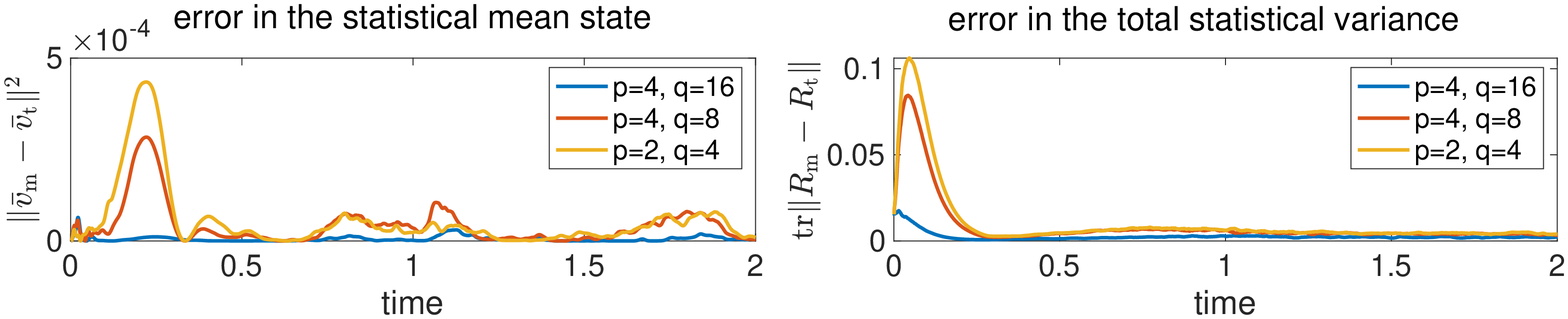}}

\caption{Prediction of the mean and averaged variance in both large and small
scale modes from full RBM model with different batch sizes $\left(p,q\right)$
compared with the full MC model with $JL=256$ modes in small scales and sample size $M=5\times10^5$.
$M_{1}=500$ samples are used for all the RBM models. The trajectory
errors in the mean and variance are also compared.\label{fig:tseries_2layer}}
\end{figure}

\begin{figure}
\centering
\includegraphics[scale=0.4]{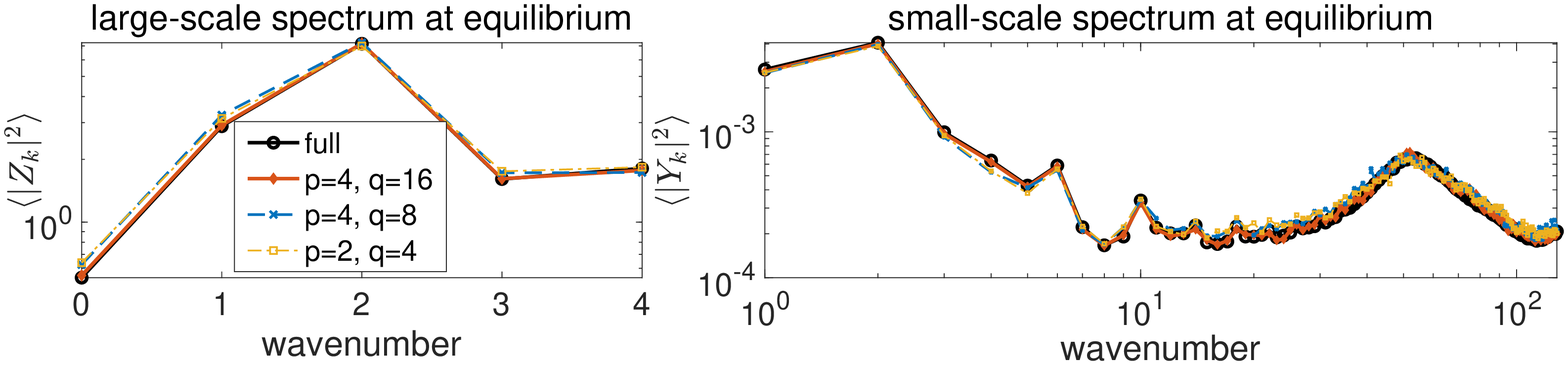}

\caption{Prediction of the variance spectra for large and small scale modes
in full RBM model with different batch sizes $\left(p,q\right)$ compared
with the full MC model with $JL=256$ modes in small scales  and sample size $M=5\times10^5$.\label{fig:spect_2layer}}

\end{figure}

In Figure \ref{fig:tseries_2layer}, we first show the RBM model prediction
for the mean $\bar{u},\bar{v}$ and average total variance $\sum r_{k}^{u}/J,\sum r_{l}^{v}/JL$
with three different batch sizes, as well as the pointwise errors
in mean and total varianc at each time step. It shows accurate recovery
of the key statistics in both the mean and variance, and in both the
starting transient stage and final equilibrium. The high accuracy
is maintained with almost indistinguishable results as we reduce the
batch sizes $\left(p,q\right)$. Relatively larger errors are observed
in the starting transient stage when a small batch size is applied.
Notice that significant computational reduction is achieved where the smallest
batch size refers using only $q=4$ nonlinear coupling terms for each
$Y_{l}$ out of the total $JL=256$ terms. In addition, the detailed
model prediction for the variance spectra in each large and small
scale modes is shown in Figure \ref{fig:spect_2layer}. A wide range
of the small-scale modes with high wavenumbers are excited due to
the instability in small scales, as illustrate in Figure \ref{fig:Growth-rates-2layer}.
This presents a highly challenging scenario as we only resolve a very
small batch of modes for the nonlinear coupling term to stabilize
the large number of unstable modes. Again, the energy spectrum is
recovered exactly with uniformly high skill using the three extremely
small batch sizes. This highlights the robust skill in the RBM model
to achieve computational efficiency and maintain high prediction accuracy
at the same time.

\begin{figure}
\centering
\subfloat[large-scale modes]{\includegraphics[scale=0.32]{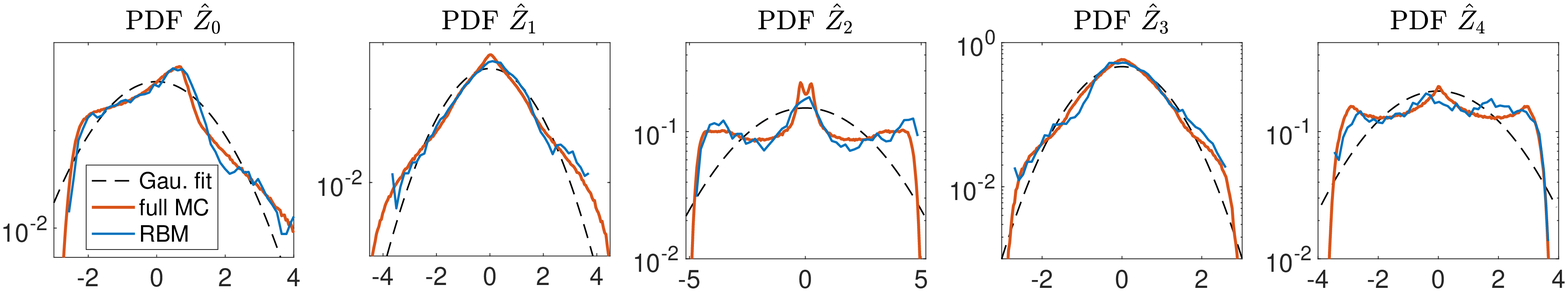}

}

\vspace{-1em}

\subfloat[small-scale modes]{\includegraphics[scale=0.32]{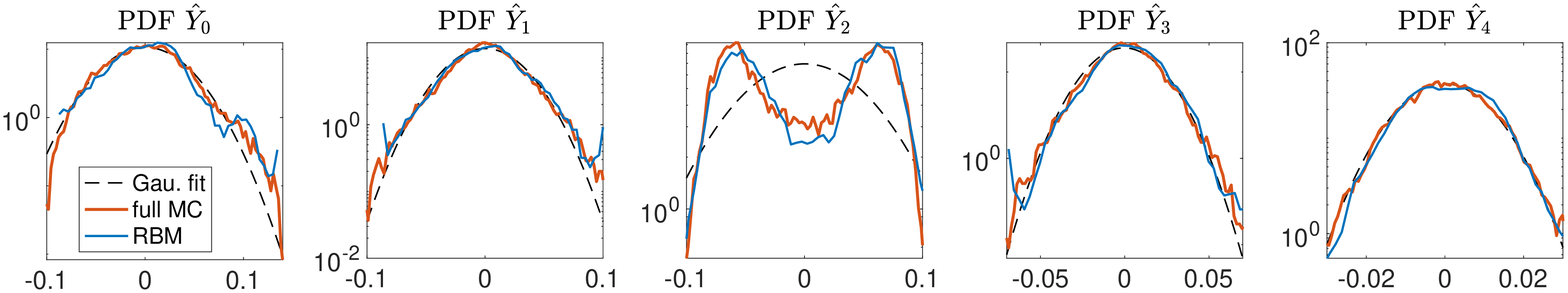}}

\caption{1D marginal PDFs of the leading large and small scale stochastic modes
$Z_{k},Y_{l}$ of the two-layer L-96 system. The RBM predictions are
compared with the truth from direct MC simulation. The Gaussian fits
of the PDFs with the same mean and variance are also plotted in dashed
lines.\label{fig:1D-marginal-PDFs-2layer}}

\end{figure}

\begin{figure}
\subfloat[MC large scale $Z_{k}$]{\includegraphics[scale=0.31]{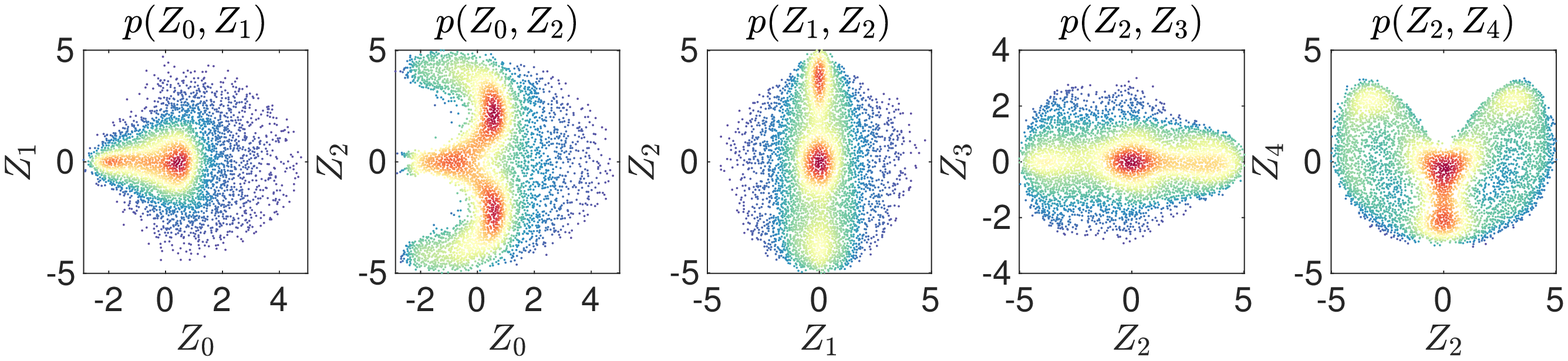}

}\hspace{.5em}\subfloat[MC small scale $Y_{l}$]{\includegraphics[scale=0.31]{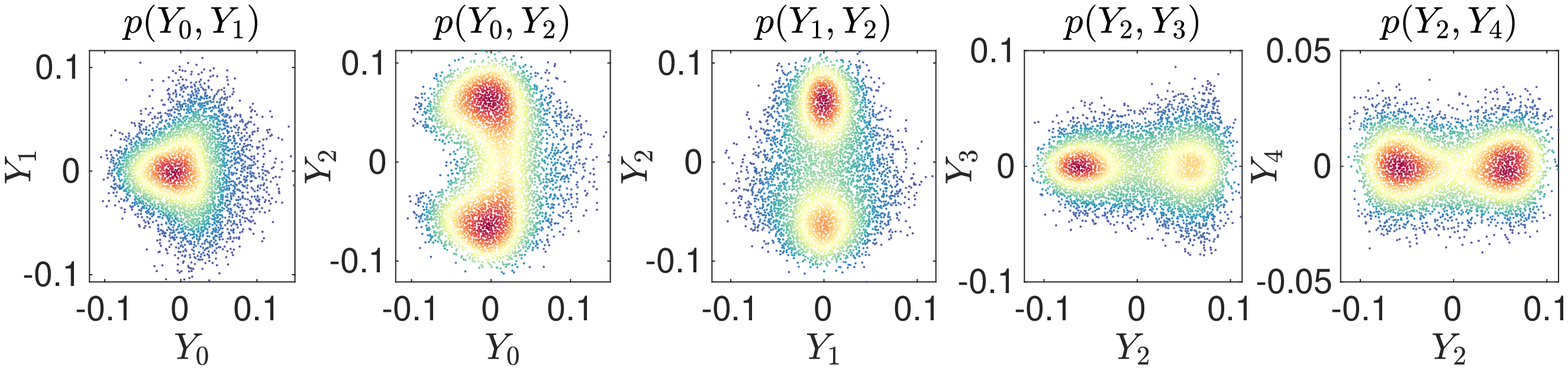}

}

\vspace{-1.em}

\subfloat[RBM large scale $Z_{k}$]{\includegraphics[scale=0.31]{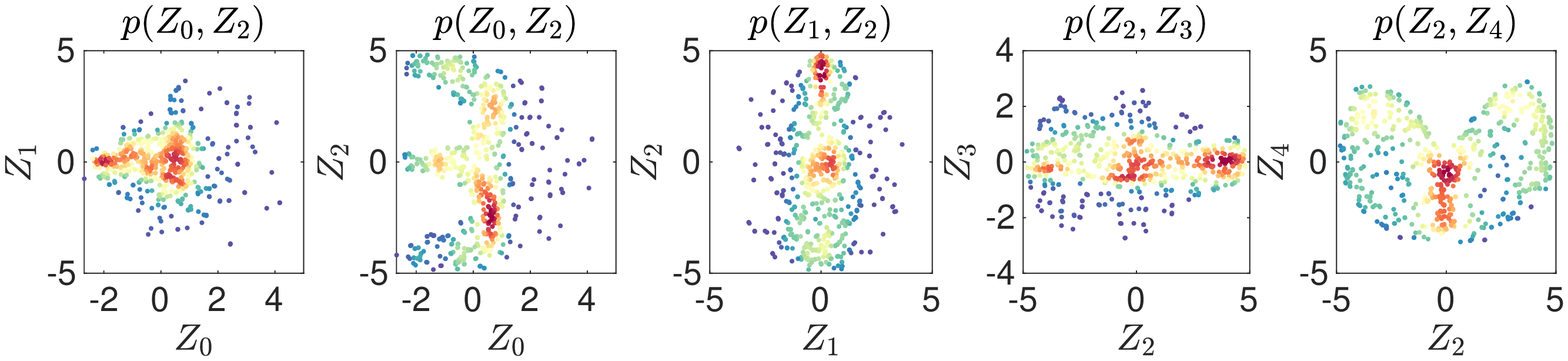}

}\hspace{.5em}\subfloat[RBM small scale $Y_{l}$]{\includegraphics[scale=0.31]{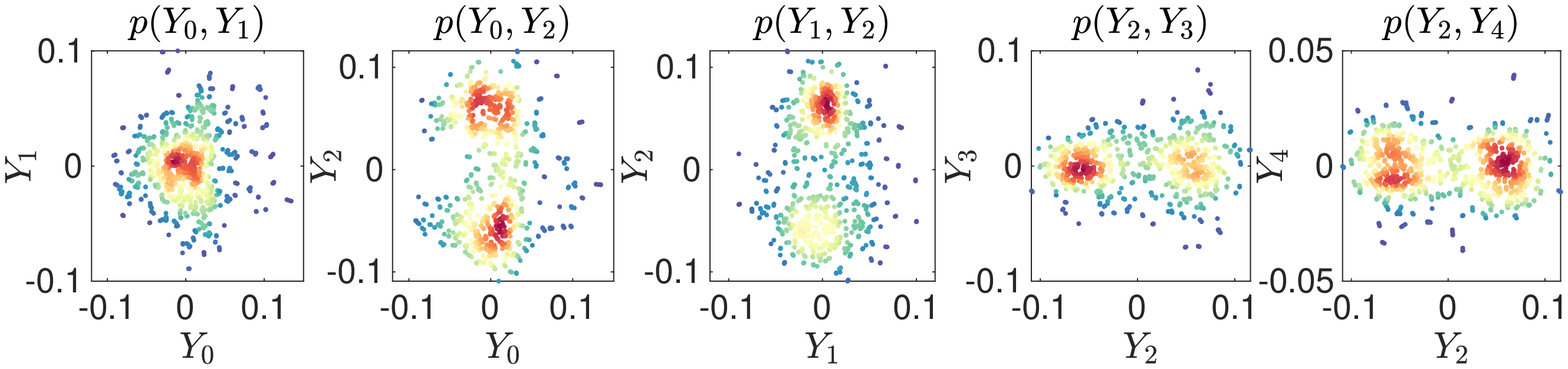}

}

\caption{2D joint PDFs of the leading large and small scale stochastic modes
$Z_{k},Y_{l}$ of the two-layer L-96 system. The RBM predictions are 
compared with the truth from direct MC simulation.\label{fig:2D-joint-PDFs-2layer}}
\end{figure}

Next, we proceed to check the ability of the RBM model to  capture the PDFs of
the dominant modes. Figure \ref{fig:1D-marginal-PDFs-2layer} shows
the marginal PDFs of the leading stochastic modes in both large and
small scales. In this case of the two-layer L-96 system, even stronger
non-Gaussian statistics are displayed with highly skewed and fat-tailed
PDFs. Accurate characterization of such non-Gaussian PDFs becomes even more critical
issue in statistical prediction of such high-dimensional systems. Usually,
the small number of samples tend to concentrate near the central part
of the PDF and miss the crucial extreme events featuring the edge regions of the large phase space. Using the RBM model
with a very small ensemble size, the strongly non-Gaussian PDFs are
successfully captured confirming the very high skill of the RBM approach
in recovering the true model statistics not only in the leading moments
but also in the more challenging higher-order statistics. In addition,
we also compare the joint PDFs in the predicted leading modes in Figure
\ref{fig:2D-joint-PDFs-2layer}. The non-Gaussian structures become
more pronounced in the 2D distributions, revealing bimodal PDFs as already
indicated in the marginal PDFs. The RBM predictions maintain
very accurate in recover the highly non-Gaussian statistics despite using
a very small ensemble size.

\begin{figure}
\centering
\includegraphics[scale=0.35]{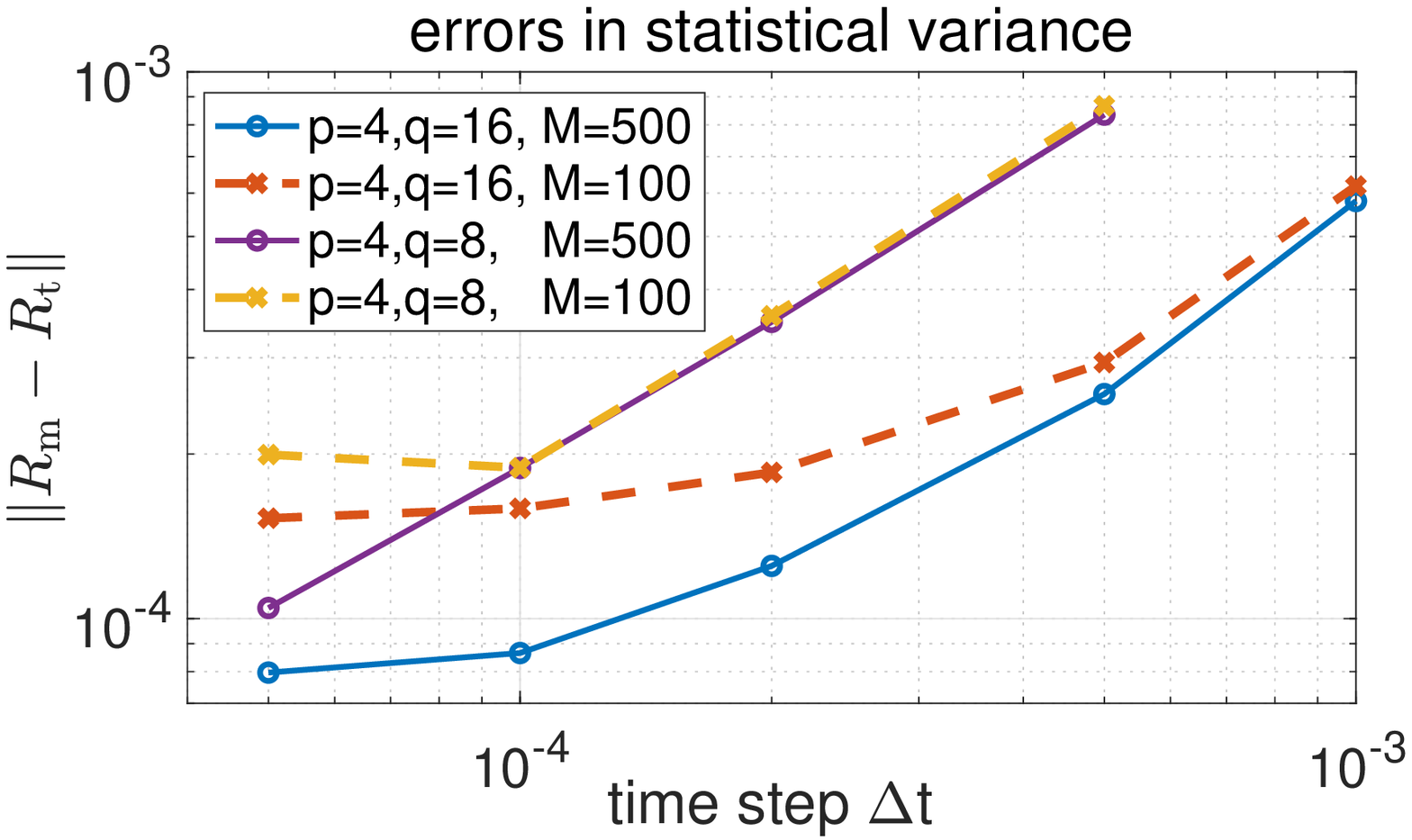}\includegraphics[scale=0.35]{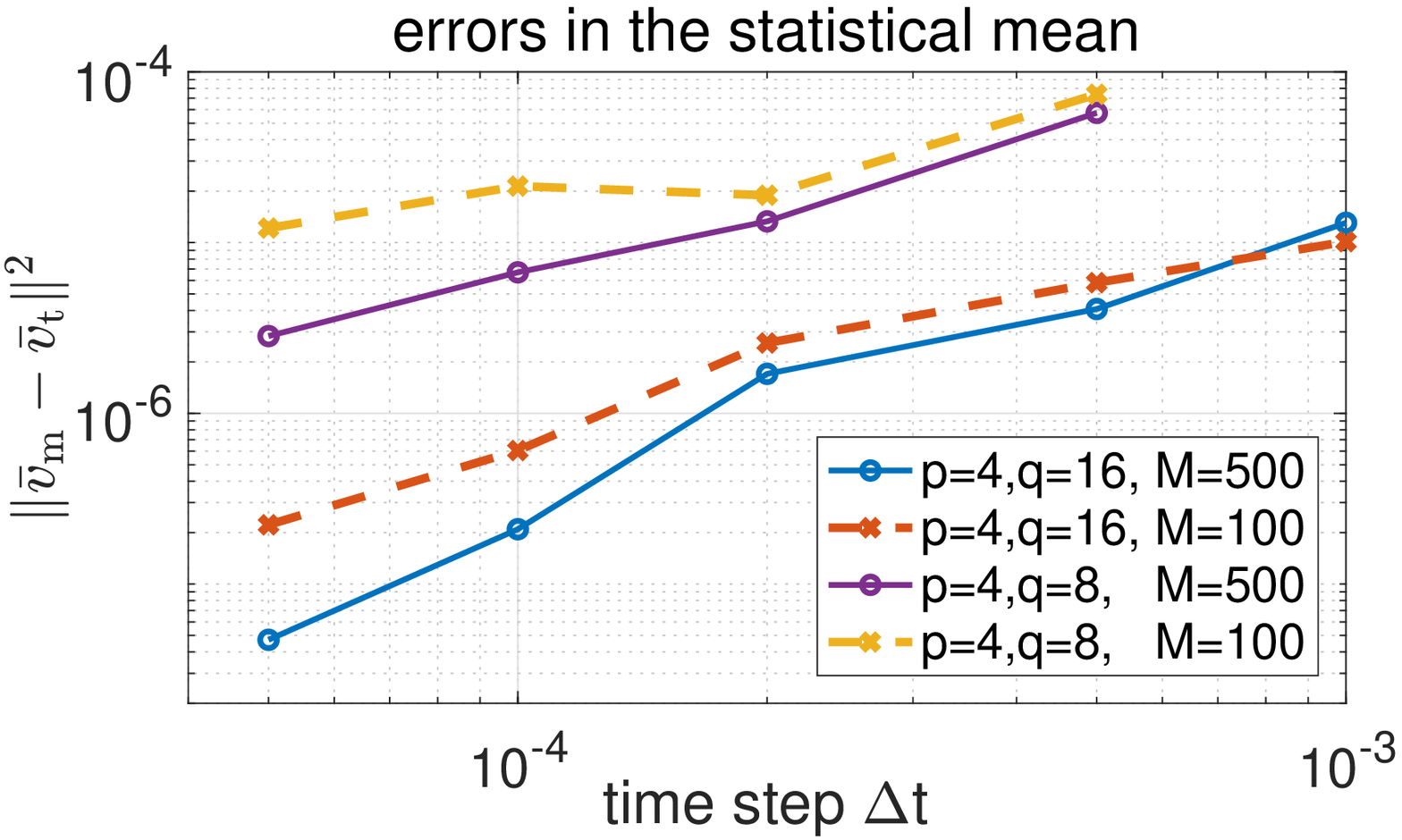}\caption{Errors in the statistical mean and total variance from the RBM prediction
with different batch and ensemble sizes for the two-layer L-96 system.\label{fig:Errors-2layer}}

\end{figure}

Finally, we provide a quantitative comparison of the RBM model prediction
errors as the key model parameters change. In Figure \ref{fig:Errors-2layer},
we show the averaged errors in the total variance $\left\Vert R_{\mathrm{m}}-R_{\mathrm{t}}\right\Vert =\frac{1}{N}\sum_{s=1}^{N}\sum_{k}\left|r_{\mathrm{m},k}\left(t_{s}\right)-r_{\mathrm{t},k}\left(t_{s}\right)\right|$
and the mean $\left\Vert \bar{v}_{\mathrm{m}}-\bar{v}_{\mathrm{t}}\right\Vert ^{2}=\frac{1}{N}\sum_{s=1}^{N}\left|\bar{v}_{\mathrm{m}}\left(t_{s}\right)-\bar{v}_{\mathrm{t}}\left(t_{s}\right)\right|^{2}$
according to the time integration step $\Delta t=t_{s+1}-t_{s}$ with
two batch sizes $\left(p,q\right)=\left(4,16\right),\left(4,8\right)$
and two ensemble sizes $M_{1}=500,100$. The errors grow with larger
time step size and with smaller ensemble and batch sizes agreeing with the intuition.
Among the range of large to moderate time step sizes, reducing the time step can
effectively improve the prediction accuracy since it is corresponding
to more frequent resampling (thus more equivalent samples in the time average) in the RBM approximation. With very small
time step size though, the error will saturate because the small ensemble
size will account for the major error in the statistical estimates.
It also implies that it might be useful to use a smaller batch size
$q$ with a relatively larger ensemble to achieve accurate prediction
with minimum computational cost.

\subsubsection{Numerical results for the reduced-order RBM model}

Next, we consider a further reduction of computation cost by applying
the reduced-order RBM model as described in Algorithm \ref{alg:red-RBM}.
In the reduced RBM model, we still keep the statistical mean and covariance
equations \eqref{eq:mean_2layer} and \eqref{eq:cov_2layer_rbm} the
same as the full RBM model case. In the two-layer L-96 system, most of the computational
cost comes from the ensemble simulation for the stochastic coefficients
in \eqref{eq:fluc_2layer} for the wide range of small-scale
modes $Y_{l}$. In most situations, we are mostly interested in the
statistics in the main large-scale modes $Z_{k}$ while the  small-scale modes give crucial combined feedback to affect the large-scale
state, thus their contributions cannot be simply ignored. The unstable growth (shown in Figure
\ref{fig:Growth-rates-2layer}) and excited small scales in equilibrium
energy spectrum (in Figure \ref{fig:spect_2layer}) emphasize the
non-negligible role of these small-scale modes. This sets the inherent
obstacle in developing effective reduced-order models for multiscale
systems.

To enable further computational reduction for the small-scale modes,
following the idea in \eqref{eq:rbm_stoch2} we group the ensemble
members in the large scale modes $\mathbf{Z}^{\left(i\right)},i=1,\cdots,M_{1}$
with one single small-scale state $\mathbf{Y}$. In applying the RBM
approximation, the same partition $\mathcal{I}_{k}^{s}$ is applied to the small number of large-scale modes. The model reduction strategy is applied
to the expensive small-scale equations. The large number of small-scale
modes are then divided into the batches $\mathcal{J}_{i}^{s}$, satisfying $\cup_{i=1}^{M_{1}}\mathcal{J}_{i}^{s}=\left\{ l:\left|l\right|\leq JL/2\right\} $.
According to the modes in the batches, the full system is decomposed
into much smaller subsystems with one large-scale sample $\mathbf{Z}^{\left(i\right)}$
together with a portion of the small-scale modes $Y_{l},l\in\mathcal{J}_{i}^{s}$.
The model reduction is made possible by the very large number of fast-mixing
small-scale modes. This leads to the coupled equations for the stochastic
coefficients in subsets $\left\{ \mathbf{Z}^{\left(i\right)},Y_{l}\right\} _{l\in\mathcal{J}_{i}^{s}}$
during the time interval $t\in\left(t_{s},t_{s+1}\right]$
\begin{equation}
\begin{aligned}\frac{\mathrm{d}Z_{k}^{\left(i\right)}}{\mathrm{d}t}= & c_{p}\sum_{m\in\mathcal{I}_{k}^{s}}\left(Z_{m}^{\left(i\right)}Z_{m-k}^{\left(i\right)*}-r_{m}^{u}\delta_{m,m-k}\right)\gamma_{m,m-k}^{u}-d_{u}Z_{k}^{\left(i\right)}-c_{L}\sum_{k+sJ\in\mathcal{J}_{i}^{s}}\tilde{\lambda}_{k+sJ}^{*}Y_{k+sJ},\\
\frac{\mathrm{d}Y_{l}}{\mathrm{d}t}= & c_{q}\sum_{p\in\mathcal{J}_{i}^{n}}\left(Y_{p}Y_{p-l}^{*}-r_{p}^{v}\delta_{p,p-l}\right)\gamma_{p,p-l}^{v}-d_{v}Y_{l}+\tilde{\lambda}_{l}Z_{\mathrm{mod}\left(l,J\right)}^{\left(i\right)},\qquad l\in\mathcal{J}_{i}^{n}.
\end{aligned}
\label{eq:fluc_2layer_rbm2}
\end{equation}
Above, the ensemble is used only to sample the low-dimensional large-scale
state $\mathbf{Z}^{\left(i\right)}\in\mathbb{R}^{J}$, while only
a small portion of small-scale modes $\left\{ Y_{l}\right\} _{l\in\mathcal{I}_{i}}$
are grouped with the $i$-th sample $\mathbf{Z}^{\left(i\right)}$
together for the time update at $t=t_{s}$. The union of all the randomly
sampled groups form the entire spectrum of small-scale modes. In this
way, we no longer need to run a very large ensemble for the small
scales $Y_{l}$ by exploiting its wide spectrum of modes acting at
different large-scale samples. Notice that the samples $\mathbf{Z}^{\left(i\right)}$
will no longer stay independent and will be linked by the small-scale
modes through the random batches. Still, the important correlations
between the small and large scale modes are maintained through this
splitting of small-scale modes. 
\begin{rem*}
In practice, we may still want to run a small number of small-scale
modes $\mathbf{Y}^{\left(j\right)},j=1,\cdots,M_{2}$. This is equivalent
to introduce an ensemble of size $M_{2}$ to the above block model
\eqref{eq:fluc_2layer_rbm2}, so that the model parameters have the
relation $\kappa=\frac{M_{1}}{M_{2}}=\frac{JL}{2q}$. For consistent
statistics, the large and small scale coupling coefficient needs to
be updated as $\tilde{\lambda}_{l}=\sqrt{\kappa}\lambda_{l}$. This
leads to the computational reduction from $O\left(M_{1}J\left(p+Lq\right)\right)$
in the full RBM model to $O\left(M_{1}Jp+M_{2}JLq\right)=O\left(M_{1}\left(q^{2}+pJ\right)\right)$
which is only dependent on the dimension $J$ of the large-scale state
and independent of the small scale dimension $JL$.
\end{rem*}
In the numerical test of the reduced-order model, we consider the
parameter regime of the two-layer L-96 system with $c=4,b=10,h=1,F=20$.
This regime gives a weaker scale separation between the large and
small scales thus sets a more important role for the small-scale processes
with a non-negligible contribution to the large-scale modes. We fix
the batch size as $\left(p,q\right)=\left(4,32\right)$. The relatively
large batch size $q=32$ (compared with the full dimension $JL=256$)
is used to ensure a larger allowed sample size $M_{2}=\frac{2qM_{1}}{JL}$
for the small-scale modes. We pick a moderate sample size $M_{1}=500$
to sample the large-scale mode of dimension $J$ and the large number
of small-scale modes only requires a small sample size $M_{2}=125$.
Another extreme sample size $M_{1}=100,M_{2}=25$ is also tested.
The reduced-order model enables even more efficient computation compared
with the very high dimension of the full system $d=J\left(L+1\right)=264$.

\begin{figure}
\includegraphics[scale=0.3]{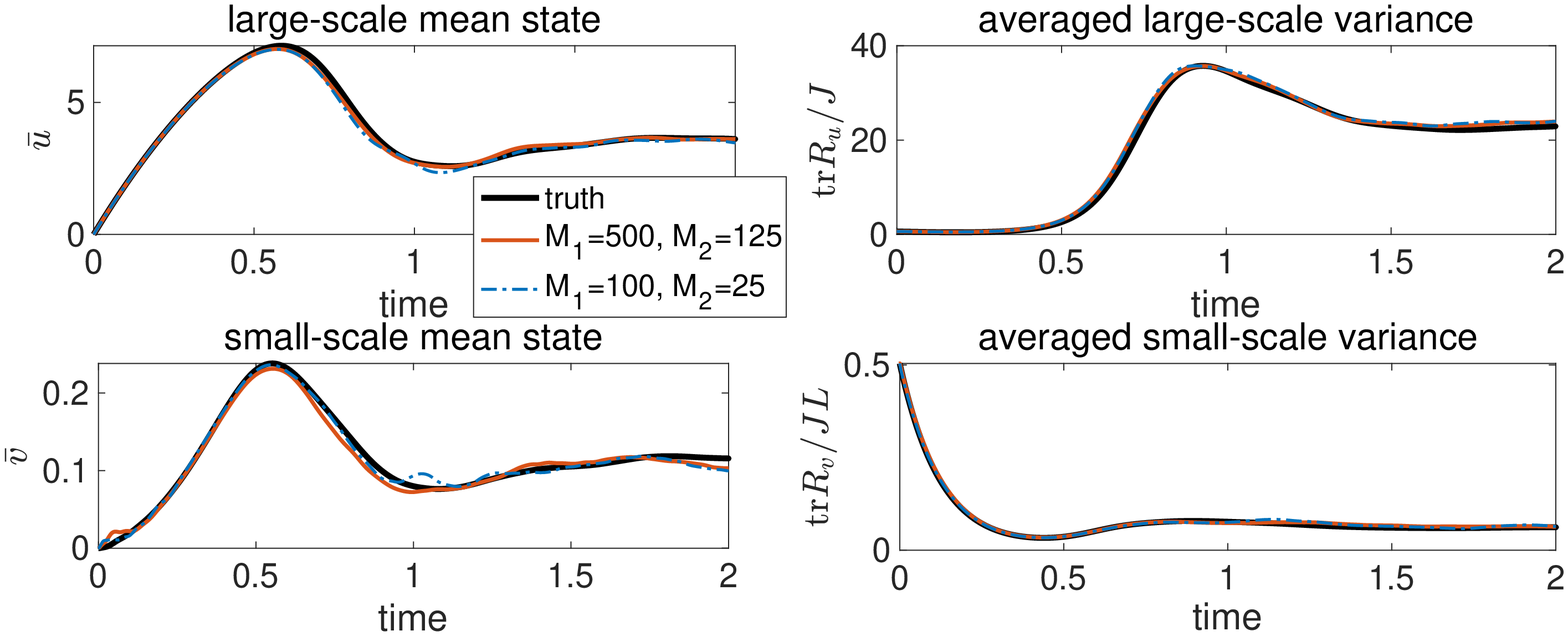}\includegraphics[scale=0.3]{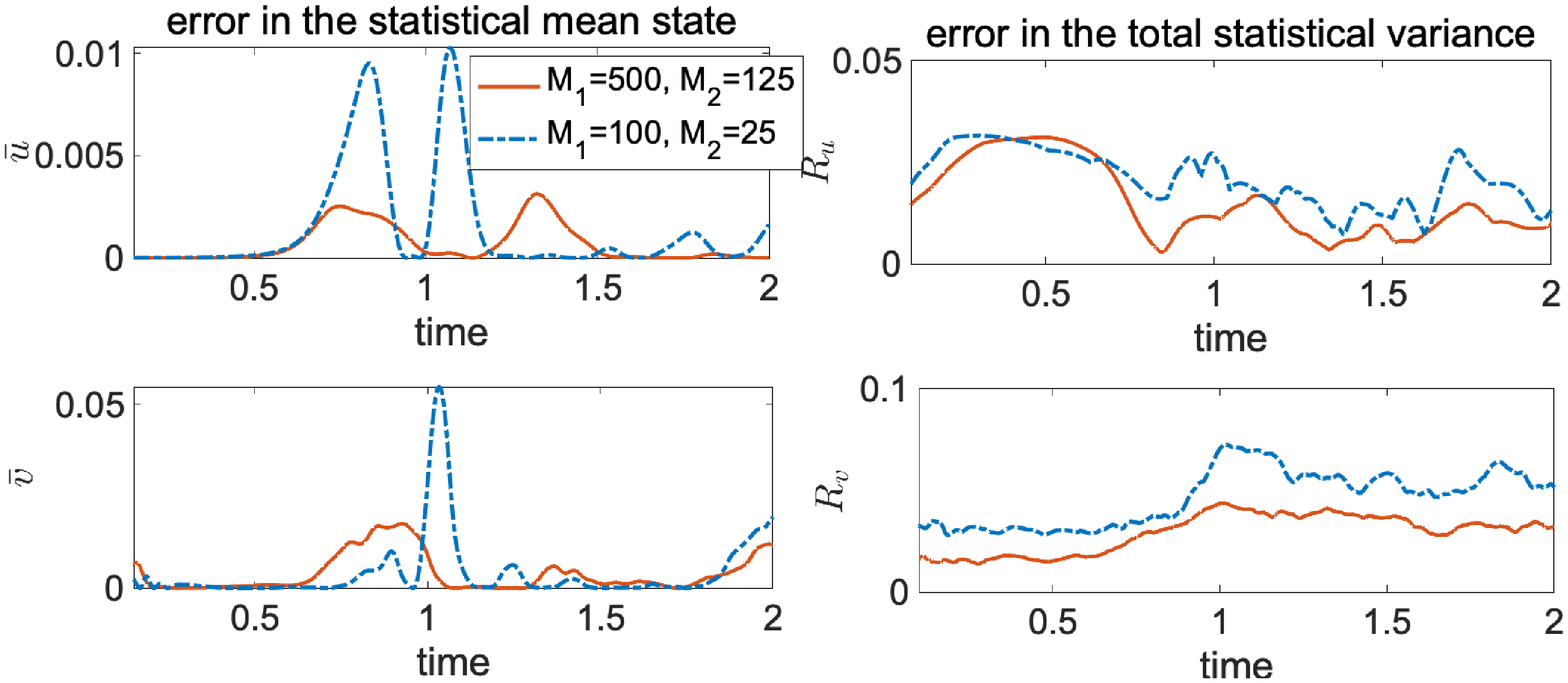}

\caption{Prediction of the mean and variances in the reduced-order RBM model
with batch size $\left(p,q\right)=\left(4,32\right)$ and sample size $(M_1,M_2)=(500,125), (100,25)$ with for $M_1$ samples for large-scale state and $M_2$ for small-scale state. The time evolutions
of the errors in mean and variance are also compared. \label{fig:Prediction-red-2layer}}

\end{figure}

\begin{figure}
\centering
\subfloat{\includegraphics[scale=0.32]{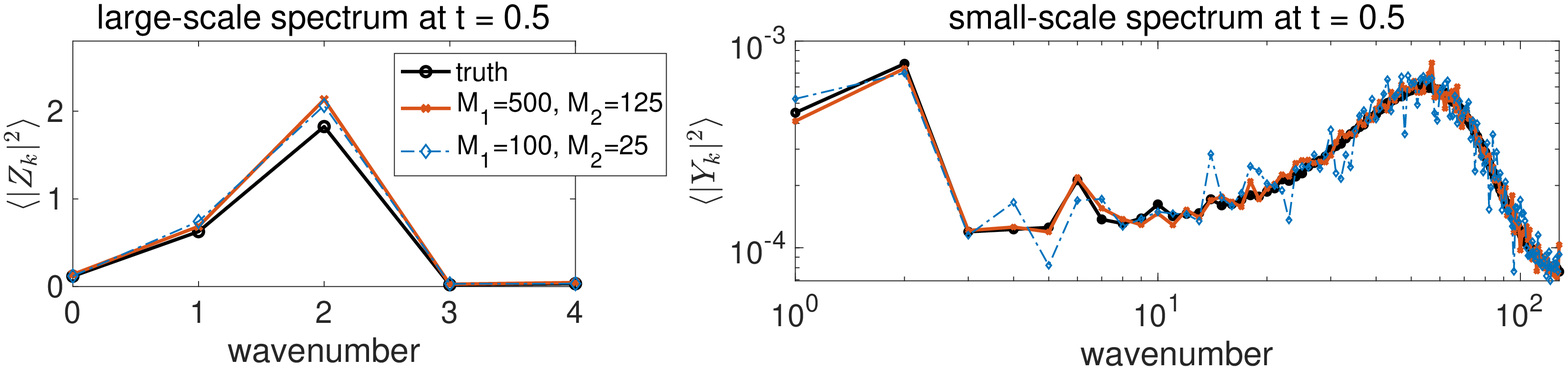}}

\vspace{-1em}

\subfloat{\includegraphics[scale=0.32]{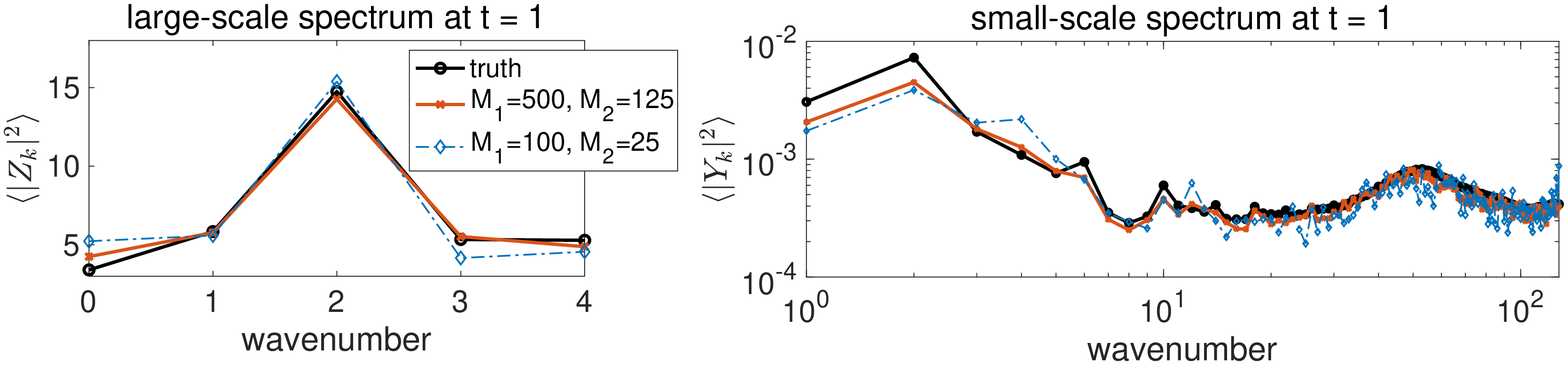}}

\vspace{-1em}

\subfloat{\includegraphics[scale=0.32]{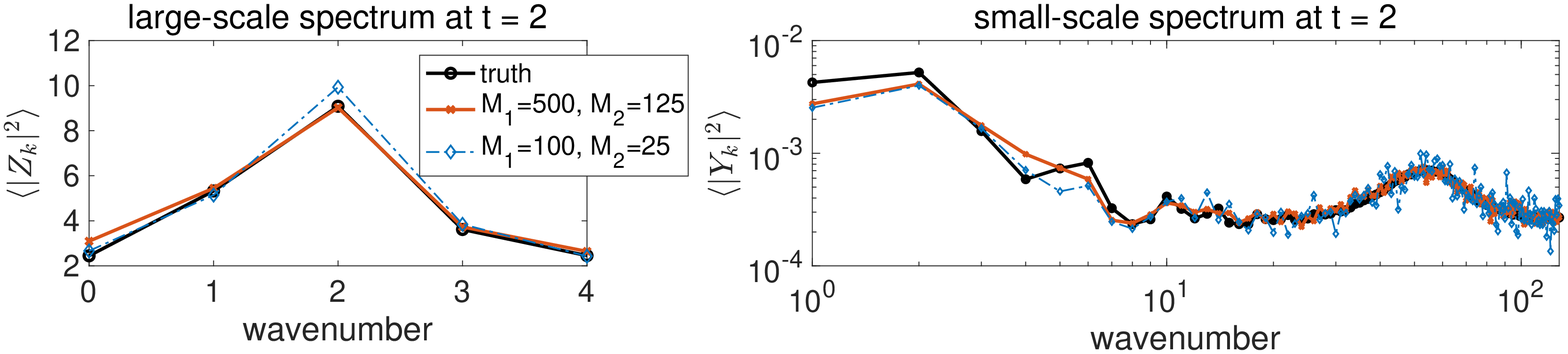}}

\caption{Reduced-order RBM model prediction of the variance spectra in both
large and small scale modes at several time instants $t=0.5,1,2$.\label{fig:Reduced-order-RBM-model}}

\end{figure}

We show the prediction results in the reduced-order RBM model in Figure
\ref{fig:Prediction-red-2layer}. It can be seen that the reduced
model maintains the high skill to capture the leading statistics in
both the large and small scale states while save additional computational
cost by avoiding running the large ensemble simulation. The pointwise
errors in mean and variance further confirm the accurate statistical
prediction even with a small and reduced sample size. Accordingly,
we show the detailed prediction of the energy spectra at several time
instants in both large and small scale modes in Figure \ref{fig:Reduced-order-RBM-model}.
In the case of reduced-order model due to the very small ensemble
size, it is expected the small-scale modes will become very noisy
and have larger fluctuation errors especially for the large wavenumber
modes. Still, the reduced-order RBM model successfully captures the
structure of spectra in both large and small scales. These encouraging results suggest further applications of the RBM models to more practical and realistic systems
exhibiting stronger multiscale coupling and a very wide spectrum covering
a large number of scales.


\section{Summary\label{sec:Summary}}

We present a systematic closure modeling framework that enables efficient
ensemble prediction of leading-order statistics and non-Gaussian PDFs
in complex turbulent systems, which are characterized by strong internal
instability and interactions of coupled spatio-temporal scales. A
general stochastic-statistical formulation is established derived
from a generic multiscale nonlinear system, capable of modeling a
wide range of complex phenomena observed in natural and engineering
systems. A mean and fluctuation decomposition of the original model
state is introduced to deal with the irreducible dynamics in the high-dimensional
equations. The fluctuation modes, which capture large uncertainty
in multiscale modes, are modeled using stochastic equations that depend
on the mean state and covariance of fluctuation modes. The leading-order
statistical equations for the mean and covariance incorporate higher-order
moment feedback from different scales. A precise high-order closure
is introduced using empirical average of the ensemble prediction of
the stochastic equation solution, which explicitly captures detailed
multiscale interactions. The coupling between stochastic and statistical
equations is further reinforced by a relaxation term for statistical
consistency. This approach effectively balances the strong internal
instability often encountered in turbulent systems, mitigating the
inherent problem of numerical divergence in the statistical solution.
The model achieves high skill in capturing non-Gaussian statistics
and extreme outliers by explicitly resolving crucial higher-order
moments instead of relying on insufficient low-order parameterizations.
The stochastic-statistical formulation offers a flexible approach
for recovering essential model statistics, making it applicable to
a wide range of problems in uncertainty quantification and data assimilation
\cite{stuart2010inverse,reich2015probabilistic,calvello2022ensemble}.

To enable efficient computation using a small ensemble size for the
stochastic equation and limit the number of high-order nonlinear coupling
terms during each time update, we generalize the idea in RBM approximation
for mean-fluctuation equations \cite{qi2023random}. The approach
decomposes the wide spectrum involving large number of multiscale
modes into small random batches. A much smaller ensemble size is made
possible by just sampling the nonlinear coupling terms involving modes
in a low-dimensional subspace inside one batch. Simutaneously, the
contribution from all the other modes is fully modeled by resampling
the batches at each time step, exploiting the ergodicity of the stochastic
modes. Consequently, high prediction accuracy concerning all high-order
feedbacks is achieved. For systems with exceptionally high dimensionality,
a model reduction strategy is proposed to further reduce the computational
cost by linking the large number of small-scale fluctuation modes
to the ensemble samples of large-scale state. The resulting algorithms
are straightforward to implement and are well-suited for a wide range
of multiscale turbulent systems. We evaluate the performance of the
proposed RBM models using the representative one-layer and two-layer
L-96 systems, which exhibit strong multiscale coupling and a wide
spectrum of energetic unstable modes. The models demonstrate uniform
high prediction skill for the leading order mean and variance, and
more notably, capture the highly non-Gaussian PDFs and extreme events
using a very small ensemble. As a result, the computational cost is
reduced to an affordable level for genuinely high dimensional systems.
In future research directions, we aim to develop a complete theory
to analyze the general stochastic-statistical modeling framework building
upon the preliminary estimates presented in this paper. The promising
results also suggest potential applications to more realistic high-dimensional
systems. The reduced-order RBM model shows potential in overcoming
the curse of dimensionality and providing an effective tool for a
wide range of practical problems related to prediction and data assimilation.

\section*{Acknowledgements}

The research of D.Q. was partially supported by the start-up funds
and the PCCRC Seed Funding provided by Purdue University. The research
of J.-G. L. is partially supported under the NSF grant No. DMS-2106988.

\section*{Data Availability}

The data that support the findings of this study are available from the corresponding author upon reasonable request.

\appendix
\renewcommand\theequation{A\arabic{equation}}
\setcounter{equation}{0}

\section{Proof of theorems in Section \ref{subsec:Error-analysis}\label{sec:Proof-of-theorems}}
\begin{proof}
[Proof of Theorem \ref{thm:converg_coeff}]Define the following functions
according to solutions to \eqref{eq:coeff_rbm} and \eqref{eq:coeff_full}
\begin{equation}
\begin{aligned}\tilde{w}\left(\mathbf{x},t\right)= & \frac{1}{K}\sum_{k=1}^{K}\mathbb{E}_{\mathbf{x}}\varphi\left(\tilde{Z}_{k}\left(t\right)\right),\\
w\left(\mathbf{x},t\right)= & \frac{1}{K}\sum_{k=1}^{K}\mathbb{E}_{\mathbf{x}}\varphi\left(Z_{k}\left(t\right)\right),
\end{aligned}
\label{eq:test_funcs}
\end{equation}
with the test function $\varphi\in C_{b}^{2}$ and initial state $\mathbf{Z}\left(0\right)=\tilde{\mathbf{Z}}\left(0\right)=\mathbf{x}\in\mathbb{R}^{K}$.
The the functions \eqref{eq:test_funcs} defined in the
time interval $t\in\left(t_{s},t_{s+1}\right]$ are governed by the backward
Kolmogorov equation \cite{varadhan2007stochastic}  as
\begin{equation}
\begin{aligned}\partial_{t}\tilde{w}=\tilde{\mathcal{L}}_{\mathcal{I}^{s}}\tilde{w} & =\sum_{k=1}^{K}\left[\sum_{m,n}I_{k}^{s}\left(m\right)\tilde{L}_{v,km}\left(\bar{u}\right)x_{m}+I_{k}^{s}\left(m\right)I_{k}^{s}\left(n\right)\tilde{\gamma}_{mnk}\left(x_{m}x_{n}-R_{mn}\right)\right]\cdot\partial_{x_{k}}\tilde{w}+\frac{1}{2}\sum_{k=1}^{K}\sigma_{k}^{2}\partial_{x_{k}}^{2}\tilde{w},\\
\partial_{t}w=\mathcal{L}w & =\sum_{k=1}^{K}\left[\sum_{m,n}L_{v,km}\left(\bar{u}\right)x_{m}+\gamma_{mnk}\left(x_{m}x_{n}-R_{mn}\right)\right]\cdot\partial_{x_{k}}w+\frac{1}{2}\sum_{k=1}^{K}\sigma_{k}^{2}\partial_{x_{k}}^{2}w.
\end{aligned}
\label{eq:backward}
\end{equation}
Above in \eqref{eq:backward}, the backward equation for the RBM model
$\tilde{w}$ is subject to the additional randomness due to the partition
$\mathcal{I}^{s}=\left\{ \mathcal{I}_{k}^{s}\right\} $ during the
time updating interval. We introduce the index function defined in
time $t\in\left(t_{s},t_{s+1}\right]$
\[
I_{k}^{s}\left(m\right)=\begin{cases}
1, & \mathrm{if}\:m\in\mathcal{I}_{k}^{s},\\
0, & \mathrm{otherwise.}
\end{cases}
\]
Notice that the functions $I_{k}^{s}$ is kept constant during each
time interval $\left(t_{s},t_{s+1}\right]$ and will change values
subject to the resampling of the random batches. According to the
conclusion in \cite{qi2023random}, we have the expectation of the
partition functions by counting the ordered combinations of the random
batches
\[
\mathbb{E}I_{k}^{s}\left(m\right)=\frac{p}{K},\quad\mathbb{E}I_{k}^{s}\left(m\right)I_{k}^{s}\left(n\right)=\frac{p}{K}\frac{p-1}{K-1}.
\]
This leads to the consistent condition on the expectation with the
partition using the proper choice of the coupling coefficients
\[
\mathbb{E}^{\mathcal{I}^{s}}\tilde{\mathcal{L}}_{\mathcal{I}^{s}}=\mathcal{L}.
\]
Using the semigroup operator $\tilde{\mathcal{S}}$ acting on the
function $w\left(\mathbf{x},t_{s}\right)$ and the above identity,
we have
\begin{align*}
\tilde{\mathcal{S}}w\left(\mathbf{x},t_{s}\right)-w\left(\mathbf{x},t_{s+1}\right) & =\mathbb{E}^{\mathcal{I}^{s}}e^{\Delta t\tilde{\mathcal{L}}_{\mathcal{I}^{s}}}w\left(\mathbf{x},t_{s}\right)-w\left(\mathbf{x},t_{s+1}\right)\\
 & =\int_{0}^{\Delta t}\left(\Delta t-\tau\right)\left[\mathbb{E}^{\mathcal{I}^{s}}\left(\tilde{\mathcal{L}}_{\mathcal{I}^{s}}\right)^{2}e^{\tau\tilde{\mathcal{L}}_{\mathcal{I}^{s}}}-\mathcal{L}^{2}e^{\tau\mathcal{L}}\right]w\left(\mathbf{x},t_{s}\right)\mathrm{d}\tau.
\end{align*}
By the assumptions \eqref{eq:assump1}, the residual terms on the last equality are uniformly
bounded
\[
\left\Vert \left(\tilde{\mathcal{L}}_{\mathcal{I}^{s}}\right)^{2}e^{\tau\tilde{\mathcal{L}}_{\mathcal{I}^{s}}}w\left(\cdot,t\right)\right\Vert _{\infty}<C,\quad\left\Vert \mathcal{L}^{2}e^{\tau\mathcal{L}}w\left(\cdot,t\right)\right\Vert _{\infty}<C.
\]
Therefore, the one-step error for between the RBM solution $\tilde{\mathcal{S}}w\left(\mathbf{x},t_{s}\right)$
and the full model $w\left(\mathbf{x},t_{s+1}\right)$ can be estimated
as
\[
\left\Vert \tilde{\mathcal{S}}w\left(\cdot,t_{s}\right)-w\left(\cdot,t_{s+1}\right)\right\Vert _{\infty}\leq C\Delta t^{2}.
\]
Finally, by applying $\tilde{\mathcal{S}}$ on the initial function
$w\left(\mathbf{x},0\right)=\varphi\left(\mathbf{x}\right)$ $s$
times and recurrently using the above contraction property, we compute
the total error at $t=t_{s}$ as
\begin{align*}
\left\Vert \tilde{\mathcal{S}}^{\left(s\right)}\varphi\left(\cdot\right)-w\left(\cdot,t_{s}\right)\right\Vert _{\infty} & \leq\left\Vert \tilde{\mathcal{S}}\left[\tilde{\mathcal{S}}^{\left(s-1\right)}\varphi\left(\cdot\right)-w\left(\cdot,t_{s-1}\right)\right]\right\Vert _{\infty}+\left\Vert \tilde{\mathcal{S}}w\left(\cdot,t_{s-1}\right)-w\left(\cdot,t_{s}\right)\right\Vert _{\infty}\\
 & \leq\left\Vert \tilde{\mathcal{S}}^{\left(s-1\right)}\varphi\left(\cdot\right)-w\left(\cdot,t_{s-1}\right)\right\Vert _{\infty}+\left\Vert \tilde{\mathcal{S}}w\left(\cdot,t_{s-1}\right)-w\left(\cdot,t_{s}\right)\right\Vert _{\infty}\\
 & \leq\sum_{i=1}^{s}\left\Vert \tilde{\mathcal{S}}w\left(\cdot,t_{i-1}\right)-w\left(\cdot,t_{i}\right)\right\Vert _{\infty}\leq C\left(t_{s}\right)\Delta t.
\end{align*}
This completes the proof of the theorem. 
\end{proof}

\

\begin{proof}
[Proof of Theorem \ref{thm:converg_mean}]Under the structure assumptions
of the bilinear term \eqref{eq:assump2}, the total statistical energy
$E=\bar{u}^{2}+\frac{1}{K}\mathrm{tr}R$ from the solutions of the
mean and covariance equations in \eqref{coupled_model} satisfies
\[
\frac{\mathrm{d}E}{\mathrm{d}t}=-2dE+\bar{u}\cdot F.
\]
Similarly, the RBM model \eqref{rbm_model} also satisfies the corresponding
energy equation with $\tilde{E}=\tilde{\bar{u}}^{2}+\frac{1}{K}\mathrm{tr}\tilde{R}$
with the consistent structure symmetry
\[
\frac{\mathrm{d}\tilde{E}}{\mathrm{d}t}=-2d\tilde{E}+\tilde{\bar{u}}\cdot F.
\]
By taking the difference between the above two equations, we can write
the solution for the difference between the statistical energy in
truth and RBM approximation formally as
\[
\delta E\left(t\right)=\delta\bar{u}^{2}+\frac{1}{K}\delta\mathrm{tr}R=\int_{0}^{t}e^{-2d(t-s)}\delta\bar{u}\left(s\right)\cdot F\left(s\right)\mathrm{d}s,
\]
where we assume the initial states are the same, $E\left(0\right)=\tilde{E}\left(0\right)$,
and $\delta E=E-\tilde{E}$, $\delta\bar{u}=\bar{u}-\tilde{\bar{u}}$.
Then using the statistical estimation for the total variance derived
in \eqref{eq:converg_coeff}, we find
\begin{align}
M\left|\delta\bar{u}\right|\leq\left|\bar{u}^{2}-\tilde{\bar{u}}^{2}\right| & \leq\frac{1}{K}\left|\mathrm{tr}R-\mathrm{tr}\tilde{R}\right|+\int_{0}^{t}e^{-2d(t-s)}\left|\delta\bar{u}\right|\cdot\left|F\right|\mathrm{d}s\nonumber \\
 & \leq C\left(T\right)\Delta t+\int_{0}^{t}\left|\delta\bar{u}\right|\cdot\left|F\right|e^{-2d(t-s)}\mathrm{d}s.\label{eq:bound_mean}
\end{align}
Above on the left hand side, we use the assumption that the mean states
$\left|\bar{u}+\tilde{\bar{u}}\right|>M$ have a uniformly common
bound from below. This is observed in the numerical simulations. Finally,
applying Gr\"{o}nwall's inequality to \eqref{eq:bound_mean}, we
have the uniformly bound for the convergence of the statistical mean
state
\[
\left|\delta\bar{u}\right|\left(t\right)\leq C\left(T\right)\Delta te^{\int_{0}^{t}\left|F\right|e^{-2d(t-s)}\mathrm{d}s}\leq C_{1}\left(T\right)\Delta t.
\]
The final result \eqref{eq:converg_mean} is achieved by taking the
supremum for $t\in\left[0,T\right]$. 
\end{proof}

\renewcommand\theequation{B\arabic{equation}}
\setcounter{equation}{0}

\subsection{Explicit equations for the one-layer L-96 system\label{subsec:Explicit-RBM-1layer}}

According to the general equations \eqref{coupled_model}, the \emph{full stochastic-statistical formulation} for the L-96 system \eqref{eq:l96_homo}
can be derived as
\begin{equation}
\begin{aligned}\frac{\mathrm{d}\bar{u}}{\mathrm{d}t}= & \frac{1}{J^{2}}\sum_{k}\gamma_{k}r_{k}-\bar{u}+F,\\
\frac{\mathrm{d}r_{k}}{\mathrm{d}t}= & \frac{1}{J}\sum_{m-n=k}\left[\left\langle Z_{m}Z_{n}^{*}Z_{k}^{*}\right\rangle _{p}\gamma_{mn}+\left\langle Z_{m}^{*}Z_{n}Z_{k}\right\rangle _{p}\gamma_{mn}^{*}\right]\\
 & -2\left(\mathrm{Re}\gamma_{k}\bar{u}+1\right)r_{k}+\epsilon^{-1}\left(\left\langle \left|Z_{k}\right|^{2}\right\rangle _{p}-r_{k}\right),\\
\frac{\mathrm{d}Z_{k}}{\mathrm{d}t}= & \frac{1}{J}\sum_{m-n=k}\left(Z_{m}Z_{n}^{*}-r_{m}\delta_{mn}\right)\gamma_{mn}-\left(\gamma_{k}^{*}\bar{u}+1\right)Z_{k}.
\end{aligned}
\label{eq:full_l96}
\end{equation}
Above we have the coupling coefficients $\gamma_{mn}=\exp\left(i2\pi\frac{m+n}{J}\right)-\exp\left(-i2\pi\frac{n-2m}{J}\right)$
and $\gamma_{k}\coloneqq\gamma_{kk}$. The first two equations are
deterministic providing the statistical mean and variance dynamics.
The third equation characterizes the stochastic evolution of the coefficients
$Z_{k}$ (with randomness from the initial ensemble). The higher-order
moments in the variance equations are recovered by the stochastic
equation containing non-Gaussian information, $\left\langle f\right\rangle _{p}=\int f\mathrm{d}p$.
Therefore, the above equations \eqref{eq:full_l96} provide a closed
system by coupling the statistical and stochastic equations including
all the higher-order feedbacks. Especially, additional relaxation
term is added to the variance equation for $r_{k}$ to guarantee the
consistent statistical mean and variance in the model prediction.
In one time update of the above full equations, the mean equation requires $J$ operations for the summation of all the variance; Each component of the $J$ variances requires $J$ operations in the summation together with $M$ operations needed in computing the empirical average  of the third moments; Finally, the ensemble simulation of the stochastic coefficients requires $J^2$ operations for each of the $M$ samples.

Correspondingly, we can derive the \emph{RBM model} for
the one-layer L-96 system \eqref{eq:full_l96} in the time updating interval $t\in\left(t_{s},t_{s+1}\right]$ as
\begin{equation}
\begin{aligned}\frac{\mathrm{d}\bar{u}}{\mathrm{d}t}= & \frac{1}{J^{2}}\sum_{k}\gamma_{k}r_{k}-\bar{u}+F,\\
\frac{\mathrm{d}r_{k}}{\mathrm{d}t}= & c_{p}\sum_{m\in\mathcal{I}_{k}^{s}}\left\langle Z_{m}Z_{m-k}^{*}Z_{k}^{*}\right\rangle _{p_{M}}\gamma_{m,m-k}+\left\langle Z_{m}^{*}Z_{m-k}Z_{k}\right\rangle _{p_{M}}\gamma_{m,m-k}^{*}\\
 & -2\mathrm{Re}\gamma_{k}\bar{u}r_{k}-2r_{k}\;+\epsilon^{-1}\left(\left\langle \left|Z_{k}\right|^{2}\right\rangle _{p_{M}}-r_{k}\right),\\
\frac{\mathrm{d}Z_{k}^{\left(i\right)}}{\mathrm{d}t}= & c_{p}\sum_{m\in\mathcal{I}_{k}^{s}}\left(Z_{m}^{\left(i\right)}Z_{m-k}^{\left(i\right)*}-r_{m}\delta_{mn}\right)\gamma_{m,m-k}-\left(\gamma_{k}^{*}\bar{u}+1\right)Z_{k}^{\left(i\right)}.
\end{aligned}
\label{eq:rbm_l96}
\end{equation}
Above, the deterministic solutions of the statistical mean $\bar{u}$
and variances $r_{k}$ are closed by the stochastic equation
for the coefficients $Z_{k}$, which are sampled by an ensemble simulation
using  a  small sample size of trajectories $i=1,\cdots,M_{1}$. The empirical statistics in the
dynamics are computed from the ensemble average of the sample realizations
at each time updating step
\[
\left\langle f\right\rangle _{p_{M}}=\frac{1}{M_{1}}\sum_{i=1}^{M_{1}}f\left(Z^{\left(i\right)}\right).
\] 
Importantly, it is crucial
to use the consistent scaling factor $c_{p}=\frac{1}{J}\frac{J-1}{p-1}$
according to the batch size $p$ in the summation for higher-order feedback.
Through the RBM approach, the set of total spectral modes $\left\{ Z_{k}\right\} _{\left|k\right|\leq J/2}$
are randomly divided into small batches $\mathcal{I}_{k}^{s}$ with
size $p$ at each time step $t=t_{s}$. This enables the large computational reduction from $O(J^2)$ to $O(Jp)$ for each single sample trajectory and reduces the sample size from a very large $M$ to $M_1\ll M$ sampling only the $p$-dimensional batch subspace.

\subsection{Explicit equations for the two-layer L-96 system\label{subsec:RBM-equations-2layer}}

Following the same procedure using the general formulation \eqref{coupled_model},
we first derive the \emph{mean equations for the large and small scale
states} of the two-layer L-96 system \eqref{eq:l96_2layer}
\begin{equation}
\begin{aligned}\frac{\mathrm{d}\bar{u}}{\mathrm{d}t} & =\frac{1}{J^{2}}\sum_{\left|k\right|\leq J/2}\gamma_{k}^{u}r_{k}^{u}-\bar{u}+F-\frac{hc}{b}L\bar{v},\\
\frac{\mathrm{d}\bar{v}}{\mathrm{d}t} & =\frac{cb}{J^{2}L^{2}}\sum_{\left|l\right|\leq JL/2}\gamma_{l}^{v}r_{l}^{v}-c\bar{v}+\frac{hc}{b}\bar{u}.
\end{aligned}
\label{eq:mean_2layer}
\end{equation}
Next, we have the set of \emph{covariance equations} associated with
the large and small scale stochastic coefficients from the Fourier
decomposition in \eqref{eq:modes_2layer}
\begin{equation}
\begin{aligned}\frac{\mathrm{d}r_{k}^{u}}{\mathrm{d}t}= & \frac{1}{J}\sum_{m-n=k}\left[\left\langle Z_{m}Z_{n}^{*}Z_{k}^{*}\right\rangle \gamma_{mn}^{u}+\left\langle Z_{m}^{*}Z_{n}Z_{k}\right\rangle \gamma_{mn}^{u*}\right]-2\left(1+\mathrm{Re}\gamma_{k}^{u}\bar{u}\right)r_{k}^{u}\\
 & -\frac{hc}{b}\frac{1}{L}\sum_{s=-L/2+1}^{L/2}\lambda_{k+sJ}^{*}r_{k+sJ}^{x}+\epsilon^{-1}\left(\left\langle \left|Z_{k}\right|^{2}\right\rangle -r_{k}^{u}\right),\quad\left|k\right|\leq J/2,\\
\frac{\mathrm{d}r_{l}^{v}}{\mathrm{d}t}= & \frac{cb}{JL}\sum_{p-q=l}\left[\left\langle Y_{p}Y_{q}^{*}Y_{l}^{*}\right\rangle \gamma_{pq}^{v}+\left\langle Y_{p}^{*}Y_{q}Y_{l}\right\rangle \gamma_{pq}^{v*}\right]-2c\left(1+b\mathrm{Re}\gamma_{l}^{v}\bar{v}\right)r_{l}^{v}\\
 & +\frac{hc}{b}\lambda_{l}r_{l}^{x}+\epsilon^{-1}\left(\left\langle \left|Y_{l}\right|^{2}\right\rangle -r_{l}^{v}\right),\\
\frac{\mathrm{d}r_{l}^{x}}{\mathrm{d}t}= & \frac{1}{J}\sum_{m-n=k}\gamma_{mn}^{u}\left\langle Z_{m}Z_{n}^{*}Y_{l}^{*}\right\rangle +\frac{cb}{JL}\sum_{p-q=l}\gamma_{pq}^{v*}\left\langle Y_{p}^{*}Y_{q}Z_{\mathrm{mod}\left(l,J\right)}\right\rangle -\left(\gamma_{\mathrm{mod}\left(l,J\right)}^{u*}\bar{u}+cb\gamma_{l}^{v}\bar{v}+1+c\right)r_{l}^{x}\\
 & +\frac{hc}{b}\lambda_{l}^{*}\left(r_{\mathrm{mod}\left(l,J\right)}^{u}-\frac{1}{L}r_{l}^{v}\right)+\epsilon^{-1}\left(\left\langle Z_{\mathrm{mod}\left(l,J\right)}Y_{l}^{*}\right\rangle -r_{l}^{x}\right),\quad\left|l\right|\leq JL/2,
\end{aligned}
\label{eq:cov_2layer}
\end{equation}
where $r_{k}^{u}=\left\langle \left|Z_{k}\right|^{2}\right\rangle ,r_{l}^{v}=\left\langle \left|Y_{l}\right|^{2}\right\rangle $
are the variances for the large and small scales, and $r_{l}^{x}=\left\langle Z_{\mathrm{mod}\left(l,J\right)}Y_{l}^{*}\right\rangle $
gives the cross-covariance. The coupling coefficients are $\gamma_{mn}^{u}=e^{2\pi i\frac{m+n}{J}}-e^{-2\pi i\frac{2m-n}{J}}$,
$\gamma_{pq}^{v}=e^{-2\pi i\frac{p+q}{JL}}-e^{2\pi i\frac{2p-q}{JL}}$,
$\lambda_{l}=\frac{1-e^{-2\pi i\frac{l}{J}}}{1-e^{-2\pi i\frac{l}{JL}}}$,
and $\gamma_{k}^{u}=e^{4\pi i\frac{k}{J}}-e^{-2\pi i\frac{k}{J}}$,
$\gamma_{l}^{v}=e^{-4\pi i\frac{l}{JL}}-e^{2\pi i\frac{l}{JL}}$.

Correspondingly, the RBM approximation for the covariance equations
are derived according to the full and reduced-order stochastic equations
proposed in \eqref{eq:fluc_2layer-rbm1} and \eqref{eq:fluc_2layer_rbm2}.
The explicit \emph{RBM model equations} during the time interval $t\in\left(t_{s},t_{s+1}\right]$
can be found as
\begin{equation}
\begin{aligned}\frac{\mathrm{d}r_{k}^{u}}{\mathrm{d}t}= & c_{p}\sum_{m\in\mathcal{I}_{k}^{s}}\left[\left\langle Z_{m}Z_{m-k}^{*}Z_{k}^{*}\right\rangle \gamma_{m,m-k}^{u}+\left\langle Z_{m}^{*}Z_{m-k}Z_{k}\right\rangle \gamma_{m,m-k}^{u*}\right]-2\left(1+\mathrm{Re}\gamma_{k}^{u}\bar{u}\right)r_{k}^{u}\\
 & -c_{L}\frac{hc}{b}\sum_{k+sJ\in\mathcal{I}_{k}^{s}}\lambda_{k+sJ}^{*}r_{k+sJ}^{x}+\epsilon^{-1}\left(\left\langle \left|Z_{k}\right|^{2}\right\rangle -r_{k}^{u}\right),\\
\frac{\mathrm{d}r_{l}^{v}}{\mathrm{d}t}= & c_{q}cb\sum_{p\in\mathcal{J}_{l}^{s}}\left[\left\langle Y_{p}Y_{p-l}^{*}Y_{l}^{*}\right\rangle \gamma_{p,p-l}^{v}+\left\langle Y_{p}^{*}Y_{p-l}Y_{l}\right\rangle \gamma_{p,p-l}^{v*}\right]-2c\left(1+b\mathrm{Re}\gamma_{l}^{v}\bar{v}\right)r_{l}^{v}\\
 & +\frac{hc}{b}\lambda_{l}r_{l}^{x}+\epsilon^{-1}\left(\left\langle \left|Y_{l}\right|^{2}\right\rangle -r_{l}^{v}\right),\\
\frac{\mathrm{d}r_{l}^{x}}{\mathrm{d}t}= & c_{p}\sum_{m\in\mathcal{I}_{k}^{s}}\gamma_{m,m-k}^{u}\left\langle Z_{m}Z_{m-k}^{*}Y_{l}^{*}\right\rangle +c_{q}cb\sum_{p\in\mathcal{J}_{l}^{s}}\gamma_{p,p-l}^{v*}\left\langle Y_{p}^{*}Y_{p-l}Z_{\mathrm{mod}\left(l,J\right)}\right\rangle -\left(\gamma_{\mathrm{mod}\left(l,J\right)}^{u*}\bar{u}+cb\gamma_{l}^{v}\bar{v}+1+c\right)r_{l}^{x}\\
 & +\frac{hc}{b}\lambda_{l}^{*}\left(r_{\mathrm{mod}\left(l,J\right)}^{u}-\frac{1}{L}r_{l}^{v}\right)+\epsilon^{-1}\left(\left\langle Z_{\mathrm{mod}\left(l,J\right)}Y_{l}^{*}\right\rangle -r_{l}^{x}\right),
\end{aligned}
\label{eq:cov_2layer_rbm}
\end{equation}
with the important rescaling parameters $c_{p}=\frac{1}{J}\frac{J-1}{p-1}$, $c_{q}=\frac{1}{JL}\frac{JL-1}{q-1}$,
and $c_{L}=\frac{1}{L}\frac{L-1}{q-1}$. In the same way as the one-layer case, the higher-order moments are
computed through the empirical average of the computed samples.

\bibliographystyle{plain}
\bibliography{refs}

\begin{thebibliography}{10}

\bibitem{arnold2013stochastic}
HM~Arnold, IM~Moroz, and TN~Palmer.
\newblock Stochastic parametrizations and model uncertainty in the {L}orenz '96
  system.
\newblock {\em Philosophical Transactions of the Royal Society A: Mathematical,
  Physical and Engineering Sciences}, 371(1991):20110479, 2013.

\bibitem{bishop2001adaptive}
Craig~H Bishop, Brian~J Etherton, and Sharanya~J Majumdar.
\newblock Adaptive sampling with the ensemble transform {K}alman filter. part
  i: Theoretical aspects.
\newblock {\em Monthly weather review}, 129(3):420--436, 2001.

\bibitem{calvello2022ensemble}
Edoardo Calvello, Sebastian Reich, and Andrew~M Stuart.
\newblock Ensemble {K}alman methods: a mean field perspective.
\newblock {\em arXiv preprint arXiv:2209.11371}, 2022.

\bibitem{chen2020predicting}
Nan Chen and Andrew~J Majda.
\newblock Predicting observed and hidden extreme events in complex nonlinear
  dynamical systems with partial observations and short training time series.
\newblock {\em Chaos: An Interdisciplinary Journal of Nonlinear Science},
  30(3):033101, 2020.

\bibitem{cousins2015unsteady}
Will Cousins and Themistoklis~P Sapsis.
\newblock Unsteady evolution of localized unidirectional deep-water wave
  groups.
\newblock {\em Physical Review E}, 91(6):063204, 2015.

\bibitem{daum2003curse}
Fred Daum and Jim Huang.
\newblock Curse of dimensionality and particle filters.
\newblock In {\em 2003 IEEE Aerospace Conference Proceedings (Cat. No.
  03TH8652)}, volume~4, pages 4\_1979--4\_1993. IEEE, 2003.

\bibitem{dewar2007zonal}
Robert~L Dewar and Raden~Farzand Abdullatif.
\newblock Zonal flow generation by modulational instability.
\newblock In {\em Frontiers in Turbulence and Coherent Structures}, pages
  415--430. World Scientific, 2007.

\bibitem{Diamond2010}
Patrick~H. Diamond, Sanae-I. Itoh, and Kimitaka Itoh.
\newblock {\em {Modern Plasma Physics}}.
\newblock Cambridge University Press, Cambridge, 2010.

\bibitem{donoho2000high}
David~L Donoho et~al.
\newblock High-dimensional data analysis: The curses and blessings of
  dimensionality.
\newblock {\em AMS math challenges lecture}, 1(2000):32, 2000.

\bibitem{friedman1997bias}
Jerome~H Friedman.
\newblock On bias, variance, 0/1---loss, and the curse-of-dimensionality.
\newblock {\em Data mining and knowledge discovery}, 1(1):55--77, 1997.

\bibitem{Frisch1995}
Uriel Frisch.
\newblock {\em {Turbulence}}.
\newblock Cambridge University Press, nov 1995.

\bibitem{gao2023transition}
Yuan Gao, Tiejun Li, Xiaoguang Li, and Jian-Guo Liu.
\newblock Transition path theory for {L}angevin dynamics on manifolds: Optimal
  control and data-driven solver.
\newblock {\em Multiscale Modeling \& Simulation}, 21(1):1--33, 2023.

\bibitem{gao2022master}
Yuan Gao, Wuchen Li, and Jian-Guo Liu.
\newblock Master equations for finite state mean field games with nonlinear
  activations.
\newblock {\em arXiv preprint arXiv:2212.05675}, 2022.

\bibitem{gao2023data}
Yuan Gao, Jian-Guo Liu, and Nan Wu.
\newblock Data-driven efficient solvers for {L}angevin dynamics on manifold in
  high dimensions.
\newblock {\em Applied and Computational Harmonic Analysis}, 62:261--309, 2023.

\bibitem{giggins2019stochastically}
Brent Giggins and Georg~A Gottwald.
\newblock Stochastically perturbed bred vectors in multi-scale systems.
\newblock {\em Quarterly Journal of the Royal Meteorological Society},
  145(719):642--658, 2019.

\bibitem{houtekamer1996system}
Peter~L Houtekamer, Louis Lefaivre, Jacques Derome, Harold Ritchie, and
  Herschel~L Mitchell.
\newblock A system simulation approach to ensemble prediction.
\newblock {\em Monthly Weather Review}, 124(6):1225--1242, 1996.

\bibitem{jin2020random}
Shi Jin, Lei Li, and Jian-Guo Liu.
\newblock Random batch methods ({RBM}) for interacting particle systems.
\newblock {\em Journal of Computational Physics}, 400:108877, 2020.

\bibitem{jin2021convergence}
Shi Jin, Lei Li, and Jian-Guo Liu.
\newblock Convergence of the random batch method for interacting particles with
  disparate species and weights.
\newblock {\em SIAM Journal on Numerical Analysis}, 59(2):746--768, 2021.

\bibitem{kalnay2003atmospheric}
Eugenia Kalnay.
\newblock {\em Atmospheric modeling, data assimilation and predictability}.
\newblock Cambridge university press, 2003.

\bibitem{lesieur1987turbulence}
Marcel Lesieur.
\newblock {\em Turbulence in fluids: stochastic and numerical modelling},
  volume 488.
\newblock Nijhoff Boston, MA, 1987.

\bibitem{leutbecher2008ensemble}
Martin Leutbecher and Tim~N Palmer.
\newblock Ensemble forecasting.
\newblock {\em Journal of computational physics}, 227(7):3515--3539, 2008.

\bibitem{lorenz1996predictability}
Edward~N Lorenz.
\newblock Predictability: A problem partly solved.
\newblock In {\em Proc. Seminar on predictability}, volume~1. Reading, 1996.

\bibitem{luo2022stability}
Yushuang Luo, Xiantao Li, and Wenrui Hao.
\newblock Stability preserving data-driven models with latent dynamics.
\newblock {\em Chaos: An Interdisciplinary Journal of Nonlinear Science},
  32(8):081103, 2022.

\bibitem{majda2006nonlinear}
Andrew Majda and Xiaoming Wang.
\newblock {\em Nonlinear dynamics and statistical theories for basic
  geophysical flows}.
\newblock Cambridge University Press, 2006.

\bibitem{majda2016introduction}
Andrew~J Majda.
\newblock {\em Introduction to turbulent dynamical systems in complex systems}.
\newblock Springer, 2016.

\bibitem{majda2019statistical}
Andrew~J Majda, MNJ Moore, and Di~Qi.
\newblock Statistical dynamical model to predict extreme events and anomalous
  features in shallow water waves with abrupt depth change.
\newblock {\em Proceedings of the National Academy of Sciences},
  116(10):3982--3987, 2019.

\bibitem{majda2016improving}
Andrew~J Majda and Di~Qi.
\newblock Improving prediction skill of imperfect turbulent models through
  statistical response and information theory.
\newblock {\em Journal of Nonlinear Science}, 26:233--285, 2016.

\bibitem{majda2018strategies}
Andrew~J Majda and Di~Qi.
\newblock Strategies for reduced-order models for predicting the statistical
  responses and uncertainty quantification in complex turbulent dynamical
  systems.
\newblock {\em SIAM Review}, 60(3):491--549, 2018.

\bibitem{majda2019linear}
Andrew~J Majda and Di~Qi.
\newblock Linear and nonlinear statistical response theories with prototype
  applications to sensitivity analysis and statistical control of complex
  turbulent dynamical systems.
\newblock {\em Chaos: An Interdisciplinary Journal of Nonlinear Science},
  29(10):103131, 2019.

\bibitem{maulik2019subgrid}
Romit Maulik, Omer San, Adil Rasheed, and Prakash Vedula.
\newblock Subgrid modelling for two-dimensional turbulence using neural
  networks.
\newblock {\em Journal of Fluid Mechanics}, 858:122--144, 2019.

\bibitem{nicholson1983introduction}
Dwight~Roy Nicholson and Dwight~R Nicholson.
\newblock {\em Introduction to plasma theory}, volume~1.
\newblock Wiley New York, 1983.

\bibitem{orrell2003model}
David Orrell.
\newblock Model error and predictability over different timescales in the
  {L}orenz '96 systems.
\newblock {\em Journal of the atmospheric sciences}, 60(17):2219--2228, 2003.

\bibitem{palmer2019ecmwf}
Tim Palmer.
\newblock The {ECMWF} ensemble prediction system: Looking back (more than) 25
  years and projecting forward 25 years.
\newblock {\em Quarterly Journal of the Royal Meteorological Society},
  145:12--24, 2019.

\bibitem{pedlosky2013geophysical}
Joseph Pedlosky.
\newblock {\em Geophysical fluid dynamics}.
\newblock Springer Science \& Business Media, 2013.

\bibitem{qi2022machine}
Di~Qi and John Harlim.
\newblock Machine learning-based statistical closure models for turbulent
  dynamical systems.
\newblock {\em Philosophical Transactions of the Royal Society A},
  380(2229):20210205, 2022.

\bibitem{qi2023data}
Di~Qi and John Harlim.
\newblock A data-driven statistical-stochastic surrogate modeling strategy for
  complex nonlinear non-stationary dynamics.
\newblock {\em Journal of Computational Physics}, 485:112085, 2023.

\bibitem{qi2023random}
Di~Qi and Jian-Guo Liu.
\newblock A random batch method for efficient ensemble forecasts of multiscale
  turbulent systems.
\newblock {\em Chaos: An Interdisciplinary Journal of Nonlinear Science},
  33(2):023113, 2023.

\bibitem{qi2022anomalous}
Di~Qi and Eric Vanden-Eijnden.
\newblock Anomalous statistics and large deviations of turbulent water waves
  past a step.
\newblock {\em AIP Advances}, 12(2):025016, 2022.

\bibitem{reich2015probabilistic}
Sebastian Reich and Colin Cotter.
\newblock {\em Probabilistic forecasting and {B}ayesian data assimilation}.
\newblock Cambridge University Press, 2015.

\bibitem{sapsis2013statistically}
Themistoklis~P Sapsis and Andrew~J Majda.
\newblock Statistically accurate low-order models for uncertainty
  quantification in turbulent dynamical systems.
\newblock {\em Proceedings of the National Academy of Sciences},
  110(34):13705--13710, 2013.

\bibitem{smith2001disentangling}
Leonard~A Smith.
\newblock Disentangling uncertainty and error: On the predictability of
  nonlinear systems.
\newblock {\em Nonlinear dynamics and statistics}, pages 31--64, 2001.

\bibitem{stuart2010inverse}
Andrew~M Stuart.
\newblock Inverse problems: a {B}ayesian perspective.
\newblock {\em Acta numerica}, 19:451--559, 2010.

\bibitem{tong2021extreme}
Shanyin Tong, Eric Vanden-Eijnden, and Georg Stadler.
\newblock Extreme event probability estimation using {PDE}-constrained
  optimization and large deviation theory, with application to tsunamis.
\newblock {\em Communications in Applied Mathematics and Computational
  Science}, 16(2):181--225, 2021.

\bibitem{vanden2001non}
Eric Vanden~Eijnden.
\newblock Non-{G}aussian invariant measures for the {M}ajda model of decaying
  turbulent transport.
\newblock {\em Communications on Pure and Applied Mathematics: A Journal Issued
  by the Courant Institute of Mathematical Sciences}, 54(9):1146--1167, 2001.

\bibitem{varadhan2007stochastic}
SR~Srinivasa Varadhan.
\newblock {\em Stochastic processes}, volume~16.
\newblock American Mathematical Soc., 2007.

\bibitem{wilks2005effects}
Daniel~S Wilks.
\newblock Effects of stochastic parametrizations in the {L}orenz '96 system.
\newblock {\em Quarterly Journal of the Royal Meteorological Society: A journal
  of the atmospheric sciences, applied meteorology and physical oceanography},
  131(606):389--407, 2005.

\end{thebibliography}

\end{document}